\documentclass{amsart}
\title{$L^2$-index theorems, KK-theory, and connections}
\author{Thomas Schick}
\email{schick@uni-math.gwdg.de, http://www.uni-math.gwdg.de/schick}
\address{FB Mathematik, Universit{\"a}t G{\"o}ttingen\\ Bunsenstr.~3, 37073 G{\"o}ttingen, Germany}
\keywords{K-theory, index theory,curvature, Chern character,
  $C^*$-algebra, $L^2$-index, Mishchenko-Fomenko index theorem, trace}
        
\subjclass{19K35, 19K56, 46M20, 46L80, 58J22}

\usepackage{amsmath,amsthm,latexsym,amscd,amssymb}
\usepackage[backref,bookmarks=true,bookmarksnumbered=true,hypertex]{hyperref}
\usepackage[backrefs]{amsrefs}

\usepackage{amsthm}

\swapnumbers
\theoremstyle{plain}
\newtheorem{theorem}{Theorem}[section]
\newtheorem{lemma}[theorem]{Lemma}
\newtheorem{corollary}[theorem]{Corollary}
\newtheorem{proposition}[theorem]{Proposition}

\theoremstyle{definition}
\newtheorem{definition}[theorem]{Definition}
\newtheorem{example}[theorem]{Example}
\newtheorem{notation}[theorem]{Notation}

\theoremstyle{remark}
\newtheorem{remark}[theorem]{Remark}

\usepackage{amsmath,amsfonts,latexsym}

\newcommand{\reals}{\mathbb{R}}

\newcommand{\complexs}{\mathbb{C}}

\newcommand{\naturals}{\mathbb{N}}

\newcommand{\integers}{\mathbb{Z}}

\newcommand{\rationals}{\mathbb{Q}}

\DeclareMathOperator{\id}{id}

\newcommand{\boundedops}{\mathcal{B}}

\newcommand{\abs}[1]{\left\lvert#1\right\rvert} %absolute value
\newcommand{\norm}[1]{\left\lVert#1\right\rVert}

\newcommand{\tensor}{\otimes}
\newcommand{\into}{\hookrightarrow}
\newcommand{\onto}{\twoheadrightarrow}
\newcommand{\iso}{\cong}

\DeclareMathOperator{\im}{im}      %image
\DeclareMathOperator{\End}{End}    %Endomorphisms
\DeclareMathOperator{\Hom}{Hom}    %Homomorphisms
\DeclareMathOperator{\Tor}{Tor}    %Tor
\DeclareMathOperator{\tr}{tr}
\DeclareMathOperator{\pr}{pr}
\DeclareMathOperator{\coker}{coker}

\DeclareMathOperator{\ind}{ind}

\DeclareMathOperator{\Td}{Td}
\DeclareMathOperator{\ch}{ch}

\newcommand{\forget}[1]{}

\newcommand{\innerprod}[1]{\langle #1 \rangle}

{\catcode`@=11\global\let\c@equation=\c@theorem}

\allowdisplaybreaks[2]

\newcommand{\NeumannN}{\mathcal{N}}

\begin{document}
%\date{Last compiled \today; last edited  \heuteIst or later}

%\title{$L^2$-index theorems, KK-theory, and connections}
%\author{Thomas Schick\thanks{
%e-mail: {schick@uni-math.gwdg.de}\protect\\
%%~\protect\href{mailto:thomas.schick@math.uni-muenster.de}{thomas.schick@math.uni-muenster.de}
%%www:~\protect\href{http://www.math.uni-muenster.de/u/schickt}
%%{http://www.math.uni-muenster.de/u/schickt} 
%FB Mathematik, Universit\"{\"a}t G{\"o}ttingen, Bunsenstr.~3, 37073 G{\"o}ttingen, Germany}\\
%}
        
%\maketitle

\begin{abstract} \parindent0pt
Let $M$ be a compact manifold and $D$ a Dirac type differential
operator on $M$. Let $A$ be a $C^*$-algebra. Given a bundle $W$ (with connection) of
$A$-modules over $M$,
the operator $D$ can be twisted with this bundle. One can then use a
trace on $A$ to define numerical indices of this twisted operator. We
prove an explicit formula for these indices. Our result does complement
the Mishchenko-Fomenko index theorem valid in the same situation.
We establish generalizations of these explicit index formulas if the trace is only
defined on a dense and holomorphically closed subalgebra $\mathcal{B}$.

As a corollary, we prove a generalized Atiyah $L^2$-index theorem if
the twisting bundle is flat.

\smallskip
There are actually many different ways to define these numerical
indices. From their construction, it is not clear at all that they
coincide. An
substantial part of the paper is a complete proof of their equality. In
particular, we establish the (well known but not well documented)
equality of Atiyah's definition of the $L^2$-index with a K-theoretic definition.

In case $A$ is a von Neumann algebra of type 2, we put special
emphasis on the calculation and interpretation of the center valued
index. This completely contains all the K-theoretic information about
the index of the twisted operator.

\smallskip
Some of our calculations are done in the framework of bivariant
KK-theory.

%MSC 2000: 19K35, 19K56, 46M20, 46L80, 58J22
\end{abstract}

\maketitle
\tableofcontents

\section{Introduction}

Let $M$ be a closed smooth manifold, $D\colon \Gamma(E)\to \Gamma(E)$
a generalized Dirac operator on the finite dimensional (graded) Dirac
bundle $E$ over $M$.

Assume that $A$ is a $C^*$-algebra and $W$ a smooth bundle of finitely
generated projective modules over $A$ equipped with a connection
$\nabla_W$. In this situation, one can define the twisted Dirac
operator $D_W$ (compare \eqref{eq:def_of_twisted_dirac}). The
resulting operator is an elliptic $A$-operator in the sense of
Mishchenko-Fomenko \cite{MIshchenko-Fomenko}. In particular, its index
$\ind(D_W)\in K_0(A)$
as an element of the K-theory of $A$ is defined. Mishchenko and
Fomenko prove a
formula for this index (or rather its rationalization). Improvements
to a K-theoretic level, and even equivariant generalizations, can be
found e.g.~in \cite{MR1787114,MR1865952,MR1266008}. However, for
certain purposes these formulas are rather inexplicit.

The main goal of this paper is an explicit index formula in this context, in
terms of the curvature of the twisting bundle $W$. This can not be
done directly for the index. However, whenever we have a trace
$\tau\colon A\to Z$ with values in any commutative $C^*$-algebra
(e.g.~the complex numbers), it induces a homomorphism $\tau\colon
K_0(A)\to Z$, and we get the explicit formula 
\begin{equation}\label{eq:explict_index_formula}
\tau(\ind(D_W)) = \innerprod{ \ch(\sigma(D)) \cup \Td(T_\complexs
      M) \cup
      \ch_\tau(W) , [TM]} \in Z
\end{equation}
in Theorem \ref{theo:index_theorem}, where the crucial term
$\ch_\tau(W)$ can be calculated directly from the curvature of $W$.

Such a Chern-Weil approach to higher index theorems can replace
heat equation proofs. This is e.g.~remarked by Mathai in
\cite[p.~14]{MR2004e:19002} and can be used to simplify
his proof of the Novikov conjecture for low dimensional cohomology
classes. Applications of the explicit index formulas to the study of
manifolds with metrics of positive scalar curvature are obtained by
Hanke and Schick in
\cite{math.GT/0403257}. There, we also use
generalizations of Equation \eqref{eq:explict_index_formula} to situations where a trace is
only defined on a dense subalgebra $\mathcal{B}$ of $A$ which is
closed under holomorphic functional calculus, as explained in Section
\ref{sec:trace-class-subalg}.

\begin{corollary}
Assume that $M$ is connected. Then the index formula
\eqref{eq:explict_index_formula} implies that for a 
  \emph{flat} bundle $W$
  \begin{equation*}
    \tau(\ind(D_W)) = \ind(D) \cdot d,
  \end{equation*}
  where $d:=\dim_\tau(W_x)$ is the ``fiber dimension'', the trace of
  the projection onto the (finitely generated projective) fiber $W_x$
  of $W$ over an arbitrary point $x\in M$.
\end{corollary}
If $A$ is a finite von Neumann
algebra and $\tau$ is its center valued trace, $\tau(\ind(D_W))$
contains as much information as $\ind(D_W)$.

There are several other ways to define an index for $D$ twisted with
$W$. The most direct is probably given by the Kasparov product of a
KK-element defined by $W$ with the index element $D$ defines in $KK(
C(M),\complexs)$. In Theorem \ref{theo:twist_and_Kasparov product} we
show that this coincides with the index defined directly using
the Mishchenko-Fomenko calculus.

If $A$ is a finite von Neumann algebra, it is more popular to twist with
$A$-Hilbert space bundles, where the fibers are ordinary Hilbert
spaces, but with an appropriate action of the von Neumann algebra
$A$. We show that we can assign such a bundle $l^2(W)$ to $W$ as above
(and
vice versa), and that the twisted indices obtained both ways are
essentially equal (compare Theorem \ref{theo:general_Atiyah} and
Corollary \ref{corol:l2-equals ordinary index}). Here we need the
additional assumption that the trace is normal.

A special situation occurs if $A=C^*\Gamma$ is the $C^*$-algebra of a
discrete group and $W$ is the flat bundle associated to a
unitary representation $\pi_1(M)\to C^*\Gamma$ induced from a
group homomorphism $\pi_1(M)\to \Gamma$. Associated to this
homomorphism is a $\Gamma$-covering space $\tilde M\to M$, and we can
lift $D$ to an operator $\tilde D$ on $\tilde M$. Atiyah defines the $L^2$-index
$\ind_{(2)}(\tilde D)$
of $\tilde D$ in terms of sections on $\tilde M$ and proves his
$L^2$-index theorem in \cite{Atiyah(1976)}. We show in Theorem
\ref{theo:covering_equal_twist} that there is a direct correspondence
between this $L^2$-index (and generalizations hereof) and the index of
$D$ twisted with the flat $C^*\Gamma$-module bundle $W$ as
above. In particular, 
\begin{equation*}
  \ind_{(2)}(\tilde D) = t(\ind(D_W)),
\end{equation*}
where $W=\tilde M\times_{\Gamma} C^*\Gamma$ is the flat bundle with
fiber $C^*\Gamma$ associated to the $\Gamma$-covering $\tilde M$ and
$t\colon C^*\Gamma\to\NeumannN\Gamma\to\complexs$ is the canonical
trace (producing the coefficient of the trivial element and factoring
through the group von Neumann algebra $\NeumannN\Gamma$).

Finally, we consider the situation where $\mathcal{B}$ is a dense
subalgebra of the $C^*$-algebra $A$ with a trace
$\tau\colon\mathcal{B}\to Z$. The prototypical situation is the
algebra of trace class operators in the algebra of compact operators
on a separable Hilbert space, with the ordinary trace. If $\mathcal{B}$ is closed under
holomorphic functional calculus, then $\tau$ induces a homomorphism
$\tau\colon K_0(A)\to Z$. In this situation, if the Hilbert $A$-module
bundle $W$ is induced up from a bundle $\mathcal{W}$ of finitely
generated projective $\mathcal{B}$-modules, then we can define and use the
curvature of $\mathcal{W}$ to give an explicit expression for
$\ch_\tau(W)=\ch_\tau(\mathcal{W})$. We then prove the index formula
\begin{equation*}
\tau(\ind(D_W)) = \innerprod{ \ch(\sigma(D)) \cup \Td(T_\complexs
      M) \cup
      \ch_\tau(\mathcal{W}) , [TM]} 
\end{equation*}
Our proof of the  index formula for Hilbert $A$-module bundles works by just using a number of
crucial properties of the K-theory of $A$ and of $C(M)\tensor A$. Since
$\mathcal{B}$ is closed under holomorphic functional calculus, its
K-theory shares these properties. We will therefore only briefly
describe where changes in the first prove are necessary to obtain the
second result.

Along the way, we solve a number of related questions, in particular
the following.
\begin{enumerate}
\item We develop a Chern-Weil calculus for connections on Hilbert
  $A$-module bundles.
\item We proof existence and uniqueness of smooth structures on
  Hilbert $A$-module bundles, and show how K-theory of a manifold with coefficients
  in $A$ is described using smooth bundles.
\item The index $\ind(D_W)\in K_0(A)$ has to be defined in a
  complicated way, since kernel and cokernel of $D_W$ are in general
  not finitely generated projective over $A$. If $A$ is a von Neumann
  algebra, we prove that this caution is not necessary and that one
  can use the naive definition of the index.
\item We prove that for a finite von Neumann algebra, Hilbert
  $A$-modules and $A$-Hilbert spaces are equivalent categories, and
  that the same is true for bundles with corresponding fibers.
\item We establish a one-to-one correspondence between section of bundles on a
  $\Gamma$-covering space and of the associated flat $\NeumannN\Gamma$-bundle.
\end{enumerate}

\section{Notation and conventions}

Throughout this paper, $A$ denotes a unital $C^*$-algebra. Much of the
theory can be carried out for non-unital $C^*$-algebras, but for quite
a few statements, the existence of a unit is crucial, and they would
have to be reformulated considerable in the non-unital case. In our
applications, we are interested mainly in the reduced $C^*$-algebra
and the von Neumann algebra of a discrete group, which are unital.

For some of our constructions, we will have to restrict to the case
where $A$ is a von Neumann algebra.

\section{Hilbert modules and their properties}
\label{sec:hilb-modul-their}

In this section, we recall the notion of a Hilbert $C^*$-module and
its basic properties. A good and more comprehensive introduction to
this subject is e.g.~\cite{Lance(1995)} or \cite[Chapter 15]{Wegge-Olsen}. 

\begin{definition}

  A \emph{Hilbert $A$-module} $V$ is a right $A$-module $V$ with an
  $A$-valued ``inner product'' $\innerprod{\cdot,\cdot}\colon
  V\times V\to A$ 
  with the following properties:
  \begin{enumerate}
  \item $\innerprod{v_1,v_2a} =
    \innerprod{v_1,v_2}a\qquad\forall v_1,v_2\in V,\; a\in A$
  \item $\innerprod{v_1+v_2,v_3} =
    \innerprod{v_1,v_3}+\innerprod{v_2,v_3}\quad\forall v_1,v_2,v_3\in V,$
  \item $    \innerprod{v_1,v_2}=(\innerprod{v_2,v_1})^*\qquad\forall
    v_1,v_2\in V$,
  \item $\innerprod{v,v}$ is a non-negative self-adjoint element of
    the $C^*$-algebra 
    $A$ for each $v\in V$, and $\innerprod{v,v}=0$ if and only if
    $v=0$.
  \item The map $v\mapsto \abs{\innerprod{v,v}}_A^{1/2}$ is a norm on
    $V$, and $V$ is a Banach space with respect to this norm.
  \end{enumerate}
  Given two Hilbert $A$-modules $V$ and $W$, a \emph{Hilbert $A$-module
  morphism} $\Phi\colon V\to W$ is a continuous (right) $A$-linear map
  which has an adjoint $\Phi^*\colon W\to V$,
  i.e.~$\innerprod{\Phi(v),w}_W=\innerprod{v,\Phi^*(w)}_V$ for all
  $v\in V$, $w\in W$. The
  vector space of all such maps is denoted $\Hom_A(V,W)$.

  $\Hom_A(V,W)$ is an $\End_A(W)$-left-$\End_A(V)$-right
  module (but is not equipped with an inner product in general). The
  Hilbert $A$-module $V$ itself is an $\End_A(V)$-$A$-bimodule.
\end{definition}

\begin{example}
  The most important example of a Hilbert $A$-module is $A^n$ with
  inner product $\innerprod{(a_i),(b_i)} = \sum_{i=1}^n a_i^*b_i$.

  In this case, $\Hom_A(A^n,A^m) \iso M(m\times n, A)$, where the
  matrices act by multiplication from the left.The adjoint homomorphism is given by taking the
  transpose 
  matrix and the adjoint of each entry.
  In particular, $\End_A(A)\iso A$ as $C^*$-algebra. 

  We also consider $H_A$, the standard countably generated Hilbert
  $A$-module. It is the completion of $\bigoplus_{i\in\naturals} A$
  with respect to the norm $\abs{(a_i)}=\abs{\sum_{i\in\naturals}
    a_i^*a_i}_A^{1/2}$ and with the corresponding $A$-valued inner product.  

  Given two Hilbert $A$-modules $V$ and $W$, their direct sum $V\oplus
  W$ is a Hilbert $A$-module with $\innerprod{(v_1,w_1),(v_2,w_2)} =
  \innerprod{v_1,v_2}_V+\innerprod{w_1,w_2}_W$. 
\end{example}

In \cite[page 8]{Lance(1995)} the following result is proved.
\begin{lemma}\label{lem:Banach_space_Hom_A}
  Assume that $V$ and $W$ are Hilbert $A$-modules. Then $\Hom_A(V,W)$ is  a
  Banach space with the operator norm, and $\End_A:=\Hom_A(V,V)$ is a $C^*$-algebra.
\end{lemma}

If $A$ is a von Neumann algebra we get the following stronger
result:
\begin{proposition}\label{prop:End_A_is_von_Neumann_if_A_is}
  If $A$ is a von Neumann algebra then the same is true for $\End_A(H_A)$.
\end{proposition}
\begin{proof}
  This follows from the isomorphism $\End_A(H_A)\iso
  \boundedops(H)\tensor A$ (spacial tensor product), since
  $\boundedops(H)$ is a von Neumann 
  algebra, and (spacial) tensor products of von Neumann algebras are von Neumann
  algebras.
\end{proof}

\begin{example}
  Assume that $V=A^n$ and $W=A^m$. Then we can identify $\Hom_A(V,W)$
  with $M(n\times m,A)$, matrices acting by multiplication from the
  left. On the other hand, $M(n\times m,A)= A^{nm}$ is itself a
  Hilbert $A$-module (if $A$ is not commutative, this $A$-module
  structure is of course not compatible with the action of
  $\Hom_A(V,W)$ on the $A$-modules $V$ and $W$).

  However, as Hilbert $A$-module $\Hom_A(V,W)$ inherits the structure
  of a Banach space. The corresponding Banach norm $\abs{\cdot}$ is in
  general not equal to the 
  operator Banach norm $\norm{\cdot}$ from Lemma \ref{lem:Banach_space_Hom_A}. But it is
  always true that the two norms are equivalent. For
  $\Phi\in\Hom_A(A^n,A^m)$, represented by the matrix $(a_{ij})\in
  M(n\times m,A)$, with $e_i=(0,\dots,0,1,0,\dots,0)$ ($1$ at
  the $i$th position), and for arbitrary $v\in V$
  \begin{equation*}
    \begin{split}
      \abs{\Phi(v)} &=\abs{\sum_{j=1}^m e_j \innerprod{\Phi_j,v}}\le
      \sum_{j=1}^m \abs{\innerprod{\Phi_j,v}}\\
      &\le \sum_{j=1}^m
      \abs{\Phi_j}\cdot\abs{v} \le\sqrt{m}\abs{\Phi}\cdot\abs{v},
  \end{split}
\end{equation*}
where $\Phi_j$ is the adjoint of the $j$th row of $\Phi$.
  Since this holds for arbitrary $v\in V$,
  \begin{equation*}
    \norm{\Phi}\le \abs{\Phi}.
  \end{equation*}
  On the other hand
  \begin{equation*}
%    \begin{split}
            \abs{\Phi}^2 =\abs{\sum_{i=1}^n\sum_{j=1}^m a_{ij}^*a_{ij}}
 %     \le \left(\sum_{i=1}^n\abs{\sum_{j=1}^m a_{ij}}^2\right)^{1/2}\\
      \le \sum_{i=1}^n \abs{\sum_{j=1}^m a_{ij}^*a_{ij}}
      =\sum_{i=1}^n \abs{\Phi(e_i)} ^2\le n\norm{\Phi}^2.
%  \end{split}
  \end{equation*}
\end{example}

\begin{remark}
  In particular, when we are looking at functions defined on a smooth
  manifold with values in $\Hom_A(V,W)$, the smooth ones are
  inambiguously defined, using either of the two norms to define a
  Banach space structure on $\Hom_A(V,W)$.
\end{remark}

\begin{lemma}
  Assume that $V$ is a Hilbert $A$-module. The map
  \begin{equation*}
   \alpha\colon V\to \Hom_A(V,A); v\mapsto (x\mapsto \innerprod{v,x})
  \end{equation*}
  is an $A$-sesquilinear isomorphism. $A$-sesquilinear means that
  $\alpha(v a) = a^* \alpha(v)$ for all $v\in V$ and $a\in A$. Recall
  that $\Hom_A(V,A)$ is a left $A$-module (even a left Hilbert
  $A$-module) because of the
  identification $\End_A(A)\iso A$.
\end{lemma}
\begin{proof}
  \cite[page 13]{Lance(1995)}.
\end{proof}

\begin{definition}
  A \emph{finitely generated projective} Hilbert $A$-module $V$ is a
  Hilbert $A$-module which is isomorphic as Hilbert $A$-module to a
  (closed) orthogonal direct summand of $A^n$ for suitable $n\in\naturals$. In other
  words, there is a Hilbert $A$-module $W$ such that $V\oplus W\iso
  A^n$. The corresponding projection $p\colon A^n\to A^n$ with range
  $V$ and kernel $W$ is a projection in $M(n\times n,A)$,
  i.e.~satisfies $p=p^2=p^*$. On the other hand, the range of each
  such projection is a finitely generated projective Hilbert $A$-module.
\end{definition}

We will also consider tensor products of the modules we are
considering. Assume e.g.~that $V$ is a Hilbert $A$-module, and that
$W$ is a left $A$-module. Then we consider the algebraic tensor
product $V\tensor_A W$, still an $\End_A(V)$ left module. In general,
it would not be appropriate to consider only the algebraic tensor product,
but we would have to find suitable completions. However, we will
apply this construction only to finitely generated projective
modules, where no such completions are necessary.

\begin{example}\label{ex:V_tensor_Vdual}
  Let $V$ be a Hilbert $A$-module. Then $\Hom_A(V,A)$ is an
  $A$-$\End_A(V)$ bimodule (since $\End_A(A)\iso A$). Consequently, we
  can consider $V\tensor_A \Hom_A(V,A)$ as an $\End_A(V)$
  bimodule. It is even an algebra, with multiplication map
  \begin{equation*}
    \begin{split}
    (V\tensor_A\End_A(V,A)) \tensor_{\End_A(V)} (V\tensor_A\End_A(V,A)) & \to
    V\tensor_A\End_A(V,A)\\
    (v_1\tensor \phi_1)\tensor(v_2\tensor
    \phi_2) &\mapsto v_1 (\phi_1(v_2))\tensor \phi_2.
  \end{split}
\end{equation*}

The map $\iota\colon V\tensor_A \Hom_A(V,A)\to \End_A(V)$ which sends $v\tensor
\phi$ to the endomorphism $x\mapsto v\phi(x)$ is a ring homomorphism
which respects the $\End_A(V)$ bimodule structure.
\end{example}

\begin{definition}
  Let $X$ be a locally compact Hausdorff space. A \emph{Hilbert $A$-module
  bundle} $E$ over $X$ is a topological space $E$ with projection
  $\pi\colon E\to X$ such that each fiber $E_x:=\pi^{-1}(x)$ ($x\in
  X$) has the structure of a Hilbert $A$-module, and with local
  trivializations  $\phi\colon E|_U\xrightarrow{\iso} U\times V$ which are
  fiberwise Hilbert $A$-module isomorphisms.

  If $X$ is a smooth manifold, a \emph{smooth structure} on a Hilbert
  $A$-module bundle $E$ is an atlas of local trivializations such that
  the transition functions 
  \begin{equation*}
x\mapsto \phi_2\circ\phi_1^{-1}(x)\colon U_1\cap U_2 \to \Hom_A( V_1, V_2 )
\end{equation*}
are smooth maps with values in the Banach space $\Hom_A(V_1,V_2)$.

Given two smooth Hilbert $A$-module bundles $W$ and ${W_2}$ on $X$, then
$\Hom_A(W,{W_2})$ (constructed fiberwise) also carries a canonical smooth structure.

  A Hilbert $A$-module bundle is called \emph{finitely generated projective},
  if the fibers are finitely generated projective Hilbert
  $A$-modules, i.e.~if they are direct summands in finitely generated
  free Hilbert $A$-modules.

  We also define \emph{finitely generated projective $A$-module
    bundles} (not Hilbert $A$-module bundles!), which are locally
  trivial bundles of left $A$-modules 
  which are direct summands in $A^n$. Using a partition of unity and
  convexity of the space of $A$-valued inner products, we
  can choose a Hilbert $A$-module bundle structure on each such
  finitely generated projective $A$-module bundle.
\end{definition}

\begin{definition}
   The smooth sections of a bundle $W$ on a smooth manifold $M$ are denoted by
  $\Gamma(W)$. If $V$ is a Hilbert $A$-module, then we sometimes write
  $C^\infty(M,V):=\Gamma(M\times V)$ for the smooth sections of the
  trivial bundle $M\times V$.

  For the continuous sections we write $C(M,V)$. 

  The space of smooth differential forms is denoted
  $\Omega^*(M)=\Gamma(\Lambda^*T^*M)$. By definition, differential
  forms with values in a Hilbert $A$-module bundle $W$ are the sections
  of $\Lambda^*T^*M\tensor W$. We sometimes write $\Omega^*(M;W):=
  \Gamma(\Lambda^*T^*M\tensor W)$. Note that the wedge product of
  differential forms induces a map
  \begin{equation*}
    \Omega^p(M;W)\tensor \Omega^q(M;{W_2})\to \Omega^{p+q}(M;W\tensor {W_2}).
  \end{equation*}
\end{definition}

\begin{lemma}\label{lem:iso_with_without_innerprod}
  Given two finitely generated projective Hilbert $A$-module bundles
  $W$ and ${W_2}$ on a locally compact Hausdorff space $X$ which are
  isomorphic as $A$-module bundles, then there is an isomorphism which
  preserves the inner products as well.

  If $X$ is a smooth manifold and both bundles carry smooth structures
  and the given isomorphism preserves
  the smooth structure, we can arrange for the new isomorphism to
  preserve the smooth structure and the inner product at the same time.
\end{lemma}
\begin{proof}
  We use the property that the inclusion of the isometries into all
  invertible operators is a homotopy equivalence.

  More precisely,
  assume that $\Phi\in C(X,\Hom_A(W,{W_2}))$ is an isomorphism. Then we
  can decompose $\Phi = U \abs{\Phi}$ with
  $\abs{\Phi}(x)=\sqrt{\Phi(x)^*\Phi(x)}\in \End_A(W_x)$, using the
  fact that $\End_A(W_x)$ is a $C^*$-algebra by
  Lemma~\ref{lem:Banach_space_Hom_A} and $U(x)=\Phi
  \abs{\Phi}^{-1}$. Then $U$ and $\abs{\Phi}$ are continuous sections
  of the corresponding endomorphism bundles, and $U(x)$ is an isometry
  for each $x\in X$, i.e.~provides the desired isomorphism which
  preserves the inner products.

  Of course we use that multiplication, taking the adjoint, taking the
  inverse, and $a\mapsto \abs{a}$ are all continuous operations for
  $A$-linear adjointable operators.

  In case we have smooth structures, the isomorphism being smooth
  translates to $\Phi$ being a smooth section of $\Hom_A(W,{W_2})$. The
  new isomorphism will be smooth since all operations involved, namely
  multiplication, taking the adjoint, taking the inverse, and
  $a\mapsto\abs{a}=\sqrt{a^*a}$ are smooth, even analytic, operations for $A$-linear
  adjointable invertible operators.
\end{proof}

\begin{theorem}\label{theo:unique_smooth_structure}
  Let $V_1$ and $V_2$ be two smooth Hilbert $A$-module bundles on a
  paracompact manifold $M$ which
  are topologically isomorphic (but the isomorphism is not necessarily
  smooth). Then there is also a smooth isomorphism between the two
  bundles.

  In other words, up to isomorphism there is at most one smooth
  structure on a given Hilbert $A$-module bundle.
\end{theorem}
\begin{proof}
  An isomorphism between $V_1$ and $V_2$ is the same as a continuous
  section $s$ of the bundle $Hom_A(V_1,V_2)$ which takes values in the
  subset of invertible elements $Iso_A(V_1,V_2)$ of each fiber. The fact that
  $End_A(V_i)$ are $C^*$-algebra bundles (and a von Neumann series argument)
  shows that the invertible elements form an \emph{open} subset of
  $End_A(V_1,V_2)$.  

  The smooth structures
  on $V_1$ and $V_2$ induce a smooth structure on $End_A(V_1,V_2)$,
  and $s$ is a smooth section if and only if the corresponding bundle
  isomorphism is smooth.

  Observe that the inverse morphism $s^{-1}$ is obtained by taking
  fiberwise the inverse: $s^{-1}(x)=s(x)^{-1}$. The map
  $Iso_A(V_1,V_2)\to Iso_A(V_2,V_1);\;s\mapsto s^{-1}$ is smooth (even
  analytic), in particular
  continuous. This is the reason why it suffices to consider $s$
  alone.

  Assume for the moment that $M$ is compact. Then, to the given $s$ we
  find $\epsilon>0$ such that $\abs{s(x)-y}<\epsilon$ implies that
  $y\in Iso_A((V_1)_x,(V_2)_x)$. Using the continuity of $s$ we can
  find a finite collection $\{x_i\}\subset M$ of points, and a smooth partition of
  unity $\phi_i$ with support in some neighborhood $U_i$ of $x_i$,
  with smooth trivialization $\psi_i$ of our bundles over $U_i$, such that
  \begin{equation*}
  t(x):=\sum_i \phi_i(x) \psi_i^{-1}s(x_i)
\end{equation*}
  satisfies $\abs{t(x)-s(x)}<\epsilon$ for all $x\in M$. Observe that
  $s(x_i)$ is mapped to nearby fibers (on $U_i$) using the
  trivializations. The section $t\in \End_A(V_1,V_2)$ is by its
  definition smooth, and invertible by the choice of $\epsilon$.

  This method generalizes to paracompact manifolds in the usual way,
  replacing $\epsilon$ by  a function $\epsilon(x)>0$, and the finite
  partition of unity by a locally finite partition of unity.
\end{proof}

\subsection{Structure of finitely generated projective bundles}
\label{sec:struct-finit-gener}

By definition, finitely generated projective Hilbert
$A$-modules are direct
summands of modules of the form $A^n$. We know that, on compact
spaces, complex vector bundles are direct summands of trivial vector
bundles. We now put these two observations together.

\begin{theorem}\label{theo:structure_of_fgp_bundles}
  Let $X$ be a compact Hausdorff space and $\pi\colon W\to X$ a finitely generated
  projective Hilbert $A$-module bundle.
  \begin{enumerate}
  \item \label{item:iso_to_subbundle}
    Then $W$ is isomorphic (as Hilbert $A$-module bundle preserving
    the inner product) to a direct summand of a trivial bundle
    $X\times A^n$ for suitable $n$ (with orthogonal complement bundle
    $W^\perp$ such that $W\oplus W^\perp= X\times A^n$).
   \item \label{item:iso_to_proh_image}
    In other words, there is a projection valued function
    $\varepsilon\colon X\to M(n\times n,A)$ such that $W$ is
    isomorphic to the fiberwise image of $\varepsilon$.
   \item \label{item:proj_image_is_bundle}
    Vice versa, the image of every such projection valued function is
    a finitely generated projective Hilbert $A$-module bundle.
  \item  \label{item:proj_image_iso}
    If $\varepsilon_1$ and $\varepsilon_2$ are two projection valued
    functions as above, then, for some $\delta>0$ determined by
    $\varepsilon_1$, if $\abs{\varepsilon_1(x)-\varepsilon_2(x)}<\delta$
    for each $x\in X$, then the two image bundles are isomorphic.
  \item \label{item:smooth_proj_image}
     If $X$ is a smooth manifold and $W$ is a smooth bundle,
    then the function $\varepsilon$ can be chosen smooth. The image
    bundle inherits a canonical smooth structure, and $W$ is
    isomorphic to this bundle as a smooth bundle.
 \item \label{item:existence_of_smooth_structures}
    Every finitely generated projective Hilbert $A$-module bundle over
    a smooth compact base manifold admits
    a smooth structure. It is unique upto isomorphism.
\end{enumerate}
\end{theorem}
\begin{proof}
Assume that the situation of the theorem is given.

%\begin{itemize}
%\item
\noindent\ref{item:iso_to_subbundle}:
 Choose a
  covering $U_1,\dots, U_k$ of $X$ with trivializations
  $\alpha_i\colon W|_{U_i}\xrightarrow{\iso} U_i\times V_i$, and $\hat
  V_i$ with
  $V_i\oplus \hat V_i\iso A^n$ (of course, if $X$ is connected, all the
  $V_i$ are isomorphic). Choose a partition of unity $\phi_i^2\ge 0$
  subordinate to the covering $\{U_i\}$. Define the (isometric!) embedding 
  \begin{equation*}
j\colon W\to
  X\times (A^n)^k\colon v\mapsto  \left( \sum
    \phi_i(\pi(v))\alpha_i(v)\right)_{i=1,\cdots,n} .
\end{equation*}
Claim: the fiberwise orthogonal complements to $W$ in $X\times A^{nk}$
form a Hilbert $A$-module bundle ${W_2}$ such that
$W\oplus {W_2}=X\times A^{nk}$. To prove the claim, first of all, we can study
${W_2}$ for each component of $X$ separately, and therefore assume that
all $V_i$ are equal (to $V$ with complement $\hat V$). Secondly, it
suffices to find ${W_2}$ such that
$W\oplus {W_2}\iso X\times V^k$; then $W\oplus {W_2}\oplus
(X\times \hat V^k)\iso X\times A^{nk}$. Observe that the embedding $j$
factors through an embedding (also called $j$) 
\begin{equation*}
  j\colon W\into X\times V^k.
\end{equation*}
We claim that this embedding has an orthogonal complement $W^\perp$ with $j(W)\oplus
W^\perp\iso X\times V^k$. Therefore we can use ${W_2}:=W^\perp$ to
conclude that $W$ has a complementary Hilbert $A$-module bundle.

In contrast to Hilbert spaces, not every Hilbert $A$-submodule does have
an orthogonal complement. Therefore, we have to prove the above
assertion. Observe that there is no problem in defining the
complementary bundle $W^\perp:=\{ (x,v)\in X\times V^k\mid v\perp
j(W_x)\}$. Positivity of the inner product implies $j(W)\cap W^\perp=X\times
\{0\}$. It remains to prove that for each fiber $j(W_x)+ W^\perp_x =
V^k$. For this, observe that $j(W_x)= \{(\phi_1\alpha_1
(v),\dots,\phi_k\alpha_k(v))\mid v\in W_x\}$, with
$\phi_1,\dots,\phi_k\in \reals$ and not all $\phi_k=0$, and
$\alpha_i\colon W_x\to V$ Hilbert $A$-module isometries. Without loss
of generality, $\phi_1\ne 0$. Then
\begin{equation*}
  j(W_x) = \{ (v, \beta_2(v),\cdots,\beta_k(v))\mid v\in V\},
\end{equation*}
with $\beta_i=\phi_1^{-1}\phi_i\alpha_i\circ\alpha_1^{-1}\in \End_A(V)$. More precisely,
they are real multiples
(zero is possible) of Hilbert $A$-module isometries. Observe that an
isometry is automatically adjointable, the inverse being the adjoint.

We claim that $W_x^\perp$ is the Hilbert $A$-submodule $U_x$ of $V^k$
generated by the elements
\begin{equation*}
(-\beta_i^*(v),0,\cdots,0,v,0,\cdot,0),\quad\text{with
$v\in V$ at the $i$th position ($i=2,\dots,k$)}.
\end{equation*}
Because of the calculation of inner products
\begin{multline*}
\innerprod{(v,\beta_2(v),\dots,\beta_k(v)),(-\beta^*_i(w),0,\dots,0,w,0,\dots,0)}
= \innerprod{v,-\beta_i^*(w)}_V + \innerprod{\beta_i(v),w}_V\\
=
  -\innerprod{\beta_i(v),w}_V + \innerprod{\beta_i(v),w}_V =0
\end{multline*}
each of these elements is indeed contained in $W_x^\perp$. To show
that the sum satisfies $j(W_x)+U_x=V^k$ we have for arbitrary
$(v_1,\dots,v_k)\in V^k$ to find $w_1,\dots,w_k\in V^k$ with 
\begin{equation*}
  \begin{split}
    w_1-\beta_2^*(w_2)-\dots-\beta_k^*(w_k) & = v_1\\
    \beta_2(w_1) +w_2 & = v_2\\
    \ldots\\
    \beta_k(w_1) + w_k &= v_k.
  \end{split}
\end{equation*}
Equivalently (adding $\beta_i^*$ of the lower equations to the first one),
\begin{equation*}
  \begin{split}
    w_1 + \beta_2^*\beta_2(w_1) +\dots + \beta_k^*\beta_k(w_1) &= v_1+
    \beta_2^*(v_2)+\dots+ \beta_k^*(v_k)\\
    w_2 & = v_2-\beta_2(w_1)\\
    \ldots\\
    w_k &= v_k - \beta_k(w_1).
  \end{split}
\end{equation*}
Since $1+\beta_2^*\beta_2+\dots\beta_k^*\beta_k\ge 1$ is an invertible
element of the $C^*$-algebra $\End_A(V)$, there is indeed a (unique)
solution $(w_1,\dots,w_k)$ of our system of equations.

It remains to check that $W^\perp$ (with the $A$-valued inner product
given by restriction) is really a locally trivial bundle of Hilbert
$A$-modules. Because of our description of $W^\perp$,
$W^\perp|_{\{\phi_1\ne 0\}}\to V^{k-1}\colon
(v_1,\dots,v_k)\mapsto (v_2,\dots,v_k)$ is an isomorphism of right
$A$-modules and therefore gives a trivialization of a right $A$-module
bundle (finitely generated projective).

The transition functions (here between $\{\phi_1\ne 0\}$ and
$\{\phi_i\ne 0\}$) are given by
\begin{multline*}
  (v_2,\dots,v_k)\mapsto
  \left(-\phi_1^{-1}\phi_2(\alpha_2\circ\alpha_1^{-1})^*(v_2)\dots-\phi_1^{-1}\phi_k(\alpha_k\alpha_1^{-1})^*(v_k),
  v_2,\dots, v_k\right) \\
  \mapsto
\left(-\phi_1^{-1}\phi_2(\alpha_2\circ\alpha_1^{-1})^*(v_2) \dots-\phi_1^{-1}\phi_k(\alpha_k\circ\alpha_1^{-1})^*(v_k),
v_2,\dots,\hat{v_i},\dots, 
v_k\right).
\end{multline*}
Here, $\hat{v_i}$ means that this entry is left out. 

In particular, we observe that in the case where $X$ is a smooth
manifold and $W$ is a smooth bundle, if we choose a smooth partition
of unity, the complementary bundle
$W^\perp$ obtains a canonical smooth structure, as well. Moreover, the
inclusions of $W$ and $W^\perp$ into $X\times V^k$ are both smooth.

Then $W^\perp\oplus (X\times \hat V^k)$ also has a smooth structure,
and again the inclusions are smooth.

By Lemma
\ref{lem:iso_with_without_innerprod}, from the non-inner product
preserving trivialization of $W^\perp$ we produce
trivializations which respect the given inner product.

%\item
\noindent\ref{item:iso_to_proh_image}:
Define now $\varepsilon\colon X\to M(nk\times nk,A)=\Hom_A(A^{nk},A^{nk})$ such that
$\varepsilon(x)$ is the matrix representing the orthogonal projection from
$A^{nk}$ onto $j(W_x)$. $\varepsilon$ can be written as the composition
of three maps: the inverse
of the isomorphism $W\oplus {W_2}\to X\times A^{nk}$ which is continuous,
the projection $W\oplus {W_2}\onto W$ (which, in a local trivialization is
constant, and therefore depends continuously on $x\in X$), and the inclusion of $W$ into
$X\times A^{nk}$, which is continuous. Altogether,
$x\mapsto\varepsilon(x)$ is a continuous map.

Moreover, if $X$ and the bundle $W$ are smooth and we perform our
construction using the smooth structure, then the above argument
implies that $\epsilon$ is a smooth map.

%\item
\noindent \ref{item:proj_image_is_bundle} and \ref{item:smooth_proj_image}:
We now have to show that the images $W_x:=\im(\varepsilon(x))$ of a
projection valued map $\varepsilon\colon X\to
M(n\times n,A)$ form a finitely generated projective Hilbert
$A$-module bundle $W$, with a canonical smooth structure if $X$ and
$\varepsilon$ are smooth. Evidently, each fiber is a finitely
generated projective Hilbert $A$-module. But one still has to check
(as before for $W^\perp$) that this is a locally trivial bundle.

Fix $x_0\in X$. We claim that $\varepsilon(x_0)|_{W_x}\colon W_x\to W_{x_0}$
defines a trivialization of $W|_U$, for $U$ a sufficiently small open
neighborhood of $x_0$. To see this, precompose $\varepsilon(x_0)$ with
$\varepsilon(x)$. For $x=x_0$, this is the identity map, and it depends
continuously on $x$. Therefore
$\varepsilon(x_0)\circ\varepsilon(x)\colon W_{x_0}\to W_x\to W_{x_0}$
is an isomorphism for $x$ sufficiently small (the invertible endomorphisms
being an open subset of all endomorphisms). More precisely, if
$\abs{\varepsilon(x_0)-\varepsilon(x)}<1$, then
$\abs{\varepsilon(x_0)-\varepsilon(x_0)\varepsilon(x)}<1$
and then, since $\varepsilon(x_0)$ is the identity on $W_{x_0}$, by the
von Neumann series argument
$\varepsilon(x_0)\varepsilon(x)$ is an isomorphism. In
the same way, under the same assumption 
\begin{equation*}
\varepsilon(x)\circ\varepsilon(x_0)\colon W_x\to
W_{x_0}\to W_x
\end{equation*}
is an isomorphism. This shows that we have indeed constructed local
trivializations of $W$, which therefore is a Hilbert $A$-module bundle
(we obtain other trivialization which preserve the inner product by Lemma
\ref{lem:iso_with_without_innerprod}).

Let $\alpha_x:= (\varepsilon(x_0)\colon W_x\to W_{x_0})^{-1}\colon
W_{x_0}\to W_{x}$ be the inverse of the trivialization isomorphism
(where defined). We want to show that our trivializations define a
smooth structure on $W$ if $\varepsilon(x)$ is a smooth function. We have to show that
$\varepsilon(x_1)\circ \alpha_x\colon W_{x_0}\to W_{x_1}$ depends
smoothly on $x$ (where defined). To do this, we precompose with the
isomorphism $\varepsilon(x_0)\circ\varepsilon(x)\colon W_{x_0}\to
W_x\to W_{x_0}$. By assumption, this and therefore automatically also
its inverse depend smoothly on $x$. But the composition is
$\varepsilon(x_1)\circ \varepsilon(x)$, which again is a smooth
function of $x$. This establishes smoothness of $W$.

If we construct the subbundle $W$ and the projection $\varepsilon$ as
in \ref{item:iso_to_subbundle} and \ref{item:proj_image_is_bundle}, we
still have to check that the smooth structures coincide. The map
$\alpha_i^{-1}\colon U_i\times V\to W|_{U_I}$ is (by definition of the
smooth structure of $W$) a smooth map, the embedding $i\colon W\to
X\times A^n$ is a smooth map, and the projection
$\varepsilon(x_0)\colon W\to X\times W_{x_0}$ is a bounded linear map
which (in the coordinates just constructed) does not depend on $x$ and
therefore is also smooth. The
composition of these maps is therefore also smooth, and it is the map
which changes from the old smooth bundle chart to the new smooth
bundle chart. Therefore the inclusion gives an isomorphism of smooth
bundles between $W$ and the subbundle $i(W)$ which is the image bundle
of $\varepsilon$.

%\item
\noindent\ref{item:proj_image_iso}: Given two projection valued
  functions $\varepsilon_1$ and $\varepsilon_2$ with image bundle
  $W_1$ and $W_2$, respectively, $\varepsilon_1\colon W_2\to W_1$ will
  be an isomorphism (not preserving the inner products) if
  $\epsilon_1$ and $\epsilon_2$ are close enough because of
  exactly the
  same argument which showed in \ref{item:proj_image_is_bundle} that
  the projections can be used to get local trivializations.

%\item
\noindent \ref{item:existence_of_smooth_structures}: By Theorem
  \ref{theo:unique_smooth_structure} there is up to isomorphism at most one smooth
  structure on a Hilbert $A$-module bundle $W$. Therefore it suffices
  to prove that one smooth structure exists. To do this, embed a
  finitely generated projective Hilbert $A$-module bundle $W$ into
  $X\times A^n$ as in \ref{item:iso_to_subbundle}. Let
  $\varepsilon\colon X\to M(n,A)$ be the projection valued function
  such that the image bundle is (isomorphic to) $W$. Choose a smooth
  approximation $\varepsilon'$ to $\varepsilon$, sufficiently close
  such that the image bundles are isomorphic by
  \ref{item:proj_image_iso}. Observe that we can appoximate continuous
  functions to Banach spaces arbitrarily well by smooth function, and
  we can achieve that the new smooth function is projection valued by application
  of the holomorphic functional calculus (because of the analyticity,
   smoothness is preserved). Because $\varepsilon'$ is smooth, the
  image bundle obtains a smooth structure by
  \ref{item:smooth_proj_image}, and this does the job.
%\end{itemize}
\end{proof}

\begin{remark}
  The usual approximation results work for the infinite dimensional
  bundles we are considering: if $M$ is a compact manifold and $W$ is
  a finitely generated projective Hilbert $A$-module bundle on $X$,
  then the space of smooth sections is dense for the $C^k$-topology in
  the space of $k$-times differentiable sections.
\end{remark}

\subsection{K-theory with coefficients in a $C^*$-algebra}
\label{sec:k-theory-with}

\begin{definition}
  Let $X$ be a compact Hausdorff space and $A$ a $C^*$-algebra. The
  K-theory of $X$ with coefficients in $A$, $K(X;A)$, is defined as
  the Grothendieck group of isomorphism
classes of finitely generated projective Hilbert $A$-module bundles
over $X$.
\end{definition}

\begin{proposition}\label{prop:K_theory_of_compacts_via_bundles}
 Let $X$ be a compact Hausdorff space. Then
  $$K(X;A)\iso K_0(C(X,A)),$$ i.e.~the K-theory group of Hilbert $A$-module
  bundles is isomorphic to the K-theory of the $C^*$-algebra of
  continuous $A$-valued functions on $X$. The isomorphism is
  implemented by the map which assigns to a Hilbert bundle the module of
  continuous sections of this bundle.

  Observe also that $C(X,A)\iso C(X)\tensor A$, where we use the
  (minimal) $C^*$-algebra tensor product. (Actually, since $C(X)$ is continuous
  and therefore nuclear, there is only one option for the tensor
  product.) 
\end{proposition}
\begin{proof}
  By Theorem \ref{theo:structure_of_fgp_bundles}, every finitely
  generated projective Hilbert $A$-module bundle $W$ has a complement
  ${W_2}$ such that $W\oplus {W_2}\iso X\times A^n$ for a suitable
  $n$. Moreover, 
  \begin{equation*}
      C(X,W)\oplus C(X,{W_2})\iso C(X,W\oplus {W_2})\iso C(X,A^n) \iso (C(X,A))^n,
\end{equation*}
  i.e.~$C(X,W)$ is a direct summand in a finitely
generated free $C(X,A)$-module
and therefore is finitely generated projective.

An isomorphism $W\to {W_2}$ of Hilbert $A$-module bundles induces an
isomorphism $C(X,W)\to C(X,{W_2})$ of $C(X,A)$-modules. Moreover,
$C(X,W\oplus {W_2})\iso C(X,W)\oplus C(X,{W_2})$ as
$C(X,A)$-modules. It follows that
\begin{equation}\label{eq:def_s}
  s\colon K(X;A)\to K_0(C(X,A));\; W\mapsto C(X,W)
\end{equation}
is a well defined group homomorphism. 

We now explain how to construct
the inverse homomorphism.
Assume therefore that $L$ is a finitely generated projective
$C(X,A)$-module with complement $L'$, i.e.~$L\oplus L'=
C(X,A)^n=C(X,X\times A^n)$. Define the set
\begin{equation*}
  W:= \{ (x,v)\in X \times A^n\mid \exists s\in L; s(x)=v \}.
\end{equation*}
We claim that $W$ is a finitely generated Hilbert $A$-module bundle
with $C(X,W)\iso L$, where $\pi\colon W\to X$ is given by
$\pi(x,v)=x$. Let $p\colon C(X,X\times A^n)\to L$ be the 
projection along $L'$. We have to prove that $W$ is a locally trivial
bundle. Fix $x\in X$. Define 
\begin{equation*}
\alpha_x\colon X\times W_x\to W;\; (x,v)\mapsto (x,p(c_v)(x))  
\end{equation*}
where $c_v\in C(X,A^n)$ is the constant section with value $v\in
W_x\subset A^n$. Restricted to a sufficiently small neighborhood $U\subset
X$ of $x$, this map is an isomorphism $U\times W_x\to W|_U$. This can
be seen as follows: we compose $\alpha_x$ with the map $\beta\colon
W\to X\times W_x$ with $(y,v)\mapsto (y,p(c_v)(x))$. Then
$\beta\circ\alpha_x\colon X\times W_x\to X\times W_x$ is a continuous
section of $\End_A(X\times W_x)$ and its value at $x$ is
$\id_{W_x}$. By continuity, and since the set of invertible elements
of the $C^*$-algebra $\End_A(W_x)$ is open, $\beta\circ\alpha_x(y)$ is
invertible if $y$ is close enough to $x$. By Lemma
\ref{lem:iso_with_without_innerprod}, we can turn this into an
isomorphism $W|_U\to U\times W_x$ which preserves the inner products.

Consequently, $W$ is a Hilbert $A$-module bundle. Using the local
trivializations constructed above, we conclude also that indeed $C(X,W)=L$.

The same reasoning applies to show that $L'=C(X,W')$ with $W'$ defined
in the same way as $W$ is defined, and $C(X,W)\oplus
C(X,W')=C(X,A^n)$. From this, we conclude that $W\oplus W'=X\times
A^n$, i.e.~$W$ is a finitely generated projective Hilbert $A$-module
bundle.

Assume that $\rho\colon L\to N$ is an isomorphism of finitely generated
projective $C(X,A)$-modules. Assume that $L\oplus L'\iso A^n$ and $N\oplus
N'\iso A^m$. We can assume that there is an isomorphism $\rho'\colon
L'\to N'$ (simply replace $L'$ by $L'\oplus (N\oplus N')$ and $N'$ by
$N'\oplus (L\oplus L')$).
Then our construction shows that $\rho$ induces an isomorphism between
the Hilbert $A$-module bundles associated to $L$ and to $N$,
respectively. Similarly, the Hilbert $A$-module bundle associated to
$L\oplus N$ is the direct sum of the bundles associated to $L$ and to
$N$. Consequently, the construction defines a homomorphism
\begin{equation}\label{eq:def_t}
  t\colon  K_0(C(X,A))\to K^0(X;A).
\end{equation}

The maps $s$ of \eqref{eq:def_s} and $t$ of \eqref{eq:def_t} are by
their construction inverse to each other. This concludes the proof of
the proposition.
\end{proof}

%\begin{definition}\label{def:K_infty} Let $A$ be a $C^*$-algebra and
%  $M$ a smooth compact
%  manifold.

%  Let $K^0_\infty(M;A)$ be the subgroup of $K^0(M;A)$ consisting of
%  elements $x$ represented by a formal difference $[V_1]-[V_2]$ of
%  \emph{smooth} Hilbert $A$-module bundles over $M$. Because of
%  Theorem \ref{theo:structure_of_fgp_bundles}
%  \ref{item:existence_of_smooth_structures},
%  $K^0_\infty(M;A)=K^0(M;A)$. 
%\end{definition}

For several reasons, in particular to be able to discuss Bott
periodicity conveniently, it is useful to extend the definition of
K-theory from compact to locally compact spaces. For the latter ones,
we will restrict ourselves to compactly supported K-theory (which is
the usual definition).

\begin{definition}
  Let $X$ be a locally compact Hausdorff space. Denote its one-point
  compactification $X_+$. Then 
  \begin{equation*}
    K^0_c(X;A):= \ker( K^0(X_+;A) \to K^0(\{\infty\};A)),
  \end{equation*}
  where the map is induced by the inclusion of the additional point
  $\infty\into X_+$.
\end{definition}

\begin{proposition}
  Assume that $X$ is a locally compact Hausdorff space. Then
  \begin{equation*}
    K_c^0(X;A) \iso K_0(C_0(X;A)).
  \end{equation*}
\end{proposition}
\begin{proof}
  The split exact sequence of $C^*$-algebras
  \begin{equation*}
    0\to C_0(X;A) \to C(X_+;A)\xrightarrow{ev_\infty} A\to 0
  \end{equation*}
  gives rise to the split exact sequence in K-theory
  \begin{equation*}
    0\to K_0(C_0(X;A))\to K_0(C(X_+;A)) \to K_0(A)\to 0.
  \end{equation*}
  We know by the proof of Proposition \ref{prop:K_theory_of_compacts_via_bundles}
  that $\ker(K_0(C_0(X_+;A)\to K_0(A)))$ is given by the Grothendieck group
  of finitely generated projective Hilbert $A$-module bundles over
  $X_+$, where the fiber over $\infty$ formally is zero.
\end{proof}

As in the case of a compact space $X$, we now show that $K^0_c(X;A)$
can be described in terms of finitely generated projective bundles
over $X$.

\begin{proposition}\label{prop:K_theory_of_non_compacts_via_bundles}
  Assume that $X$ is a locally compact Hausdorff space. The group
  $K^0_c(X;A)$ is isomorphic to the group of stable isomorphism
  classes of tuples 
  \begin{equation*}
(W,{W_2}, \phi_W,\phi_{W_2})
\end{equation*}
where $W$ and ${W_2}$ are finitely
  generated projective Hilbert $A$-module bundles on $X$ and
  $\phi_W\colon W_{X\setminus K} \to (X\setminus K) \times P$,
  $\phi_{W_2}\colon W_{X\setminus K}\to (X\setminus K)\times P$ are
  trivializations of the restriction of $W$ and ${W_2}$ to the complement
  of a compact subset $K$ of $X$ where the range for both these
  trivializations is equal.

  Two such tuples $(W,\dots)$ and $(V,\dots)$  are defined to be
  stably isomorphic if there is a
  third one $(U,\dots)$ and isomorphisms $W\oplus U\to V\oplus
  U$, $W_2\oplus U_2\to V_2\oplus U_2$ such that the induced
  isomorphisms of the trivializations on the common domain of
  definition $(X\setminus K)\times (P_W\oplus
  P_U) \to (X\setminus K)\times (P_V\oplus P_U)$ extends to an isomorphism
  $(X_+\setminus K)\times (P_W\oplus P_U) \to (X_+\setminus K)\times
  (P_V\oplus P_U)$, and correspondingly for $W_2$, \ldots.
%  involving these isomorphism and the trivializations commutes on the
%  common domain of definition.

  The sum is given by direct sum, where the trivializations have to be
  restricted to the common domain of definition.
\end{proposition}
\begin{proof}
  We have shown that $K^0(X_+;A)$ is the Grothendieck group of
  finitely generated projective Hilbert $A$-module bundles over
  $X_+$. The kernel of the map to $K^0(\infty;A)$ is therefore given
  by formal differences of two Hilbert $A$-module bundles over $X_+$
  with isomorphic fibers over $\infty$. A tuple $(W,{W_2},\phi_W,\phi_{W_2})$
  gives rise to such a formal difference by extending the bundles $W$
  and ${W_2}$ to $X_+$ using the trivialization on $X\setminus K$. The
  equivalence relation on the tuples is made exactly in such a way
  that this map is well defined. On the other hand, a formal
  difference of two bundles $W$, ${W_2}$ on $X_+$ gives the first two
  entries of such a tuple by restriction to $X$, and trivializations
  $W|_{X_+\setminus K} \to (X_+\setminus K) \times W_{\infty}$,
  ${W_2}|_{X_+\setminus K} \to (X_+\setminus K)\times {W_2}_{\infty}$ together with
  an identification of $({W_2})_\infty$ with $W_\infty$ (which is possible
  since we assume that the two are isomorphic) give by restriction
  rise to the required isomorphisms. Again we see that our equivalence
  relation is made in such a way that different choices (including
  different choices of the trivializations) give rise to equivalent
  tuples.

  The maps being well defined, it is immediate from their
  definitions that they are inverse to each other.
\end{proof}

Recall that in this language it is possible to define the first
K-theory group using ``suspension'' in the following way.
\begin{definition}
  Assume that $X$ is a compact Hausdorff space. Define
  \begin{equation*}
    K^1(X;A) := K^0_c(X\times\reals;A).
  \end{equation*}
  In particular,
  \begin{equation*}
    K_1(A):= K^1(\{*\};A) = K^0_c(\reals;A).
  \end{equation*}
\end{definition}

\subsubsection{Bott periodicity}
\label{sec:bott-periodicity}

We now formulate Bott periodicity in our world of Hilbert $A$-module
bundles.

\begin{theorem}\label{theo:Bott_periodicity}
  Assume that $X$ is a compact Hausdorff space. Then there is
  an isomorphism
  \begin{equation*}
   \beta\colon K^0(X;A) \to K^0_c(X\times \reals^2;A);\; W\mapsto
    \pi_1^*W\tensor \pi_2^* B.
  \end{equation*}
  Here $B$ is the Bott generator of $K^0_c(\reals;\complexs)$. It
  corresponds under the identification with $\ker(K^0(S^2)\to
  K^0(\complexs))$ to the formal difference $H-1$ where $H$ is the
  Hopf bundle and $1$ the $1$-dimensional trivial
  bundle. $\pi_{1}\colon X\times \reals^2\to X$ and $\pi_2\colon
  X\times \reals^2\to\reals^2$ are the projections, and the tensor
  product, being a tensor product of a bundle of finitely generated
  projective Hilbert $A$-modules with a bundle of finite
  $\complexs$-vector spaces, is well defined.
\end{theorem}
\begin{proof}
  The result is of course perfectly well known. For the convenience of
  the reader we show here that the general facts about Bott
  periodicity implies that our map does the job.

  Our proof follows the idea of \cite[Exercise 9.F]{Wegge-Olsen}.
  The given map $\beta$ is functorial in $X$ and $A$. It is classical that for
  $A=\complexs$ it is the Bott periodicity isomorphism. Moreover,
  $K^0(X;A)=K_0(C(X)\tensor A) = K_0(pt; C(X)\tensor A)$, and this
  identification is compatible with $\beta$. Therefore we can
  assume that $X=\{pt\}$. Use now Morita equivalence $K_0(A)\iso
  K_0(M_n(A))$ which is induced by a (non-unital) $C^*$-algebra homomorphism
  $A\to M_n(A)$ and therefore compatible with $\beta$. For any $x\in
  K_0(A)$ we find $n\in\naturals$ and projections $p,q\in M_n(A)$ such
  that $x= [p]-[q]\in K_0(M_n(A))$, where we use Morita equivalence
  to view $x$ as an element in $K_0(M_n(A))$. Define $c_p\colon
  \complexs\to M_n(A)$ by $c_p(1)=p$. By naturality,
  $\beta(p)=c_p(\beta(1))$, i.e.~the natural transformation $\beta$ is
  determined by the specific value $\beta(1)$ for $1\in
  K_0(\complexs)$.

  Since the usual Bott periodicity homomorphism coincides with $\beta$
  on $K_0(\complexs)$ and is also a natural transformation, the two
  coincide for all $C^*$-algebras, proving the claim.
\end{proof}

\begin{remark}
  Theorem \ref{theo:Bott_periodicity} extends to locally compact
  Hausdorff spaces $X$. The proof has to be slightly modified, because
  $C_0(X)\tensor A$ is not unital, such that we haven't defined
  $K^0(\{pt\};C_0(X)\tensor A)$. Since we don't need the result in
  this paper, we omit the details.
\end{remark}

\subsection{Traces and dimensions of Hilbert $A$-modules}
\label{sec:traces-hilbert-a}

\begin{proposition}\label{prop:identify_endomorphisms}
  Assume that $V$ is a finitely generated projective Hilbert
  $A$-module. Then the map 
  \begin{equation}\label{eq:def_of_i}
    \iota\colon V\tensor_A \Hom_A(V,A) \to \End_A(V)
  \end{equation}
  of Example \ref{ex:V_tensor_Vdual}
  given by $v\tensor \phi\mapsto (x\mapsto v\phi(x))$  is a canonical
  isomorphism. Since the
  isomorphism is canonical, the corresponding result holds for any
  Hilbert $A$-module bundle $W$, i.e.
  \begin{equation*}
W\tensor_A\Hom_A(W,A) \iso
  \End_A(W).
\end{equation*}
\end{proposition}
\begin{proof}
  In general, the image of $V\tensor_A \Hom_A(V,A)$ in $\End_A(V)$ is
  (after completion) by definition the algebra of compact operators
  $K(V)$. Since $V$ is finitely generated projective,
  $K(V)=\boundedops(V)=\End_A(V)$, and it is not necessary to complete.
  
  More explicitly, recall that $V$ is a direct summand in $A^n$ and let $p\in
  \End_A(A^n)$ be the projection with image $V$. Then $\End_A(V)=
  p\End_A(A^n)p$, $\Hom_A(V,A)= \Hom_A(A^n,A)p$, and $V=p(A^n)$ can
  be considered as submodules of $\End_A(A^n)$, $\End_A(A^n,A)$
  and $A^n$, respectively.

  Then 
  \begin{multline*}
V\tensor_A \Hom_A(V,A) = pA^n\tensor_A \Hom_A(A^n,A)p = p
(A^n\tensor_A A^n)p\\
     = p\End_A(A^n)p =\End_A(V).
   \end{multline*}
   % One checks immediately that the 
The identifications are given by the
   maps we have to consider.
\end{proof}

\begin{definition}
  For each algebra $A$ let $[A,A]$ be the subspace of $A$ generated by
  commutators $[a,b]:=ab-ba$ for $a,b\in A$.
\end{definition}

\begin{lemma}\label{lem:ev_not_trace}
  Given a finitely generated projective Hilbert $A$-module $V$, there
  is a canonical linear homomorphism $ev\colon \End_A(V)\to A/[A,A]$, given by the
  composition 
  \begin{equation*}
    ev\colon \End_A(V)\xleftarrow{\iso} V\tensor_A \Hom_A(V,A)
    \xrightarrow{v\tensor \phi \mapsto \phi(v)+[A,A]}A/[A,A].
  \end{equation*}
  Since this homomorphism is canonical, it extends to a finitely
  generated projective Hilbert $A$-module bundle $W$, to give rise to
  a bundle homomorphism
  \begin{equation*}
    ev\colon \End_A(W) \to M\times (A/[A,A]).
  \end{equation*}
  This homomorphism does have the trace
  property, i.e.~for all endomorphisms $\Phi_1$ and $\Phi_2$,
  \begin{equation}\label{eq:semi_trace_property}
    ev(\Phi_1\circ \Phi_2) -ev(\Phi_2\circ\Phi_1) = [A,A] = 0 \in A/[A,A].
  \end{equation}
\end{lemma}
\begin{proof}
The first
  assertion is true by Proposition \ref{prop:identify_endomorphisms}. Observe that $va\tensor \phi$ is mapped to
  $\phi(v)a+[A,A]$, whereas $v\tensor a\phi$ is mapped to $a
  \phi(v)+[A,A]$. Clearly, $\phi(v)\cdot a- a\cdot \phi(v)\in[A,A]$.

  For the trace property, observe that for
  $\phi_1,\phi_2\in \Hom_A(V,A)$ and $v_1,v_2\in A$
  \begin{equation*}
    \iota(v_1\tensor\phi_1)\circ \iota(v_2\tensor \phi_2) = \iota(v_1
    (\phi_1(v_2))\tensor\phi_2),
  \end{equation*}
  with $\iota$ of \eqref{eq:def_of_i}. It follows that 
  \begin{equation*}
    ev(\iota(v_1\tensor\phi_1)\circ\iota(v_2\tensor\phi_2)) = \phi_2(v_1)\cdot
    \phi_1(v_2) + [A,A], 
  \end{equation*}
  whereas
  \begin{equation*}
    ev(\iota(v_2\tensor\phi_2)\circ\iota(v_1\tensor\phi_1)) = \phi_1(v_2)\cdot
    \phi_2(v_1) +[A,A], 
  \end{equation*}
  i.e.~the difference of the two elements is zero in $A/[A,A]$. Because
  $\End_A(V)$ is linearly generated (using the isomorphism $\iota$
  to
  $\Hom_A(V,A)\tensor_A V$) by elements of the form $\iota(v\tensor\phi)$,
  Equation \eqref{eq:semi_trace_property} follows.
\end{proof}

An immediate consequence of Lemma \ref{lem:ev_not_trace} is the
following Lemma.
\begin{lemma}Let $Z$ be a commutative $C^*$-algebra (e.g.~$\complexs$
  or the center of $A$). 
  Let $\tau\colon A\to Z$ be a trace, i.e.~$\tau$ is linear and
  $\tau(ab)=\tau(ba)$ for each $a,b\in A$, or, in other words, $\tau$
  factors through $A/[A,A]$. Then the composition of
  $\tau$ and $ev$ is well defined and is a $Z$-valued trace on $\End_A(V)$ for
  each finitely generated projective Hilbert $A$-module $V$, and
  correspondingly for a finitely generated projective Hilbert
  $A$-module bundle $W$ on $M$. In the latter case it extends to a linear homomorphism
  \begin{equation*}
    \tau\colon \Omega^*(M;\End_A(W)) \to \Omega^*(M;Z);\;
    \eta\tensor\Phi\mapsto \eta\tensor\tau(ev(\Phi)).
  \end{equation*}
\end{lemma}

\section{Connections and curvature on Hilbert $A$-module bundles}
\label{sec:connections}

\begin{definition}\label{def:trivial_connection}
  Let $V$ be a Hilbert $A$-module. Consider the trivialized Hilbert
  $A$-module bundle $M\times V$. For a smooth section
  $f\in\Gamma(M\times V)$, define
  \begin{equation*}
    df \in\Gamma(T^*M\tensor (M\times V))
  \end{equation*}
  by the formula (locally) $df := \sum dx_i \tensor \frac{\partial
    f}{\partial x_i}$. 
\end{definition}

\begin{definition}\label{def:connections_A_module_and_Banach}
  A \emph{connection $\nabla$ on a smooth Hilbert $A$-module bundle
    $W$} is an $A$-linear map $\nabla\colon \Gamma(W)\to
  \Gamma(T^*M\tensor W)$ which is a derivation with respect to multiplication with
  sections of the trivial bundle $M\times A$, i.e.
  \begin{equation*}
    \nabla(sf ) = s df  +  \nabla(s)f\qquad\forall
    s\in\Gamma(W),\; f\in C^\infty(M;A).
  \end{equation*}
  Here we use the multiplication $W\tensor T^*M\tensor (M\times A)\to
  W\tensor T^*M\colon s\tensor \eta\tensor f\mapsto sf\tensor
  \eta$. (In particular, elements of $A$ are considered to be of
  degree zero.)

  We say that $\nabla$ is a \emph{metric connection}, if
  \begin{equation*}
    d\innerprod{s_1,s_2} = \innerprod{\nabla s_1,s_2} +
    \innerprod{s_1,\nabla s_2}
  \end{equation*}
  for all smooth sections $s_1,s_2$ of $W$. Here, we consider
  $\innerprod{s_1,s_2}$ to be a section of the trivial bundle $M\times
  A$.

  If $L$ is only a smooth bundle of Banach spaces, a
  \emph{connection} on $L$ is a $\complexs$-linear map $\nabla\colon
  \Gamma(L)\to \Gamma(T^*M\tensor L)$ which is a derivation with
  respect to multiplication with smooth functions $f\in C^\infty(M,\complexs)$.
\end{definition}

Observe that in this sense $d$ as defined in Definition
\ref{def:trivial_connection} is a connection, the so called
\emph{trivial connection} on the trivial bundle $M\times V$, which is
actually even a metric connection with respect to the pointwise
$A$-valued inner product $\innerprod{s_1,s_1}(x)=\innerprod{s_1(x),s_2(x)}_V$.

\begin{lemma}\label{lem:connection_diff}
  Given two connections $\nabla_1,\nabla_2$ on a smooth finitely generated
  projective Hilbert $A$-module bundle
  $W$, their difference $\omega:=\nabla_1-\nabla_2$ is a $1$-form with
  values in the endomorphisms $\End_A(W)$, i.e.~a section of
  $T^*M\tensor End_A(W)$. If both connections are metric connections,
  $\omega$ takes values in the skew adjoint endomorphisms of $W$.

  The difference being an endomorphism valued $1$-form means that for
  each smooth section $s$ of $W$
  and each vector field $X$
  \begin{equation*}
    (\nabla_1)_X(s)-(\nabla_2)_X(s) = \omega(X)(s),
  \end{equation*}
  where on the right hand side the endomorphism $\omega(X)$ is applied
  fiberwise to the value of the section $s$.
\end{lemma}
\begin{proof}
  We define $\omega(X)$ by the left hand side. The expression is
  $C^\infty(M)$-linear in $X$ and $A$-linear in $s$. We have to check
  that it really defines an
  endomorphism valued $1$-form, i.e.~that $\omega(X)(s)_x$ depends
  only on $s_x$ (for arbitrary $x\in M$), or
  equivalently (because of linearity), that $\omega(X)(s)$ vanishes at
  $x$ if $s$ vanishes at
  $x$.

  Observe first that, from the multiplicativity formula for
  connections, $\omega(sf) = \omega(s)f$ for every smooth section $s$
  of $W$ and every smooth $A$-valued function $f$ on $M$.

  Secondly, using a smooth cutoff function, we can write $s=s_1+s_2$
  such that $s_1$ is supported on a neighborhood $U$ of $x$ over which $W$
  is trivial, and $s_2$ vanishes in a neighborhood of $x$. Locally,
  $W|_U\subset U\times A^n$ as a direct summand. Using this
  trivialization, we can write
  $s_1=\sum  e_if_i$ with $e_i=(0,\dots,0,1,0,\dots,0)$, and $s(x)=0$
  if and only $f_i(x)=0$ for all $i$. Extending $\omega$ (arbitrarily)
  to the complement of $W$, we can conclude that $\omega(X)(s)(x) =
  \sum \omega(X)(e_i)(x) f_i(x) =0$, if $f_i(x)=0$ for all $i$. In
  other words, $\omega(X)(s)_x =0$ if $s_x=0$, i.e.~$\omega$ is a
  $1$-form.

  Assume now that $\nabla_1$ and $\nabla_2$ are metric
  connections. Then $0 =
  \innerprod{\omega(s_1),s_2}-\innerprod{s_1,\omega(s_2)} $. Since the
  inner product is taken fiberwise, the operator $\omega(X)(x)$ is
  skew adjoint for each $x\in M$ and each vector field $X$.
\end{proof}

\begin{definition}\label{def:connection_on_pullback}
  Let $f\colon M\to N$ be a smooth map between smooth manifolds and
  $W\to N$ a smooth finitely generated projective Hilbert $A$-module bundle
  with a connection
  $\nabla$. Then we define on the pull back bundle $f^*W$ a connection
  $f^*\nabla$ in the following way:
  \begin{multline}\label{eq:pullback_connection}
    (f^*\nabla)_X((f^*s)u):=  (f^*s)(du(X)) +
    f^*(\nabla_{f_*X}(s))u\\ \forall u\in C^\infty(M),\;s\in\Gamma(W),\;X\in\Gamma(TN).
  \end{multline}
  The existence of local trivializations (and the fact that the fibers
  are finitely generated $A$-module) imply that each section
  of $f^*W$ is (locally) a $C^\infty(M)$-linear combination of sections of the form
  $(f^*s) u$ as above. By linearity, we therefore define $f^*\nabla$ for
  arbitrary sections of $f^*W$. The expression is well defined because
  $\nabla$ satisfies the Leibnitz rule 
\end{definition}

\begin{lemma}\label{lem:difference_and_pullback_connection}
  Let $f\colon M\to N$ be a smooth map and $W\to N$ a smooth finitely
  generated projective Hilbert $A$-module bundle. Assume that $\nabla$ and
  $\nabla_1$ are connections on $W$ with difference
  $\omega=\nabla-\nabla_1$. Then $f^*\nabla-f^*\nabla_1= f^*\omega$
\end{lemma}
\begin{proof}
  This follows immediately from formula \eqref{eq:pullback_connection}
  for the pullback connection and the definition of the pullback of a
  differential form.
\end{proof}

\begin{definition}
 The \emph{curvature} $\Omega$ of the connection $\nabla$ on the finitely
 generated projective Hilbert
 $A$-module bundle $W$ is the operator $\nabla\circ\nabla$.

 Here, $\nabla$ is extended to differential forms with values in $W$
 using the Leibnitz rule
 \begin{equation*}
   \nabla(\omega\tensor s) = d\omega\tensor s + (-1)^{\deg(\omega)}
   \omega\nabla(s) 
 \end{equation*}
 for all differential forms $\omega$ and all sections $s$ of $W$.
\end{definition}

\begin{proposition}\label{prop:curvature_and_its_properties}
  The curvature is a
 $2$-form with values in $\End_A(W)$. If the connection is a metric
 connection, then $\Omega$ takes values in skew adjoint $2$-forms.

 Locally, we can trivialize $W|_U\iso U\times V$. Then on $W|_U$ the
 connection $\nabla$ and a trivial connection $\nabla_V$ (depending on
 the trivialization) are given. They differ by the endomorphism valued
 $1$-form $\omega$, i.e.\ $\nabla = \nabla_V+\omega$.

 Then $\Omega= d\omega+\omega\wedge \omega$. This implies
 $d\Omega=\Omega\wedge \omega-\omega\wedge \Omega$. We use the product
 \begin{equation*}
   \begin{split}
   \Gamma&(T^*M\tensor \End_A(W))\tensor\Gamma(T^*M\tensor \End_A(W))\\
   & \to \Gamma(T^*M\tensor T^*M\tensor \End_A(W)\tensor\End_A(W))\\
   & \to
   \Gamma(\Lambda^2 T^*M\tensor \End_A(W))=\Omega^2(M;\End_A(W)).
 \end{split}
\end{equation*}
\end{proposition}
\begin{proof}
  As in the proof of Lemma \ref{lem:connection_diff}, we only have to
  show that $\Omega$ is $C^\infty(M;A)$-linear. We compute for $s\in
  \Gamma(W)$ and $f\in C^\infty(M;A)$
  \begin{equation*}
    \begin{split}
      \nabla(\nabla(sf)) & = \nabla(s\tensor df)+\nabla(\nabla(s)f) =
      s\tensor d(df) + \nabla(s)df - \nabla(s)df + 
      \nabla(\nabla(s))f\\
      = \nabla(\nabla(s))f.
  \end{split}
\end{equation*}
Here we used that $d^2=0$ by Lemma
\ref{lem:trivial_connection_curvature}. The minus sign arises since
$\nabla(s)$ is a $1$-form. From $C^\infty(M;A)$-linearity, if follows
that $\Omega$ is an endomorphism valued $2$-form, since $W$ is
finitely generated projective.

Next observe that by Lemma \ref{lem:difference_of_connections_on_forms}
\begin{equation*}
  \begin{split}
    \nabla\circ\nabla &= (\nabla_V+\omega)(\nabla_V+\omega)\\
   & =
   \nabla_V\nabla_V +\omega\nabla_V +\nabla_V\circ\omega +
   \omega\wedge\omega
   = \omega\nabla_V + d\omega - \omega\nabla_V +\omega\wedge\omega\\
   &= d\omega+\omega\wedge\omega.
\end{split}
\end{equation*}
Here we use the fact that for each $s\in\Gamma(W)$
\begin{equation*}
\nabla_V(\omega\wedge s) = d\omega\wedge s - \omega\wedge \nabla_V s
\end{equation*}
(the minus arises because $\omega$ is a $1$-form, i.e.~has odd degree).
Moreover, $\nabla_V\nabla_V=0$ by Lemma
\ref{lem:trivial_connection_curvature}, since $\nabla_V$ is by
definition a trivial connection.

Then 
\begin{equation*}
d\Omega = dd\omega + (d\omega)\wedge\omega - \omega\wedge d\omega =
(d\omega+\omega)\wedge\omega - \omega\wedge (\omega+d\omega) = \Omega\wedge\omega-\omega\wedge\Omega.
\end{equation*}

If $\nabla$ is a metric connection, then $\omega$ takes values in skew
adjoint endomorphisms by Lemma \ref{lem:connection_diff} (our
trivialization $W|_U\iso U\times V$ is a trivialization of Hilbert
$A$-modules, therefore its trivial
connection is a metric connection). The same is then true for
$d\omega$, since the skew adjoint endomorphisms form a linear subspace
of all endomorphisms. Moreover, the square $\omega\wedge\omega$ is a
two form  with values in skew adjoint endomorphisms because of
the anti-symmetrization procedure involved in the square:
\begin{equation*}
  \omega\wedge\omega(X,Y) = \omega(X)\circ\omega(Y)-\omega(Y)\circ\omega(X),
\end{equation*}
whereas
\begin{multline*}
  (\omega\wedge\omega(X,Y))^* =
  \omega(Y)^*\omega(X)^*-\omega(X)^*\omega(Y)^*\\
=
  \omega(Y)\omega(X)-\omega(X)\omega(Y)= -\omega\wedge\omega(X,Y).
\end{multline*}
\end{proof}

In the proof of Proposition \ref{prop:curvature_and_its_properties} we
have used that the curvature of a trivial connection is zero, and that
the difference of two connections is known even for the extension to
differential forms. We prove both facts now.
\begin{lemma}\label{lem:difference_of_connections_on_forms}
  If $\nabla_1-\nabla_2=\omega$ for two connections on the Hilbert
  $A$-module bundle $W$, as in Lemma
  \ref{lem:connection_diff}, then the same formula holds for the
  extension of the connection to differential forms with values in $W$,
  i.e.~the action of $\omega$ is given by the following composition:
  \begin{equation*}
    \begin{split}
      \Gamma(\Lambda^p T^*M\tensor W)& \xrightarrow{\cdot\tensor\omega\tensor\cdot}
      \Gamma(\Lambda^p T^*M\tensor T^*M\tensor\End_A(W)\tensor W)\\
      & \xrightarrow{\wedge\tensor \cdot} \Gamma(\Lambda^{p+1}T^* M\tensor W).
  \end{split}
\end{equation*}
\end{lemma}
\begin{proof}
  We only have to check that $\nabla_1+\omega$ satisfies the Leibnitz
  rule. However,
  \begin{equation*}
    (\nabla_1+\omega)(\eta\tensor s) = d\eta\tensor s +
    (-1)^{deg(\eta)} \eta\wedge\nabla_1 s +
    (-1)^{deg(\eta)}\eta\wedge\omega s,
  \end{equation*}
  for each $s\in \Gamma(W)$ and each differential form $\eta$,
  since multiplication with $\omega$ is $C^\infty(M;A)$-linear and in
  particular $C^\infty(M)$-linear.
\end{proof}

\begin{lemma}\label{lem:trivial_connection_curvature}
  For the trivial connection $d$ on a trivialized bundle $M\times V$,
  $d\circ d=0$, i.e.~the curvature is zero. 
\end{lemma}
\begin{proof}
  We compute in local coordinates for a smooth section $f$ of $M\times
  V$ 
  \begin{equation*}
d(df) = d(\sum dx_i \frac{\partial
    f}{\partial x_i}) = \sum dx_jdx_i \frac{\partial^2 f}{\partial
    x_j\partial x_i} =0,
\end{equation*}
  since $dx_idx_i=0$ and $dx_idx_j=-dx_jdx_i$.
\end{proof}

\begin{definition}\label{def:connect_on_hom}
  Connections $\nabla_W$ and $\nabla_{W_2}$ on the Hilbert $A$-module
  bundles $W$ and ${W_2}$, respectively, induce by the
  Leibnitz rule a connection $\nabla$ on the smooth bundle of Banach
  spaces $\Hom_A(W,{W_2})$ with
  \begin{equation*}
    \nabla_{W_2}(\Phi(s)) = (\nabla\Phi)(s) + \Phi(\nabla_W s)
  \end{equation*}
  for each smooth section $\Phi$ of $\Hom_A(W,{W_2})$ and each smooth section $s$ of $W$.
%  Moreover, $\Hom_A(W,{W_2})$ carries a connection which satisfies the
%  following Leibnitz rule:
%  \begin{equation*}
%    \nabla_{W_2}(\Phi(s)) = \nabla(\Phi)(s) +
%    \Phi(\nabla_W(s))\qquad\forall \Phi\in\Gamma(\Hom_A(W,{W_2})),\;s\in
%    \Gamma(W). 
%  \end{equation*}
\end{definition}

\begin{lemma}\label{lem:connect_on_product_with_c_bundle}
  Assume that $E$ is a smooth finite dimensional complex Hermitian
  vector bundle and ${W_2}$ is a
  smooth Hilbert
  $A$-module bundle with connections
  $\nabla_E$ and $\nabla_{W_2}$, respectively. The fiberwise (algebraic) tensor
  product over $\complexs$ is then a Hilbert $A$-module bundle
  $E\tensor {W_2}$, since $E$ is finite dimensional and ${W_2}$ is finitely
  generated projective. By the Leibnitz rule it carries a connection
  $\nabla_\tensor$ with
  \begin{equation*}
    \nabla_\tensor(\sigma\tensor s) =\nabla_E(\sigma)\tensor s +
    \sigma\tensor \nabla_{W_2}(s)\qquad\forall \sigma\in\Gamma(E),\; s\in\Gamma({W_2}).
  \end{equation*}
  If $\Omega_E$ is the curvature of $\nabla_E$ and $\Omega_{W_2}$ the one
  of $\nabla_{W_2}$, then
  \begin{equation*}
    \Omega_\tensor = \Omega_E\tensor \id_{W_2} + \id_E\tensor\Omega_{W_2}
  \end{equation*}
  is the curvature of $\nabla_\tensor$.
\end{lemma}
\begin{proof}
  If $V$ is a finite dimensional Hermitian $\complexs$-vector space and $W$ a
  Hilbert $A$-module, then $V\tensor W\iso W^{\dim V}$ with
  isomorphism canonical up to the choice of an orthonormal basis of $V$. This implies
  that $E\tensor {W_2}$ becomes a Hilbert $A$-module bundle in a canonical
  way. It is a standard calculation that the formula defines a
  connection on the tensor product.

  For the curvature we obtain
  \begin{equation*}
    \begin{split}
      \Omega_{\tensor} =&\nabla_{\tensor} \nabla_{\tensor} =
      (\nabla_E\tensor \id_{W_2} + \id_E\tensor\nabla_{W_2})
      (\nabla_E\tensor \id_{W_2} + \id_E\tensor\nabla_{W_2}) \\
      =& (\nabla_E \nabla_E)\tensor \id_{W_2} +\id_E\tensor
      (\nabla_{W_2} \nabla_{W_2}) \\
&+ (\nabla_E\tensor \id_{W_2})
      (\id_E\tensor \nabla_{W_2}) +(\id_E\tensor\nabla_{W_2}) (\nabla_E\tensor\id_{W_2}).
  \end{split}
\end{equation*}
Observe that operators of the form $f\tensor \id$ commute with
operators of the form $\id\tensor g$ on $E\tensor {W_2}$. Consequently,
the usual sign rule when interchanging the $1$-forms $\id_{W_2}\tensor
\nabla_{W_2}$ and $\nabla_E\tensor \id_{W_2}$ applies to give $(\id_E\tensor
\nabla_{W_2})(\nabla_E\tensor \id_{W_2}) = - (\nabla_E\tensor
\id_{W_2})(\id_E\tensor\nabla_{W_2})$. This finally implies the desired
formula
\begin{equation*}
  \Omega_\tensor = \Omega_E\tensor \id_{W_2} + \id_E\tensor\Omega_{W_2}.
\end{equation*}
\end{proof}

%\begin{definition}\label{def:connections_on_products}
%  \Kommentar{What is the structure of the tensor product left? Be
%    careful with left/right module strucutres!}
%  Given two finitely generated projective Hilbert $A$-module bundles
%  $W$ and ${W_2}$ on $M$ with 
%  connections $\nabla_W$ and $\nabla_{W_2}$, the tensor product bundle
%  $W\tensor_A {W_2}$ carries a canonical tensor product connection
%  $\nabla_\tensor$ given by
%  \begin{equation*}
%    \nabla_\tensor(s_1\tensor s_2) = \nabla_W(s_1)\tensor s_2 +
%    s_1\tensor \nabla_{W_2}(s_2)\qquad \forall s_1\in\Gamma(W),\;
%    s_2\in\Gamma({W_2}). 
%  \end{equation*}
%\end{definition}

%\begin{remark}
%  Under the isomorphism $\Hom_A(W,A)\tensor_A W\xrightarrow{\iso}
%  \End_A(W)$ of Proposition \ref{prop:identify_endomorphisms} for a
%  finitely generated projective Hilbert $A$-module 
%  bundle $W$, the two ways to define a connection as of Definition
%  \ref{def:connections_on_products} coincide.
%\end{remark}

\begin{lemma}\label{lem:pullback_curvature}
  Let $f\colon M\to N$ be a smooth map and $W\to N$ a smooth finitely
  generated projective Hilbert $A$-module bundle with connection
  $\nabla$ and curvature $\Omega$. Then the curvature of the pullback
  connection $f^*\nabla$ on the pullback bundle
  $f^*W$ is $f^*\Omega$.
\end{lemma}
\begin{proof}
  By Proposition \ref{prop:curvature_and_its_properties}, locally
  $\Omega=d\omega+\omega\wedge\omega$, where $\omega$ is the
  difference between $\nabla$ and a trivial connection.

  The pullback of a trivial connection is by Definition
  \ref{def:connection_on_pullback} trivial. By Lemma
  \ref{lem:difference_and_pullback_connection}, $f^*\omega$ therefore
  is the difference between $f^*\nabla$ and a trivial
  connection. Consequently, Proposition
  \ref{prop:curvature_and_its_properties} implies that the curvature
  $\Omega^*$ of $f^*\nabla$  is given by
  \begin{equation*}
    \Omega^* = d(f^*\omega) +f^*\omega\wedge f^*\omega =
    f^*(d\omega+\omega\wedge\omega) = f^*\Omega.
  \end{equation*}
\end{proof}

\section{Chern-Weil theory}
\label{sec:chern-weyl-theory}

The prototype of the characteristic classes we want to define is the
Chern character. 
\begin{definition}\label{def:cherncharacter}
  Consider the formal power series $\exp(x)=\sum_{k=0}^\infty
  \frac{x^k}{k!}$. A differential form of degree $\ge 1$ with values
  in a ring  on a finite
  dimensional manifold can 
  be substituted for $x$.

  In particular, if $W$ is a Hilbert $A$-module bundle on a manifold
  $M$ with connection $\nabla$ and curvature $\Omega\in
  \Omega^2(M;\End_A(W))$, we define 
  \begin{equation*}
    \exp(\Omega):= \sum_{n=0}^\infty
    \frac{\overbrace{\Omega\wedge\cdots\wedge\Omega}^{k\text{
          times}}}{k!}\in \Omega^{2*}(M;\End_A(W)).
  \end{equation*}
  Given a commutative $C^*$-algebra $Z$ and a trace $\tau\colon A\to
  Z$, if $W$ is a finitely
  generated projective Hilbert $A$-module bundle, we now define
  \begin{equation*}
    \ch_\tau(\Omega):= \tau(ev(\exp(\Omega))) \in \Omega^{2*}(M;Z),
  \end{equation*}
  using the homomorphism $ev$ of Lemma \ref{lem:ev_not_trace}.
\end{definition}

\begin{lemma}\label{lem:Chern_closed}
  If $\tau$ is a trace then the characteristic class
  $\ch_\tau(\Omega)$ of Definition 
  \ref{def:cherncharacter} is closed. The cohomology class represented
  by $\ch_\tau(\Omega)$ does not depend on the connection $\nabla$ but only on the
  finitely generated projective Hilbert $A$-module bundle $W$.
\end{lemma}
\begin{proof}
  Recall that by Proposition \ref{prop:curvature_and_its_properties}
  we have 
  locally $d\Omega=\Omega\wedge \omega-\omega\wedge\Omega$ for a
  suitable endomorphism valued $1$-form $\omega$. It suffices to check
  that for each $k\in\naturals$
  \begin{equation*}
    d\tau(ev(\Omega^k)) = 0.
  \end{equation*}

  We will show that $d\tau(ev(\eta))=\tau(ev(\nabla\eta))$ for each
  $\eta\in \Omega^*(M;\End_A(W))$. This holds for
  an arbitrary connection $\nabla$, consequently we can apply it using
  (locally) the connection $d$ obtained from a trivialization. Once this is
  established, we compute (locally) and using that $\tau\circ ev$ has
  the trace property and that $\Omega$ is a form of even degree,
  \begin{equation*}
    \begin{split}
      d\tau(ev(\Omega^k)) & = \tau(ev(\nabla(\Omega^k))) = \sum_{i=0}^{k-1} \tau(ev(\Omega^i\wedge (\nabla\Omega)\wedge
      \Omega^{k-i-1}))\\
      &= \sum_{i=0}^{k-1}
      \tau(ev(\Omega^i\wedge(\Omega\wedge\omega-
      \omega\wedge\Omega)\wedge\Omega^{k-i-1}))\\
      &= \sum_{i=1}^k \tau(ev(\Omega^i\wedge\omega\wedge
      \Omega^{k-i})) - \sum_{i=0}^{k-1}
      \tau(ev(\Omega^i\wedge\omega\wedge\Omega^{k-i}))\\
      &= k \left(\tau(ev(\Omega^{k}\wedge\omega)) -
        \tau(ev(\Omega^k\wedge\omega))\right) = 0.
  \end{split}
\end{equation*}

To establish the formula $d\tau(ev(\eta))=\tau(ev(\nabla\eta))$ which
we have used above,
%Consequently, we have to show that
%  $\tau(d(ev(\Omega^k)))=0$.
%  Now, $d(ev(\eta))=ev(\nabla\eta)$ for every $\eta\in
%  \Gamma(\Lambda^*T^*M\tensor \End_A(W))$. To see this, 
it suffices to consider $\eta= \alpha \phi\tensor v$ with $\alpha\in
\Omega^*(M)$, $\phi\in\Gamma(\Hom_A(W,A))$ and $v\in \Gamma(W)$. This
is the case because such forms locally generate
$\Omega^*(M;\End_A(W))$, using the isomorphism of Proposition
\ref{prop:identify_endomorphisms}. Then,
on the one hand by Definition \ref{def:connect_on_hom}
  \begin{equation*}
    \begin{split}
      d(\tau(ev(\eta))) & = d(\tau(\alpha \phi(v)))=\tau((d\alpha) \phi(v) +
      (-1)^{\deg(\alpha)} \alpha \wedge d(\phi(v)))\\
       & = \tau((d\alpha) \phi(v) +
      (-1)^{\deg(\alpha)} \alpha \wedge ((\nabla\phi)(v)+\phi(\nabla v))).
  \end{split}
\end{equation*}
  Here, we used that the homomorphism $\tau\colon M\times A\to
  M\times Z$ is
  given by fiberwise application of $\tau\colon A\to Z$. It
  follows that $d\tau (\beta)=\tau d\beta$ for each $\beta\in
  \Omega^*(M;A)$, where we use $d\colon\Omega^*(M;A)\to
  \Omega^*(M;A)$ as defined in Definition
  \ref{def:trivial_connection}. 

  On the other hand, 
  \begin{equation*}
    \begin{split}
      \tau(ev(\nabla\eta)) & = \tau(ev((d\alpha) \phi\tensor v +
      (-1)^{\deg(\alpha)}\alpha \nabla(\phi\tensor v)))\\
      & = \tau((d\alpha) \phi(v) +      (-1)^{\deg(\alpha)} ev(\alpha\wedge
      (\nabla\phi)\tensor v+\phi\tensor \nabla(v)) )\\
      & = \tau((d\alpha) \phi(v) +
      (-1)^{\deg(\alpha)} \alpha \wedge ((\nabla\phi)(v)+\phi(\nabla v))).
  \end{split}
\end{equation*}

We now have to check that the cohomology class is unchanged if we
replace $\nabla$ by a second connection $\nabla_1$. 

Consider the projection $\pi\colon M\times [0,1]\to M$ and pull the bundle $W$
back to $M\times [0,1]$ using this projection. Using the fact that the
space of connections is convex, we equip $\pi^*W$ with a connection $\nabla_b$
which, when restricted (i.e.~pulled back) to $M\times\{0\}$ gives
$\nabla$, and when restricted to $M\times\{1\}$ gives $\nabla_1$.

By Lemma \ref{lem:pullback_curvature}, if $\Omega_b$ is the curvature
of $\nabla_b$, then its restriction to $M\times\{0\}$ is the curvature
$\Omega$ of $\nabla$, and its restriction to $M\times \{1\}$ is the
curvature $\Omega_1$ of $\nabla_1$. Application of $\tau$, $ev$ and
$\exp$ commutes with
pullback. Therefore,
\begin{equation*}
   \ch_\tau(W;\nabla)= i_0^*(\tau(ev(\exp(\Omega_b)))),\quad \text{and }
   \ch_\tau(W;\nabla_1)= i_1^*(\tau(ev(\exp(\Omega_b)))),
 \end{equation*}
 where $i_0, i_1\colon M\to M\times [0,1]$ are the inclusions
 $m\mapsto (m,0)$ and $m\mapsto (m,1)$ respectively. Since these maps
 are homotopic, the two cohomology classes represented by the two differential
 forms are equal.

This finishes the proof of the lemma.
\end{proof}

\begin{remark}
  Recall that the Chern character determines the
  rational Chern
  classes (and of course also vice versa). Therefore, the definition
  of $\ch_\tau(W)\in H^{2*}(X;Z)$ immediately gives rise also
  to Chern classes $c_{i,\tau}(W) \in H^{2i}(X;Z)$. They
  can then be used to define all other kinds of
  characteristic classes. We are not going to use this in this paper
  and therefore refrain from any further discussion.
\end{remark}

\begin{theorem}\label{theo:Chern_character_and_Bott_periodicity}
  The Chern character is compatible with Bott periodicity in the
  following sense: given a smooth finitely generated projective Hilbert
  $A$-module bundle $W$ on a compact manifold $M$ and a trace
  $\tau\colon A\to Z$, the cohomology classes 
  \begin{equation*}
\ch_\tau(W)\in
H^{2*}(M;Z)\qquad\text{and}\qquad\int_{\reals^2}\ch_\tau(\beta(W)) \in
H^{2*}(M;Z)
\end{equation*}
are equal.

  Here, $\ch_\tau(\beta(W)) = \ch_\tau(W\tensor B_+) -\ch_\tau(
  W\tensor B_-) \in H^*_c(X\times\reals^2;Z)$, where $[B_+]-[B_-]=B\in
  K^0_c(\reals^2)$ is the Bott virtual bundle on $\reals^2$ of Theorem
  \ref{theo:Bott_periodicity}. The construction of $\ch_\tau$,
  together with the proof of all its properties, immediately
  generalizes from compact base manifolds to the present case. We
  simply have to use on the two bundles two connections which coincide
  near infinity (using the given isomorphism between $B_+$ and $B_-$
  near infinity) to produce a compactly supported closed form on
  $X\times\reals^2$ representing a well defined element in compactly
  supported cohomology $H^*_c(X\times\reals^2;Z)$.

  The map $\int_{\reals^2}\colon H^*_c(X\times\reals^2;Z) \to
  H^{*-2}(X;Z)$ is the usual integration over the fiber
  homomorphism (tensored with the identity on $Z$), which in terms of
  de Rham cohomology is given by
  integration over the fibers of the product $X\times \reals^2$.
\end{theorem}
\begin{proof}
  To prove the result, on $W\tensor B_+$ and $W\tensor B_-$ we choose
  product connections. By Lemma
  \ref{lem:connect_on_product_with_c_bundle} we then obtain for the
  curvature 
  \begin{equation*}
\Omega_{W\tensor
    B_+} = \Omega_W\tensor \id_{B_+} + \id_{W}\tensor
  \Omega_{B_+}.
\end{equation*}
Since the two summands commute,
  \begin{multline}\label{eq:exp_of_tensor_product}
    \exp(\Omega_{W\tensor B_+}) = \exp(\Omega_W\tensor\id_{B_+})\wedge
   \exp(\id_W\tensor\Omega_{B_+}) \\ =
   (\exp(\Omega_W)\tensor\id_{B_+})\wedge (\id_W\tensor \exp(\Omega_{B_+})).
 \end{multline}
 Consequently, we have to study 
 \begin{equation*}
\tau\left(ev((a (\phi\tensor v)\tensor
 \id_{B_+})\wedge b \id_W \tensor(\psi\tensor u))\right)
\end{equation*}
with
 $a,b\in\Omega^*(M)$, $\phi\in\Gamma(\Hom_A(W,A))$, $v\in\Gamma(W)$,
 $\psi\in\Gamma(\Hom_\complexs(B_+,\complexs))$, $u\in
 \Gamma(B_+)$. We obtain
 \begin{equation*}
   \begin{split}
    & \tau\left(ev((a (\phi\tensor v)\tensor \id_{B_+})\wedge b \id_W
     \tensor(\psi\tensor u))\right) = \tau\bigl(ev(a\wedge b (\phi\tensor
     v)\tensor (\psi\tensor u))\bigr)\\
     & = a\wedge b \tau(\phi(v)\cdot \psi(u)) = a\tau(\phi(v)) \wedge
     b \psi(u)\qquad\text{(observe that $\psi(v)\in\complexs$)}.
\end{split}
\end{equation*}
Recall that $\psi(u)$ is the ordinary fiberwise trace $\tr$ of the
endomorphism valued section
corresponding to
\begin{equation*}
\psi\tensor u \in
\Gamma(\Hom_\complexs(B_+,\complexs))\tensor\Gamma(B_+)\iso
\Gamma(\End_\complexs(B_+)).
\end{equation*}
We obtain for general endomorphism valued forms of the form
$(\omega\tensor\id_{B_+})\wedge (\id_W\tensor \eta)$ with $\omega\in
\Omega^p(M; \End_A(W))$ and $\eta\in \Omega^q(M;\End_\complexs(B_+))$
(since the special ones considered above locally span the space of
such sections)
\begin{equation*}
  \tau\left(ev((\omega\tensor\id_{B_+})\wedge (\id_W\tensor \eta))\right) =
  \tau(ev(\omega))\wedge \tr(\eta). 
\end{equation*}
In particular, applying this formula to
\eqref{eq:exp_of_tensor_product}, we get
\begin{equation*}
  \begin{split}
    \ch_\tau(W\tensor B_+) & = \tau\left(ev(\exp(\Omega_{W\tensor B_+}))\right) =
    \tau\left(ev(\exp(\Omega_W))\wedge \tr(\exp(\Omega_{B_+}))\right)\\
    & = \ch_\tau(W) \wedge \ch(B_+),
\end{split}
\end{equation*}
where $\ch(B_+)$ is the ordinary real differential form representing
the Chern character. It follows that
\begin{equation*}
  \ch_\tau( W\tensor B_+) -\ch_\tau(W\tensor B_-) = \ch_\tau(W)\wedge(\ch(B_+)-\ch(B_-)),
\end{equation*}
where the factor $\ch(B_+)-\ch(B_-)$ is a compactly supported closed
$2$-form on $\reals^2$ representing the Chern character $\ch(B)=c_1(B)$ of
the virtual bundle $B$ (note that this is a compactly supported closed
differential form of even degree on $\reals^2$, and the $0$-degree part
is zero). Therefore, by Fubini's theorem
\begin{equation*}
  \int_{\reals^2} (\ch_\tau(W\tensor B_+)-\ch_\tau(W\tensor B_-)) =
  \ch_\tau(W)\cdot \int_{\reals^2}(\ch(B)).
  \end{equation*}
  A fundamental property of the Bott bundle is that
  $\int_{\reals^2}(\ch(B))=1$, and this concludes the proof.
\end{proof}

An important question in the classical theory of characteristic
classes is the group where the characteristic classes live in, in
particular integrality results. We know e.g.~that the degree $2n$-part
of the Chern character after multiplication with $n!$ belongs to the
image of integral cohomology in de Rham cohomology. In our situation,
the result can not be as easy as that and depends on the trace, as is
evident from the fact that the degree zero-part is equal to the
$\tau$-dimension of the fiber of the Hilbert $A$-module bundle (a
locally constant function). Only after restriction to particular
choices of bundles and particular choices of traces, meaningful
restriction can be expected. This will not be discussed in this paper.

\section{Index and KK-theory}
\label{sec:index-kk-theory}

In this section, we give our proofs of the index theorems for operators twisted
with Hilbert $A$-module bundles. We do this by using heavily the
machinery of K-theory and KK-theory to reduce to the classical
Atiyah-Singer index theorem. The main tool we will use is the
K{\"u}nneth theorem and associativity of the Kasparov product.

For this paper, we want to avoid all technicalities about Kasparov's
bivariant KK-theory for $C^*$-algebras. We will just recall a few
basic facts to be used in here. Detailed expositions can be found in
Kasparov's original paper \cite{MR81m:58075}, or in \cite{MR99g:46104}.

We consider KK to be an additive category whose objects are the $C^*$-algebras,
and with morphism sets $KK(A,B)$. There is a functor from the category
of $C^*$-algebras to the category KK which is the identity on objects,
i.e.~every $C^*$-algebra morphism $f\colon A\to B$ gives rise to an
element $[f]\in KK(A,B)$.

We define $KK_0(A,B):=KK(A,B)$ and $KK_1(A,B):= KK(SA,B)$, where $SA:=
C_0(\reals)\tensor A$ is the suspension of $A$.

We have the following properties:
\begin{proposition}\strut

  \begin{enumerate}
  \item $KK(A,\complexs)$ is the K-homology of the $C^*$-algebra
    $A$, $KK(\complexs,A)$ its K-theory (defined in terms of
    projective finitely generated modules). In particular, if $X$ is a
    compact Hausdorff space, then $KK(C(X),\complexs)=K_0(X)$ and
    $KK(\complexs,C(X))=K^0(X)$ are the usual K-homology and K-theory of the
    space $X$.
  \item An elliptic differential operator $D$ on a smooth compact manifold
    $M$ of dimension congruent to $i$ modulo $2$ defines an element $[D]$
    in $KK_i(C(M),\complexs)$. (The
    general idea is that the KK-groups are defined as equivalence classes of
    generalized elliptic operators.)
  \item On the other hand, every smooth complex vector bundle $E$ on
    an even dimensional manifold $M$ defines an element $[E]$ in
    $KK(\complexs,C(M))$. If $D$ is a (generalized) Dirac operator,
    then the composition product $[E]\circ [D]\in
    KK(\complexs,\complexs)= K_0(\complexs)=\integers$ equals the
    Fredholm index $\ind(D_E)$ of the operator $D$ twisted by the
    bundle $E$.
  \item There is an exterior product
    \begin{equation*}
  KK(A_1,B_2)\tensor
    KK(A_2,B_2)\to KK(A_1\tensor A_2, B_1\tensor B_2),
  \end{equation*}
  where we use the minimal (spacial) tensor product throughout.

  This exterior product commutes with the composition product of the category, i.e.~if
  we have $f_i\in KK(A_i,B_i)$, $g_i\in KK(B_i,C_i)$ for $i=1,2$ then
  \begin{equation*}
    (f_1\circ g_1)\tensor (f_2\circ g_2) = (f_1\tensor f_2)\circ
    (g_1\tensor g_2).
  \end{equation*}
\item Let $Z$ be a commutative $C^*$-algebra, e.g.~$Z=\complexs$. Any
  trace $\tr\colon A\to Z$, i.e.~a continuous linear
  map with $\tr(ab)=\tr(ba)$ for each $a,b\in A$ induces a
  homomorphism of abelian groups, denoted with the same letter,
  \begin{equation*}
    \tr\colon K(A)=KK(\complexs,A)\to Z.
  \end{equation*}
% \item Every manifold $M$ with a spin-c structure, in particular every
%   manifold with a spin structure, has a K-theory orientation. About
%   this K-theory orientation, we only use that it gives rise to a
%   canonical element, the fundamental class $[M]\in
%   KK_{\dim(M)}(C(M),\complexs)$ ($\dim(M)$
%   has to be interpreted modulo $2$). This
%   class is represented by the spin-c Dirac operator on $M$. 
  \end{enumerate}
\end{proposition}

\begin{definition}\label{def:ind_D_V}
  Let $D$ be an elliptic differential operator on a closed smooth
  manifold $M$, and $W$ a smooth Hilbert $A$-module bundle over $M$. We
  define the index of $D$ twisted by $W$
  \begin{equation*}
    \ind_A(D_W):=\ind(D_W):= [W] \circ ([D]\tensor [\id_A]) \in KK(\complexs, A).
  \end{equation*}
  Observe for this definition that $[D]\in KK(C(M),\complexs)$,
  $\id_A\in KK(A,A)$, $[D]\tensor [\id_A] \in KK(C(M)\tensor A,A)$ and
  $[W]\in KK(\complexs,C(M,A))$. We also use the fact that $C(M)\tensor
  A=C(M,A)$.

  Given a trace $\tau\colon A\to Z$, from this we can define a ``numerical'' index
  \begin{equation*}
    \ind_{\tau}(D_W):= \tau(\ind(D_W))\in Z.
  \end{equation*}
\end{definition}

\begin{theorem}\label{theo:Kuenneth}
 Let $X$ be a compact Hausdorff space and $A$ a
  separable $C^*$-algebra.
  There is an exact sequence
  \begin{multline*}
    0\to K^0(X)\tensor K_0(A)\oplus K^1(X)\tensor K_1(A) \to
    K^0(X;A)\\ \to \Tor(K^0(X),K_1(A))\oplus \Tor(K^1(X),K_0(A))\to 0.
  \end{multline*}
  The restriction of the first map to the summand $K^0(X)\tensor
  K_0(A)$ sends $[E]\tensor [P]$ to $[E\tensor P]$, i.e.~we tensor the
  complex finite dimensional vector bundle $E$ (over $\complexs$) with
  the Hilbert $A$-module $P$ (considered as the trivial bundle
  $X\times P$).

  The restriction of the first map to the summand $K^1(X)\tensor
  K_1(A)=K^0_c(X\times\reals)\tensor K^0_c(\reals;A)$ is given by the
  exterior tensor product as above, producing an element in
  $K^0_c(X\times\reals\times\reals;A)$, which then has to be mapped to
  $K^0(X;A)$ by the inverse of the Bott isomorphism of Theorem
  \ref{theo:Bott_periodicity}.

The sequence implies in particular that, after tensoring with
$\rationals$,  
\begin{equation}\label{eq:kuenneth}
K^0(X;A)\tensor\rationals\iso K^*(X)\tensor K_*(A)\tensor\rationals.
\end{equation}
\end{theorem}

\begin{proposition}\label{prop:Kuenneth}
  If $A$ is a finite  von Neumann algebra,
  e.g.~$A=\NeumannN\Gamma$ for a discrete group $\Gamma$, then
  $K_0(A)$ is torsion free and $K_1(A)=0$. In particular, for each
  compact Hausdorff space $X$ we have an isomorphism
  \begin{equation*}
    K^0(X)\tensor K_0(A)\stackrel{\iso}{\to} K^0(X;A).
  \end{equation*}
\end{proposition}
\begin{proof}
  By \cite[7.1.11]{MR99g:46104}, $K_1(A)=\{0\}$ for an arbitrary von
  Neumann algebra $A$. For a finite von Neumann algebra, the canonical
  center valued trace induces an injection $\tr_{Z(A)}\colon K_0(A)\to
  Z(A)$. Since the latter is a vector space, $K_0(A)$ is torsion
  free. Then we apply the exact sequence of Theorem
  \ref{theo:Kuenneth}. 
\end{proof}

\begin{remark}
  Observe that there is no explicit formula for the inverse of this
  isomorphism. Our work with connections and curvature in the previous
  sections is motivated by the attempt to overcome this difficulty.
\end{remark}

\begin{proposition}
\label{prop:explicit_formula_for_chern}  Let $\tau\colon A\to Z$ be a trace on
$A$ with values in a commutative $C^*$-algebra $Z$. If $A$ is a finite
von Neumann algebra, consider the composition
  \begin{equation*}
   \psi_\tau \colon K^0(X;A)\xleftarrow{\iso} K^0(X)\tensor
    K_0(A) \xrightarrow{\ch\tensor\tau}
    H^{ev}(X;\rationals)\tensor Z=H^{ev}(X;Z).
  \end{equation*}
  If $Z$ is the center of the finite von Neumann algebra $A$ and
  $\tau$ is the canonical center valued trace, then this map is
  rationally injective:
  \begin{equation*}
    \psi_\tau \colon K_0(X;A)\tensor\rationals \into H^{ev}(X; Z).
  \end{equation*}
  For arbitrary $A$, the map is defined at least after tensoring with $\rationals$:
  \begin{multline*}
   \psi_\tau \colon K^0(X;A)\tensor\rationals\xleftarrow{\iso} K^0(X)\tensor
    K_0(A)\tensor\rationals \oplus K^1(X)\tensor K_1(A)\tensor
    \rationals\\
    \xrightarrow{(\ch\tensor\tau) \circ \pr_1}
    H^{ev}(X;Z)\tensor\rationals.
  \end{multline*}
  If $W$ and ${W_2}$ are smooth finitely generated
  projective Hilbert $A$-module bundles on $M$ with connections
  $\nabla_W$ and $\nabla_{W_2}$, respectively, then
  \begin{equation*}
    \ch_\tau(W)-\ch_\tau({W_2}) = \psi_\tau([W]-[{W_2}]).
  \end{equation*}
\end{proposition}
\begin{proof}
  The map $W\mapsto\ch_\tau(W)$ induces a well defined homomorphism
  \begin{equation*}
\ch_\tau\colon K^0(X;A)\to 
    H^{ev}(X; Z)
  \end{equation*}
because of the following
    observations:
 
    Assume that $W_1$ and $W_2$ are finitely generated projective Hilbert
    $A$-module bundles. We can give them a (unique) smooth structure
    by Theorem \ref{theo:structure_of_fgp_bundles}. Equipping them
    with connections $\nabla_{W_1}$ and
    $\nabla_{W_2}$, respectively, then, by using on $W_1\oplus W_2$
    the connection $\nabla_{W_1}\oplus\nabla_{W_2}$, we see that
    $\ch_\tau$ is additive wit respect to direct sum. Two smooth bundles $W$, ${W_2}$ represent the same
    K-theory element if and only if  they are stably isomorphic, i.e.~if $W\oplus
    M\times V\iso {W_2}\oplus M\times V$. By Theorem
    \ref{theo:unique_smooth_structure}, we can assume this isomorphism
    to be a smooth isomorphism. By \ref{lem:Chern_closed} $\ch_\tau$
    is independent of the connection chosen. Together with additivity,
    $\ch_\tau(W)=\ch_\tau({W_2})$.
    
    Since for a finite von Neumann algebra $A$ the map $K^0(X)\tensor
    K_0(A)\to K^0(X;A)$ is an isomorphism by Proposition
    \ref{prop:Kuenneth}, and for general $A$ the map 
    \begin{equation*}
K^0(X)\tensor
    K_0(A)\tensor\rationals \oplus 
    K^1(X)\tensor K_1(A)\tensor \rationals \to
    K^0(X;A)\tensor\rationals
  \end{equation*}
is an isomorphism by Theorem
    \ref{theo:Kuenneth}, it suffices to consider the following two
    cases:
%  \begin{enumerate}
 % \item 
  
 First, we consider a bundle $E\tensor V$ where $E$ is a finite dimensional complex
    vector bundle over $M$ and $V$ is a finitely generated projective
    Hilbert $A$-module. A connection $\nabla$ on $E$ and the trivial
    connection on $M\times V$ induce the tensor product connection on
    $E\tensor V$ by Lemma \ref{lem:connect_on_product_with_c_bundle}.
    
    The calculations in the proof of Theorem
    \ref{theo:Chern_character_and_Bott_periodicity} show that 
    \begin{equation*}
      \ch_\tau(E\tensor V) =\ch(E) \cdot \tau(V),
    \end{equation*}
    since $V$ is a ``bundle'' on the one-point space and in this case
    $\ch_\tau(V)=\tau(V)\in Z$. This shows that $\psi_\tau$
    coincides with $\ch_\tau$ on $K^0(X)\tensor K_0(A)$, or on the
    summand $K^0(X)\tensor K_0(A)\tensor\rationals$, respectively.
  
 Secondly, we have to consider elements which under Bott
    periodicity correspond to
    $E\tensor V$ where $E\in K^0_c(X\times \reals)$ is a finite
    dimensional virtual vector bundle over $X\times\reals$ which is
    zero at infinity, and $V\in K^0_c(\reals;A)$ is a virtual finitely
    generated projective Hilbert $A$-module bundle which is zero at
    infinity (such virtual bundles are by definition tuples as in
    Proposition \ref{prop:K_theory_of_non_compacts_via_bundles}).

    By Theorem \ref{theo:Chern_character_and_Bott_periodicity}, we
    have to show that $\ch_\tau(E\tensor V)=0$. The proof of Theorem
    \ref{theo:Chern_character_and_Bott_periodicity} shows that 
    \begin{equation*}
      \ch_\tau(E\tensor V) =\ch(E)\wedge\ch_\tau(V),
    \end{equation*}
    with $\ch(E)\in H^{2*}_c(X\times\reals;\reals)$ and $\ch_\tau(V)\in
    H_c^{2*}(\reals;Z)$, and where the product is an ``exterior'' wedge
    product (i.e.~one first has to pull back to the product
    $X\times\reals\times\reals$). However, in even degrees the
    compactly supported cohomology of $\reals$ vanishes, therefore the
    whole expression is zero as we had to show.
%  \end{enumerate}
\end{proof}

The importance of Proposition \ref{prop:explicit_formula_for_chern}
lies in the explicit formula, where it is not necessary to invert the
isomorphism of Proposition \ref{prop:Kuenneth}. We get for instance
the following immediate corollary.

\begin{corollary}\label{corol:ch_of_flat_bundle}
  Assume that $W$ is a flat finitely generated projective Hilbert
  $A$-module bundle over the connected manifold $M$ with typical fiber
  $V$. Then
  \begin{equation*}
    \ch_\tau(W)=\psi_\tau([W]) = \psi_\tau([M\times V]) = \dim_\tau(V) \in H^0(M;Z)
  \end{equation*}
  for each trace $\tau$ on $A$, i.e.~the K-theory class represented by
  $W$ can not be distinguished from the K-theory class represented by
  the trivial bundle using these traces. $\dim_\tau(V)$ is the zero
  dimensional cohomology class represented by the (locally) constant
  function $\dim_\tau(V)$.
\end{corollary}

\subsection{The Mishchenko-Fomenko index theorem}
\label{sec:mishch-fomenko-index}

We are now ready to reprove the cohomological version of the
Mishchenko-Fomenko index theorem. Our goal is to give a
(cohomological) formula for
$\ind_\tau(D_W)$ as defined in Definition \ref{def:ind_D_V}. 

\begin{theorem}\label{theo:index_theorem}
  Assume that $M$ is a closed smooth manifold, $D$ an elliptic differential
  operator defined between sections of  finite dimensional bundles over
  $M$. Let $W$ be a finitely generated projective Hilbert $A$-module
  bundle, and $\tau\colon A\to Z$ a trace on $A$ with values in an
  abelian $C^*$-algebra $Z$. Then
  \begin{equation}\label{eq:index_formula}
    \ind_\tau (D_W) = \innerprod{ \ch(\sigma(D)) \cup \Td(T_\complexs
      M) \cup
      \ch_\tau(W) , [TM]}.
  \end{equation}
  Here, $\ch(\sigma(D))$ is the Chern character of the symbol of $D$, a
  compactly supported (real) cohomology class on the manifold $TM$,
  $\Td(T_\complexs M)$ is the Todd class of the complexified tangent
  bundle, pulled back to $TM$ and $\ch_\tau(W)$ is the pull back of
  $\ch_\tau(W)$ to $TM$. $\innerprod{\cdot,\cdot}$ stands for the
  pairing of the compactly supported cohomology class with the locally
  finite fundamental homology class $[TM]$.

  If $M$ is oriented of dimension $n$, then integration over the fibers of $\pi\colon
  TM\to M$ immediately gives the following consequence:
  \begin{equation}\label{eq:oriented_formula}
    \ind_\tau (D_W) = (-1)^{n(n-1)/2}\innerprod{ \pi_!\ch(\sigma(D))
      \cup \Td(T_\complexs M)\cup
      \ch_\tau(W) , [M]}.
  \end{equation}
  The sign compensates for the difference between the orientation of
  $TM$ induced from $M$ and its canonical orientation as a symplectic manifold.
\end{theorem}

\begin{remark}
  If, in Theorem \ref{theo:index_theorem}, $A$ is a finite von Neumann
  algebra and $\tau\colon A\to Z$ is the canonical center valued
  trace, then we can recover $\ind(D_W)$ using the right hand side of
  Equation \eqref{eq:index_formula} or \eqref{eq:oriented_formula},
  since $\tau$ induces an injection $K_0(A)\to Z$ by Proposition
  \ref{prop:Kuenneth} applied to $X=\{*\}$.
\end{remark}

\begin{proof}[Proof of Theorem \ref{theo:index_theorem}]
  By definition, $\ind(D_W)$ and in particular $\ind_\tau(D_W)$ depend
  only on the K-theory class represented by $W$. The same is true for
  $\ch_\tau(W)$ and therefore for the right hand side of Equation
  \eqref{eq:index_formula}.

By Theorem \ref{theo:Kuenneth}, there is an integer $k\in\integers$
such that 
\begin{equation*}
k[W]= \sum_{i=1}^n {\epsilon_i} [E_i\tensor V_i] +
\sum_{j=1}^m\beta^{-1} (F_j\tensor U_j)
\end{equation*}
where $n,m\in\naturals$, $\epsilon_i\in \{-1,1\}$, $E_i$ are finite dimensional complex vector bundles
  on $M$ and $V_i$ are finitely generated projective Hilbert
  $A$-modules. The $F_j$ are elements of $K^1(X)=KK_1(\complexs,C(M))$ and
  the $U_j$ are elements of $K_1(A)=KK_1(\complexs,A)$. $\beta$ is the Bott
  periodicity isomorphism of Theorem \ref{theo:Bott_periodicity}.
 Note that $[E_1\tensor V_1]\in K^0(M;A) =
  KK_0(\complexs,C(M)\tensor A)$ is obtained as the exterior Kasparov
  product of $[E_1]\in KK_0(\complexs,C(M))= K^0(M)$ and $[V_1]\in
  KK_0(\complexs, A)=K_0(A)$.

We now study the summands $\beta^{-1}(F_j\tensor U_j)$ and $\epsilon
[E_i\tensor V_i]$ separately. By definition,
\begin{equation*}
  \begin{split}
    \ind(D_{\beta^{-1}(F_j\tensor U_j)}) &= (\beta^{-1}( F_j\tensor
    U_j))\circ ([D]\tensor \id_A).
  \end{split}
\end{equation*}
Using associativity of the (exterior and interior) Kasparov product,
and the fact that $\beta^{-1}$ is also given by Kasparov product with
a certain element, we get
\begin{equation*}
   \ind(D_{\beta^{-1}(F_j\tensor U_j)})  = \beta^{-1}( F_j\circ
   [D])\tensor U_j) =0,
 \end{equation*}
 since $F_j\circ [D] \in KK_1(\complexs,\complexs)=0$.

In the same way,
  \begin{equation*}
    \ind(D_{E_i\tensor V_i}) = [D]\circ ([E_i]\circ [V_i]) = ([D]\circ
    [E_i])\circ [V_i] = \ind(D_{E_i}) \circ [V_i] = \ind(D_{E_i}) \cdot[V_i].
  \end{equation*}
  Here, for the finite dimensional bundle $[E_i]$, 
  \begin{equation*}
    [D]\circ [E_i] =\ind(D_{E_i}) \in KK_0(\complexs,\complexs) =  \integers,
  \end{equation*}
  i.e.~the Kasparov product gives the Fredholm index of the twisted
  operator.

  Moreover, by the classical Atiyah-Singer index theorem \cite[Theorem
  13.8]{Lawson-Michelsohn(1989)}
  \begin{equation*}
    \ind(D_{E_i}) = \innerprod{ \ch(\sigma(D)) \ch(E_i)\Td(T_\complexs
      M),[TM]},
  \end{equation*}
  therefore $\ind(D_{E_i\tensor V_i}) = \innerprod{ \ch(\sigma(D))
    \ch(E_i),[TM]} [V_i]$, and 
  \begin{equation*}
    \begin{split}
      \ind_\tau(D_{E_i\tensor V_i})& = \innerprod{ \ch(\sigma(D))
        \ch(E_i)\Td(T_\complexs M),[TM]} \tau([V_i])\\
      &= \innerprod{
        \ch(\sigma(D)) \psi_\tau(E_i\tensor V_i)\Td(T_\complexs
        M),[TM]}.
\end{split}
\end{equation*}

Consequently, 
\begin{equation*}
  \begin{split}
  k\ind_\tau(D_W) & = \sum_{i=1}^n\epsilon_i\ind_\tau(D_{E_i\tensor
    V_i}) + \sum_{j=1}^m \ind_\tau(D_{F_j\tensor U_j})\\
  & = \sum_{i=1}^n\epsilon_i\innerprod{\ch(\sigma(D))\psi_\tau(E_i\tensor V_i)
        \Td(T_\complexs M),[TM]}   \;+\,0\\
      & =\innerprod{\ch(\sigma(D))\psi_\tau([ \sum_{i=1}^n\epsilon_i E_i\tensor V_i] )
        \Td(T_\complexs M),[TM]} \\
      & = \innerprod{\ch(\sigma(D))\psi_\tau(k [W])
        \Td(T_\complexs M),[TM]} \\
      &= k \innerprod{\ch(\sigma(D))\ch_\tau(W)
        \Td(T_\complexs M),[TM]},
\end{split}
\end{equation*}
where we also use that $\psi_\tau$ vanishes on the summand
$K^1(X)\tensor K_1(A)\tensor\rationals$ of
$K^0(X;A)\tensor\rationals$. 
The index formula follows. 
\end{proof}

\begin{corollary}
\label{corol:flat_twisting}
  Assume that, in the situation of Theorem \ref{theo:index_theorem}, $W$
  is a flat Hilbert $A$-module bundle with typical fiber $V$. Then
  \begin{equation*}
    \ind_\tau(D_W) = \ind(D) \dim_\tau(W).
  \end{equation*}
\end{corollary}
\begin{proof}
  Combine Theorem \ref{theo:index_theorem} and Corollary
  \ref{corol:ch_of_flat_bundle} and use the classical Atiyah-Singer
  index formula for $\ind(D)$.
\end{proof}

\begin{corollary}
  If $D$ in Theorem \ref{theo:index_theorem} is the spin Dirac
  operator of a spin manifold $M$ of dimension $n=2m$, then
  \begin{equation*}
    \ind_\tau(D_W) = \innerprod{\hat{A}(M)\ch_\tau(W),[M]}.
  \end{equation*}
\end{corollary}
\begin{proof}
  Under this asumption, $\pi_!(\ch(\sigma(D)) \Td(T_\complexs M)) =
  (-1)^m \hat{A}(M)$. Compare the proof of \cite[Theorem 13.10]{Lawson-Michelsohn(1989)}.
\end{proof}

\subsection{Atiyah's $L^2$-index}
\label{sec:atiyahs-l2-index-1}

Now we are in the situation to give a proof of one of the goals of
this paper: Atiyah's $L^2$-index, and
its center valued generalization considered by L{\"u}ck in
\cite{MR1914619} 
can be obtained from the index which an operator defines in the K-theory of
a corresponding $C^*$-algebra.

  Assume that $M$ is a closed manifold, $\Gamma$ a discrete group and
  $M\to B\Gamma$ the classifying map of a $\Gamma$-covering $\tilde M$
  of
  $M$. Consider the corresponding flat bundles $V=\tilde
  M\times_\Gamma C_r^*\Gamma$ and  $H=\tilde M\times_\Gamma
  l^2(\Gamma)$. Let $t=\tau$ be the canonical trace
  $\NeumannN\Gamma\to\complexs$ or the canonical center valued trace
  $\NeumannN\Gamma \to Z$. Let $D$ be a generalized Dirac operator on $M$ with
  lift $\tilde D$ to $\tilde M$. Using $\tilde D$ and $t$, Atiyah
  \cite{Atiyah(1976)} and
  L{\"u}ck \cite{MR1914619}
  define $L^2$-indeces $\ind_{(2)}(\tilde D)\in\complexs$ or
  $\ind_{(2)}(\tilde D)\in Z$, respectively.

\begin{theorem}\label{theo:Atiyah}
  In the situation just described
  \begin{equation*}
    \ind_{(2)}(\tilde D) = t(\ind(D_V)).
  \end{equation*}
\end{theorem}
\begin{proof}
  We have $t([C^*\Gamma])=1$. By Corollary \ref{corol:flat_twisting}
  and the main result of \cite{Atiyah(1976)}
  therefore
  \begin{equation*}
    t(\ind(D_V)) = \ind(D) = \ind_{(2)}(\tilde D).
  \end{equation*}
\end{proof}

\begin{remark}
  The proof of Theorem \ref{theo:Atiyah} we have just given is far from elegant, since we
  compute two indices and then realize that the answers are
  equal. We
  will give an alternative proof in Section \ref{sec:other-poss-proofs}.
\end{remark}

\subsection{Twisted operators}
\label{sec:twisted-operators}

%\section{Operators}

%\begin{definition}
%  Elliptic differential operators on Hilbert $A$-module bundles, in
%  particular twisting ordinary operators with such a bundle.
%\end{definition}

%It is now possible to define directly the index of such an operator,
%essentially in terms of kernel and cokernel, as carried out by
%Mishchenko and Fomenke \cite{???}. However, we will use a description
%of this in terms of KK-theory.

In Definition \ref{def:ind_D_V} we cheated somewhat when defining the
index of $D_W$ without defining the operator $D_W$ itself. However, it
is well known that, at least if $D$ is a generalized Dirac operator,
$D_W$ can be defined as a differential $A$-operator in the sense of
\cite{MIshchenko-Fomenko}. 

We quickly want to review the relevant constructions. Let $D\colon
\Gamma(E)\to\Gamma(E)$ be a generalized Dirac operator on the closed
Riemannian manifold $(M,g)$, acting on the finite
dimensional graded Dirac bundle $E$ with Clifford connection
$\nabla_E$, i.e.~$D$ is the composition
\begin{equation}\label{eq:def_of_twisted_dirac}
  D\colon \Gamma(E_+)\xrightarrow{\nabla_E}\Gamma(T^*M\tensor
  E_+)\xrightarrow{g} \Gamma(TM\tensor E_+) \xrightarrow{c} \Gamma(E_-),
\end{equation}
where $c$ denotes Clifford multiplication.

Assume that $W$ is a finitely generated projective Hilbert $A$-module
bundle with connection $\nabla_W$. Then we define the twisted Dirac
operator $D_W$ in the usual way by
\begin{equation*}
  D_W\colon \Gamma(E_+\tensor W)\xrightarrow{\nabla_E\tensor
    1+1\tensor \nabla_W}\Gamma(T^*M\tensor
  E_+\tensor W)\xrightarrow{g} \Gamma(TM\tensor E_+\tensor W)
  \xrightarrow{c} \Gamma(E_-\tensor W).
\end{equation*}

This is an elliptic differential $A$-operator of order $1$ in the sense of
\cite{MIshchenko-Fomenko} with an index in $K_0(A)$ defined as
follows.

\begin{definition}\label{def:of_Sobolev}
  Given a finitely generated smooth Hilbert $A$-module bundle $E$ over
  a compact manifold $M$, Sobolev spaces $H^s(E)$ can be
  defined ($s\in\reals$), compare e.g.~\cite{MIshchenko-Fomenko}. One
  way to do this is to pick a trivializing
  atlas $(U_\alpha)$ with subordinate partition of unity $(\phi_\alpha)$
  and then define for smooth sections $u,v$ of $E$ the inner
  product
  \begin{equation*}
    (u,v)_s = \sum_\alpha \int_{U_\alpha} \innerprod{(1+\Delta_\alpha)^s
    \phi_\alpha u(x),\phi_\alpha v(x)}\;dx,
  \end{equation*}
  where $\Delta_\alpha$ is the ordinary Laplacian on $\reals^n$ acting
  on the trivialized bundle (some diffeomorphisms of the
  trivializations are omitted to streamline the notation).

  This inner product is
  $A$-valued, and the completion of $\Gamma(E)$ with respect to this
  inner product is $H^s(E)$.
\end{definition}

\begin{remark}
  Of course, the inner product on $H^s(E)$ depends on a number of
  choices, However, two different choices give rise to equivalent
  inner products and therefore isomorphic Sobolev spaces.
\end{remark}

Then $D_W$, being a first order differential operator, induces a bounded
operator $D_W\colon H^s(E_+\tensor W)\to H^{s-1}(E_-\tensor W)$ for
each $s\in\reals$.

The key point is now that the ellipticity of $D$ allows the
construction of a \emph{parametrix} $Q_W$ which induces bounded
operators $Q_W\colon H^{s-1}(E_-\tensor W)\to H^s(E_+\tensor
W)$ for each $s\in\reals$. Parametrix means that 
\begin{equation}\label{eq:parametrix}
D_WQ_W = 1-S_0\qquad Q_W D_W=1-S_1 
\end{equation}
where
$S_0$ and $S_1$ are operators of negative order, i.e.~induce bounded
operators $S_0\colon H^s(E_-\tensor W)\to H^{s+r}(E_-\tensor W)$ and
$S_1\colon H^s(E_+\tensor W)\to H^{s+r}(E_+\tensor W)$ for some $r>0$.

Of course, $S_0$ and $S_1$ in Equation \eqref{eq:parametrix} have to
be interpreted as composition of the above operators with the
inclusion $H^{s+r}\into H^s$.

We now can conclude that $D_W$ indeed gives rise to $A$-Fredholm
operators because of the appropriate version of the Rellich lemma:

\begin{theorem}\label{theo:Rellich}
  If $M$ is compact then the inclusion
\begin{equation*}
  H^{s+r}(E)\to H^s(E)
\end{equation*}
is $A$-compact for each finitely generated projective Hilbert
$A$-module bundle $E$, as long as $r>0$.
\end{theorem}
\begin{proof}
  If $E=M\times V$, $V$ a finitely generated projective Hilbert
  $A$-module, then the definition of $H^s(E)$ amounts to
  \begin{equation*}
    H^s(E) = H^s(M)\tensor V,
  \end{equation*}
  and $i\colon H^{s+r}(E)\into H^s(E)$ becomes $(i\colon H^{s+r}(M)\to
  H^s(M))\tensor \id_V$, i.e.~the tensor product of a compact operator
  (by the classical Rellich lemma) with $\id_V$. Such an operator is
  $A$-compact. The general case follows from an appropriate partition
  of unity argument. A similar argument can be found in \cite[Section
  3]{Schick(2001)}. 
\end{proof}

In particular, $S_0$ and $S_1$ in Equation \eqref{eq:parametrix} are
$A$-compact as composition of the $A$-compact inclusion of the Rellich
Lemma \ref{theo:Rellich} with a bounded operator. Therefore, if we
consider $D_W$ as bounded operator between $H^s$ and $H^{s-1}$ then
$\ind(D_W)\in K_0(A)$ is defined.

\begin{theorem}\label{theo:twist_and_Kasparov product}
  The index just defined is equal to $\ind D_W$ as
  defined in Definition \ref{def:ind_D_V}. In particular, it does not
  depend on $s\in\reals$.
\end{theorem}
\begin{proof}
  This is  a well known fact. For completeness, we want to indicate how this can be
  done. 
  We do this in several steps.

%  \begin{enumerate}
%  \item 
Mishchenko and Fomenko consider the bounded operators
    $D_W\colon  H^s(E_+\tensor  W)\to H^{s-1}(W_-\tensor W)$. These
    are genuine differential operators. We want, however, to relate
    the operators for different $s$ and show that the index is equal
    to the index of the pseudodifferential operator
    $D_W/(\sqrt{1+D_W^2}) \colon L^2(E_+\tensor W)\to L^2(E_-\tensor
    W)$. To do this, we have to observe that $\sqrt{1+D_W^2}$ defines
    bounded even invertible operators $H^s(E_{\pm}\tensor W) \to
    H^{s-1}(E_{\pm}\tensor W)$ which commute with the operator $D$ as
    described above. Note that $D_W/(\sqrt{1+D_W^2})$ usually is
    defined in terms of unbounded normal operators on Hilbert modules,
    as explained in \cite{Lance(1995)}*{Section 9}. Here, we have to
    relate this to the operators between Sobolev spaces. This is not
    quite automatic, since functional calculus for unbounded operators
    on Hilbert $A$-modules is not quite developed in the same way as
    for the case $A=\complexs$. A possible method of proof using
    integral representations (which
    explicitly includes
    some of the results needed here) can be found in
    \cite{MR96e:58148}*{Section 1}.

    Since the index does not change if we compose with a bounded
    invertible operator we conclude two facts:
    \begin{enumerate}
    \item 
      The index of $D_W\colon H^s(E_+\tensor W) \to H^{s-1}(E_-\tensor
      W)$ is the index of 
      \begin{equation*}
D_W\colon H^1(E_+\tensor W)\to
      L^2(E_-\tensor W),
    \end{equation*}
since the first operator is obtained from
      the second by conjugation with the invertible bounded operator $(1+D^2)^{s/2}$. 
    \item The index of the bounded operator $D_W
      (1+D_W^2)^{-1/2}\colon L^2(E_+\tensor W)\to
      L^2(E_-\tensor W)$ is equal to the Mishchenko-Fomenko index,
      since it is obtained from $D_W\colon H^1(E_+\tensor W)\to
      L^2(E_-\tensor W)$ by composition with the invertible bounded
      operator $(1+D^2)^{-1/2}$.
    \end{enumerate}
%  \item

%\label{item:mishfom} 
We have to relate the Mishchenko-Fomenko index to the
    KK-index. Recall from
  \cite{MR99g:46104}*{Section 17.5} that the identification of
  $K_0(A)$ with
  $KK(\complexs,A)$ identifies the index of $D_W\colon H^s(E_+\tensor
  W)\to H^{s-1}(E_-\tensor W)$ with the KK-element represented by the Kasparov tuple
  \begin{equation*}
    \left(  L^2(E_+\tensor W)\oplus L^2(E_-\tensor W) ,
      \begin{pmatrix}
      0 & \frac{D_W}{\sqrt{1+D_W^2}} \\ \frac{D_W}{\sqrt{1+D_W^2}} & 0
      \end{pmatrix}   
\right)
  \end{equation*}
  (note that on $L^2$, $D_W/\sqrt{1+D_W^2}$ is a self adjoint odd
  operator). 
%\item

%\label{item:genKasp} 
We now have to compute the Kasparov product of our first
  definition of the twisted index, and to prove that it equals the
  KK-element just described. Unfortunately, the calculation of the
  Kasparov product is somewhat complicated. We follow here an idea due to
  Ulrich Bunke. Eventually, this comes down to the construction of
  suitable connections in the sense of Kasparov.

  Recall that $\ind(D_W) = [W]\circ ([D]\tensor \id_A) \in
  KK(\complexs, A)$. To analyze the formula, we need explicit
  representatives of the ingredients. Here we have
  \begin{equation*}
    \begin{split}
      [W] =& [\Gamma(W)\oplus 0,0] \in KK(\complexs, C(M;A))\\
      [\id_A]  =& [ A\oplus 0 , 0] \in KK(A,A)\\
      [D] =& [L^2(E^+)\oplus L^2(E^-), 
      \begin{pmatrix}
        0 & D/\sqrt{1+D^2}\\ D/\sqrt{1+D^2} & 0
      \end{pmatrix} ] \in KK(C(M), \complexs);\\
      [D]\tensor [\id_A] =& [ L^2(E^+)\tensor A\oplus L^2(E^-)\tensor
      A, 
      \begin{pmatrix}
        0 & \frac{D}{\sqrt{1+D^2}}\tensor \id_A \\
        \frac{D}{\sqrt{1+D^2}}\tensor \id_A & 0
      \end{pmatrix}]\\
      &\in KK(C(M;A),A).
  \end{split}
\end{equation*}
If $W$ is a graded bundle, a second summand for the negative part
has to be added.

From this, $[W]\circ ([D]\tensor [\id_A]) = (L^2(E^+\tensor W)\oplus
L^2(E^-\tensor W), X)$ with a suitable operator $X$. 

%\item\label{item:calcProd} 
 We claim that $X=      \begin{pmatrix}
        0 & \frac{D_W}{\sqrt{1+D_W^2}}\\
        \frac{D_W}{\sqrt{1+D_W^2}} & 0
      \end{pmatrix}$ is a possible description of this Kasparov
      product. Since $[W]$ is given by a Kasparov tuple with operator
      $0$, it suffices by \cite{MR99g:46104}*{Definition 18.4.1} to
      show that $X$ is a $      \begin{pmatrix}
        0 & \frac{D}{\sqrt{1+D^2}}\tensor \id_A \\
        \frac{D}{\sqrt{1+D^2}}\tensor \id_A & 0
      \end{pmatrix}$-connection for $L^2(E\tensor W)$.
   Since $D$ and $D_W$ both are self adjoint, the connection property
   follows as soon as we show that for each $\gamma\in\Gamma(W)$ the
   operator
   \begin{equation*}
     T_\gamma\circ \frac{D}{(D^2+1)^{-1/2}}\tensor \id_A -
     \frac{D_W}{(D_W^2+1)^{-1/2}}\circ T_\gamma
   \end{equation*}
   is a compact operator from $L^2(E)\tensor A$ to $L^2(E\tensor
   W)$, with $T_\gamma s:= s\tensor\gamma$.

   To do this, we use the integral representation
   \begin{equation*}
\frac{D}{(D^2+1)^{-1/2}} = \int_0^\infty D (D^2+1+\lambda^2)^{-1}
   \;d\lambda,
 \end{equation*}
which by \cite{MR96e:58148}*{Lemma 1.8} is norm
   convergent. By definition of the twisted Dirac operator, for each
   section $s\in L^2(E\tensor A)$
   \begin{multline}\label{eq:Dirac}
    D_W (D_W^2+1+\lambda^2)^{-1} (s\tensor \gamma) \\ =
    (D_W^2+1+\lambda^2)^{-1}(Ds\tensor\gamma) - \sum_i
    (D_W^2+1+\lambda^2)^{-1} X_i\cdot s\tensor \nabla_{X_i}\gamma,
  \end{multline}
  where $\{X_i\}$ is a local orthonormal
  frame and $X_i\cdot s$ denotes Clifford
  multiplication. $(D_W^2+1+\lambda^2)^{-1}\colon L^2(E\tensor W)\to
  L^2(E\tensor W)$ is compact,
  since it factors by \cite{MR96e:58148}*{Lemma 1.5} as a bounded
  operator to $H^2$ composed with the compact inclusion $H^2(E\tensor
  W)\to L^2(E\tensor W)$ (we use here that the base manifold $M$ is compact).

   By \cite{MR96e:58148}*{Lemma 1.5}
   $\norm{(D_W^2+1+\lambda^2)^{-1}}_{\boundedops{(L^2)}} \le
   (d+\lambda^2)^{-1}$ for a suitable constant $d$. For fixed $\gamma\in\Gamma(W)$, the operator
   \begin{equation*}
     s\mapsto \int_0^\infty \left( \sum_i 
   (D_W^2+1+\lambda^2)^{-1} X_i\cdot s\tensor \nabla_{X_i}\gamma\right)\;d\lambda
  \end{equation*}
  therefore is compact as norm convergent integral of compact
  operators.

Consequently, modulo compact operators,
\begin{multline*}
     T_\gamma\circ \frac{D}{(D^2+1)^{-1/2}}\tensor \id_A -
     \frac{D_W}{(D_W^2+1)^{-1/2}}\circ T_\gamma\\
  \equiv \int_0^\infty (T_\gamma\circ (D^2+1+\lambda^2)^{-1} D\tensor\id_A-
  (D_W^2+1+\lambda^2)^{-1} T_\gamma (D\tensor\id_A))\;d\lambda,
\end{multline*}
   using Equation \eqref{eq:Dirac} to commute $D_W$ and
   $T_\gamma$. For each fixed $\lambda$, the integrand is of order
   $-1$ and therefore a compact operator on $L^2(E\tensor W)$ (the
   argument is the same as above).

 Finally,
 \begin{equation*}
   \begin{split}
    \big( T_\gamma &\circ {(D^2+1+\lambda^2)^{-1}}\tensor \id_A -
     {(D_W^2+1+\lambda^2)^{-1}}\circ T_\gamma\big) \circ (D\tensor \id_A)\\
     =& 
     (D_W^2+1+\lambda^2)^{-1} \left( (D_W^2+1+\lambda^2)T_\gamma -
       T_\gamma (D^2+1+\lambda^2)\tensor \id_A \right)\cdot\\
      &\cdot (D^2+1+\lambda^2)^{-1} D\tensor \id_A\\
       =&      (D_W^2+1+\lambda^2)^{-1} \cdot\\
       &\left( T_\gamma D^2\tensor
         \id_A + \sum_{i} T_{\nabla_{X_i}\gamma} X_i\cdot D\tensor
           \id_A +D_W \circ \sum_{i} T_{\nabla_{X_i}\gamma} X_i\cdot\quad -
           T_\gamma D^2\tensor \id_A \right)\\
         &\qquad\cdot
       (D^2+1+\lambda^2)^{-1}D\tensor \id_A
   \end{split}
 \end{equation*}
 For the last step, we use first that $T_\gamma$ commutes with
 $(1+\lambda^2)$, and then we use twice Equation \eqref{eq:Dirac} to
 commute $D_W$ and $T_\gamma$.

 This representation shows that for each fixed $\lambda$ the operator
 in question is actually of order $-2$. Moreover,
  by \cite{MR96e:58148}*{Lemma 1.5} we have
  $\norm{(D^2+1+\lambda^2)^{-1}}_{\boundedops(L^2)} \le
  (d+\lambda^2)^{-1}$. The term in the middle braces is a bounded
  operator from $H^1$ to $L^2$ and is independent of $\lambda$. The
  operator $(D^2+1+\lambda^2)^{-1} D$ is a bounded operator from $L^2$
  to $H^1$ with norm bounded independent of $\lambda$ (since this is
  an operator on a finite dimensional bundle, this is a classical
  fact, it also follows from the definition of the norm on $H^1$ as in
  \cite{MR96e:58148}*{Equation (2)}, where $\abs{s}^2_{H^1} =
  \abs{s}^2_{L^2} + \abs{Ds}^2_{L^2}$, together with the estimates
  $\norm{D^2(D^2+1+\lambda^2)^{-1}}_{\boundedops{(L^2)}} \le C$ and
    $\norm{D(D^2+1+\lambda^2)^{-1}}_{\boundedops{(L^2)}}\le C$ with $C$
    independent of $\lambda$, as given in \cite{MR96e:58148}*{Lemma
      1.5 and Lemma 1.6}.

    It follows that 
    \begin{equation*}
  \int_0^\infty      \big( T_\gamma \circ {(D^2+1+\lambda^2)^{-1}}\tensor \id_A -
     {(D_W^2+1+\lambda^2)^{-1}}\circ T_\gamma\big) \circ (D\tensor
     \id_A) \;d\lambda
    \end{equation*}
  converges in operator norm on $L^2(E\tensor W)$. Since the integrand
  consists  of compact operators and the ideal of compact operator
  is norm closed, it follows as above that the whole integral is
  compact.

  Modulo compact operators, this is equal to $T_\gamma\circ
  \frac{D}{(D^2+1)^{-1/2}}\tensor \id_A - \frac{D_W}{(D_W^2+1)^{-1/2}}
  T_\gamma$, which is therefore compact as we had to show.

%\item
  Putting the above arguments together, it follows 
%Because of \ref{item:mishfom}, the calculation in \ref{item:genKasp}
 % and \ref{item:calcProd} show 
that the Mishchenko-Fomenko index
  equals the Kasparov product, as claimed. This finishes the proof of
  the Theorem.
% \end{enumerate}

\end{proof}

In \cite{MIshchenko-Fomenko}, a ``cohomological'' formula for this
index is
derived similar to our formula \ref{theo:index_theorem}. The
underlying strategy uses similar ideas, namely the K{\"u}nneth theorem
\ref{theo:Kuenneth} to reduce to the classical Atiyah-Singer index theorem.
The original index theorem is less explicit, because it does not take
the curvature of the twisting bundle into account. In particular,
Corollary \ref{corol:flat_twisting} does not follow directly. On the other
hand, it is more precise because it gives
K-theoretic information, whereas we neglect the part of K-theory which
is not detectable by traces. Note that, if $A$ is a finite von Neumann
algebra, by Proposition \ref{prop:Kuenneth} no information is lost.

\section[Simplified index for von Neumann algebras]{A simplified
  $A$-index for von Neumann
  algebras}
\label{sec:simpl-gener-atkins}

In this section, $A$ is assumed to be a von Neumann algebra.

Let $H_A$ be the  Hilbert $A$-module which is the
completion of $\oplus_{i=1}^\infty A$. Then $\End_A(H_A)\iso
\boundedops(H)\tensor A$, where $H$ is a separable Hilbert
space. The ``compact'' operators $K_A(H_A)$ in $\End_A(H_A)$, i.e.~the
$C^*$-algebra generated by the operators of the form $x\mapsto
v\innerprod{w,x}$ for some $v,w\in H_A$ are isomorphic to
$K(H)\tensor A$. 

One can now define the $A$-Fredholm operators $F_A(H_A)$ in $\End_A(H_A)$ to be those
operators which are invertible module $K_A(H_A)$. The generalized
Atkinson theorem states that  a suitably defined index induces an
isomorphism between the set of path components of $F_A(H_A)$ (a group
under composition) and $K_0(A)$, compare \cite[Chapter
17]{Wegge-Olsen}, originally proved by Kasimov in \cite{MR688514}. We
refer to the textbook \cite{Wegge-Olsen} because of its easy availability
and because it is rather self contained.

The problem with the definition of the index is that kernel and cokernel
of a Fredholm operator as defined above are not necessarily finitely
generated projective $A$-modules. The way around this is to compactly
perturb a given Fredholm operator.

We want to show here that this is not necessary if $A$ is a von
Neumann algebra.

The main virtue of the following result is that in case $A$ is a von
Neumann algebra, the index of an $A$-Fredholm operator is determined using
spectral calculus instead of some compact perturbation which can
hardly be made explicit.
\begin{theorem}\label{theo:index_via_ker_projection}
  Assume that $A$ is a von Neumann algebra and  $f\in \End_A(H_A)$ is
  an $A$-Fredholm operator. Since
  $\End_A(H_A)$ is a von Neumann algebra, we can use the measurable
  functional calculus and define the projections
  $p_{\ker}:=\chi_{\{0\}}(f^*f)$ and $p_{\coker}:=\chi_{\{0\}}(ff^*)$,
  where $\chi_{\{0\}}$ is the characteristic function of the set
  $\{0\}$. Then  $\im(p_{\ker})$ and $\im(p_{\coker})$ are finitely
  generated projective Hilbert $A$-modules and
  $[\im(p_{\ker})]-[\im(p_{\coker})] = \ind_A(f)\in 
  K_0(A)$, with $\ind_A:=\textrm{Mindex}$ defined in \cite[Chapter
  17]{Wegge-Olsen} as $[\ker(f+k)]-[\coker(f+k)]$ for a suitable
  $A$-compact perturbation of $f$ (any $k$ such that range, kernel and
  cokernel of $f+k$ are closed will do).
\end{theorem}
\begin{proof}
  Since $f$ is invertible module $A$-compact operators and
  $fp_{\ker}=0$, $p_{\ker}$ is zero module compact operators, i.e.~a
  compact projection. The same is true for $p_{\coker}$. By
  \cite[Theorem 16.4.2]{Wegge-Olsen}, their images are finitely
  generated projective Hilbert $A$-modules, so that in particular
  $[\im(p_{\ker})] - [\im(p_{\coker})] \in K_0(A)$ is defined.

  Since $\End_A(H_A)$ is a von Neumann algebra, each operator has a
  polar decomposition (for general $A$, this is only assured for those
  with closed range, compare \cite[Theorem 15.3.8]{Wegge-Olsen}.) Write
  therefore $f= u\abs{f}$ with a partial isometry $u$. By spectral
  calculus, $1-u^*u= p_{\ker}$ and $1-uu^*= p_{\coker}$. If $g=f+k$ is an
  $A$-compact perturbation of $f$, and $g=v\abs{g}$ is its polar
  decomposition, then $u-v$ is $A$-compact, as follows from the proof
  of \cite[Corollary 17.2.5]{Wegge-Olsen} and therefore by
  \cite[Corollary 17.2.4]{Wegge-Olsen}
  \begin{multline*}
    [p_{\ker}]-[p_{\coker}] = [1-u^*u]-[1-uu^*] \\ = [1-v^*v]-[1-vv^*] =
    [\ker(g)]-[\ker(g^*)] \in K_0(A). 
  \end{multline*}
  Since the latter is by definition the $A$-index of $f$, we are done.
\end{proof}

\begin{remark}
  Occasionally, we will use the notation $[p_{\ker(f)}]\in K_0(A)$ for
  the K-theory element represented by the image of $\ker(f)$, if we
  are in the situation of Theorem
  \ref{theo:index_via_ker_projection}. Note that we have to enlarge
  the standard ``finite projective matrices'' description a little bit
  here, since the projection is only unitarily equivalent (with a
  unitary close to one) to a finite
  projective matrix, as is proved e.g.~in \cite[Lemma
  15,4.1]{Wegge-Olsen}. We have to keep in mind that not all
  constructions immediately generalize to these generalized
  projections, e.g.~when applying traces to them.
\end{remark}

\begin{definition}\label{def:general_index}
  Let $V$ and $W$ be (topologically) countably generated Hilbert
  $A$-modules and $f\in \Hom_A(V,W)$. We call $f$ \emph{Fredholm} if
  $f\oplus \id_{H_A}\colon V\oplus H_A\to W\oplus H_A$ is Fredholm. If
  this is the case, then
  \begin{equation*}
    \ind_A(f):=\ind_A(f\oplus\id_{H_A}) \in K_0(A).
  \end{equation*}
  Observe that this definition makes sense and reduces to the
  situation of Theorem \ref{theo:index_via_ker_projection} since by
  Kasparov's stabilization theorem \cite[Theorem 15.4.6]{Wegge-Olsen}
  $V\oplus H_A\iso H_A$.
\end{definition}

\begin{corollary}\label{corol:general_ind_via_ker_proj}
  If $A$ is a von Neumann algebra, $V$ and $W$ are countably generated
  Hilbert $A$-modules and $f\in\Hom_A(V,W)$ is Fredholm, then
  \begin{equation*}
    \ind_A(f) = [\chi_{\{0\}}(f^*f)] - [\chi_{\{0\}}(ff^*)] \in K_0(A).
  \end{equation*}
\end{corollary}
\begin{proof}
  This is an immediate consequence of Definition
  \ref{def:general_index} and of Theorem
  \ref{theo:index_via_ker_projection}. 
\end{proof}

We can apply this to the twisted generalized Dirac operators
considered in Section \ref{sec:twisted-operators}
\begin{corollary}\label{corol:ind_D_V_with_ker_projection}
  Let $D\colon \Gamma(E_+)\to \Gamma(E_-)$ be a generalized Dirac
  operator, acting on the sections of a finite dimensional bundle
  $E$ over the smooth compact manifold $M$ without boundary. Let $A$
  be a von Neumann algebra and $W$ a smooth finitely generated
  projective Hilbert $A$-module bundle. Then the $A$-index of the
  twisted operator $D_W$ as defined in Definition \ref{def:ind_D_V} or
  Subsection \ref{sec:twisted-operators} can be expressed as follows:
  \begin{equation*}
    \ind_A(D_W) =[\chi_{\{0\}}(D_W^*D_W)] - [\chi_{\{0\}}(D_WD_W^*)] \in K_0(A),
  \end{equation*}
  where we understand $D_W$ to be the bounded operator
  \begin{equation*}
    D_W\colon H^1(E_+\tensor W)\to H^0(E_-\tensor D_W).
  \end{equation*}
\end{corollary}

\section{A general Atiyah $L^2$-index theorem}
\label{sec:atiyahs-l2-index}

\subsection{$A$-Hilbert spaces and bundles}
\label{sec:a-hilbert-spaces}

Atiyah's $L^2$-index theorem \cite{Atiyah(1976)} and its
generalization by  L{\"u}ck \cite{MR1914619} deal with indices obtained from
an ordinary elliptic differential operator and a trace on a von
Neumann algebra $A$, but this is done in a different way compared to
the construction in Definition \ref{def:ind_D_V}. 

Atiyah is looking at coverings of a compact manifold and a lifted Dirac type
operator (this corresponds to the twist with the canonical flat bundle of the
covering of Example \ref{ex:Hilbert_Gamma}), and is proving that the
$L^2$-index (associated to a canonical trace) coincides with the ordinary index of the operator on the compact
base manifold. He is using a parametrix construction to directly show
that the two numbers coincide. L{\"u}ck, in the same situation, is studying all
the other normal traces. He proves that they don't contain additional
information. L{\"u}ck is using the heat kernel on the covering manifold. A
proof of Atiyah's original result
using heat kernel methods is given in \cite{roe}. L{\"u}ck is also giving a
K-theoretic interpretation of his result: the index in question defines an
element of $K_0(\NeumannN\Gamma)$ which is a multiple of the trivial element $1$. This is an
infinite dimensional generalization of the well known rigidity theorem
which says that for a free action of a finite group, the equivariant
index contains no more information than the ordinary index (compare
\cite[Remark after Theorem 0.4]{MR1914619}).

Despite the different definitions and methods, there
is an easy direct translation between the two aspects, which is well
known and frequently used in the literature, but seems not to be
documented with proof. Therefore, our goal here is to prove
this connection. This is inspired
by a remark of Alain Valette who missed a citable reference for the result.

In the present subsection, we will introduce the notation and concepts
necessary to give the definition of Atiyah's (and L{\"u}ck's)
$L^2$-index. We do this in a more general setting, making transparent
some of the connections to the previous parts of this paper.

We have to introduce some further notation. Unfortunately, the term
``(finitely generated projective) Hilbert $A$-module'' is used in the
literature for two different things: the objects we have introduced so
far, but also the objects on which Atiyah's definition of the
$L^2$-index is based. The latter are honest Hilbert spaces with an
action of the  $C^*$-algebra $A$. To distinguish them from the
objects introduced above, we use the term ``$A$-Hilbert space''
(deviating from the literature at this point). We will see in Section
\ref{sec:an-equiv-categ} how to translate between these two concepts.

For our construction, we use a trace on $A$ with particular
properties. This will exist in our main example, the von Neumann
algebra of a discrete group. For the following, we recall the
construction of $l^2(A)$ which is used to pass from the algebra $A$ to
an $A$-Hilbert space.

\begin{definition}\label{def:good_L2_trace}
  Let $A$ be a $C^*$-algebra and $Z$ a commutative $C^*$-algebra
  (most important is the example $Z=\complexs$).
  A \emph{trace} $\tau\colon A\to Z$ is a linear map such that
  \begin{enumerate}
  \item $\tau(ab)=\tau(ba)$ for each $a,b\in A$.
  \item It is called \emph{positive} if $\tau(a^*a)\ge 0$ for each $a\in
    A$.
  \item It is called \emph{faithful} if $\tau(a^*a)=0$ only for $a=0$.
  \item It is called \emph{normalized} if $\tau(1)=1$. 
  \item If $A$ and $Z$ are von Neumann algebras, a positive trace $\tau$ is called \emph{normal} if it is
    ultraweakly continuous.
  \end{enumerate}
\end{definition}

\begin{notation}\label{not:special_tau}
  From now on, we assume the existence and fix a positive faithful
  normalized trace $\tau\colon A\to \complexs$.
\end{notation}

\begin{lemma}\label{lem:continuity}
  Given a trace $\tau$ as in \ref{not:special_tau}, we have the following inequality:
  \begin{equation*}
    \tau(a^*xa) \le \abs{x} \tau(a^*a)\qquad\text{if $x\in A$ is
      positive,} a\in A,
  \end{equation*}
  with $\abs{x}$ the $\reals$-valued norm of $x\in A$.

  In particular, with $a=1$, the map $\tau\colon A\to\complexs$ is
  norm continuous.
\end{lemma}
\begin{proof}
  In $A$, we have $x\le \abs{x}$ and therefore $a^*xa\le a^*\abs{x} a
  =\abs{x}a^*a$. Positivity and linearity of the trace impies the
  inequality. 
\end{proof}

\begin{definition}\label{def:pre_innerprod}
  Given the positive faithful normalized trace $\tau$ on the
  $C^*$-algebra $A$ as in \ref{not:special_tau}, define a
  sesqui-linear inner
  product on a Hilbert $A$-module $V$ by
  $\innerprod{v,w}_{2}=\tau(\innerprod{v,w})$ (linear in the second 
  entry), i.e.~we compose the $A$-valued inner product with $\tau$.
\end{definition}

\begin{lemma}
\label{lem:A_acts_boundedly_on_l2A} 
  In the situation of Definition
\ref{def:pre_innerprod}, $V$ with the constructed inner product
becomes a pre Hilbert space. Its
  completion is denoted $l^2(V)$. Right multiplication
  of $A$ on $V$ induces a $C^*$-homomorphism from $A$ to the bounded
  operators on $l^2(V)$.

  In the special case $V=A^n$, left and right multiplication both
  induce $C^*$-embeddings of $A$ into the bounded operators on
  $l^2(A)^n=l^2(A^n)$.
\end{lemma}
\begin{proof}
  Since $\tau$ is faithful and positive and the same is true for
  $\innerprod{\cdot,\cdot}$, $\innerprod{\cdot,\cdot}_{2}$
   induces a norm $\norm{\cdot}$. If $a,x\in A$, $v\in V$ then by Lemma
  \ref{lem:continuity} 
  \begin{equation*}
    \begin{split}
      \norm{va} =
      \tau(\innerprod{va,va})^{1/2} & =\tau(a^*\innerprod{v,v}a)^{1/2} =
      \tau(\sqrt{\innerprod{v,v}} a a^*\sqrt{\innerprod{v,v}})^{1/2} \\
      \le \abs{a^*a}^{1/2} \tau(\innerprod{v,v})^{1/2} = \abs{a} \cdot
      \norm{v}
  \end{split}
\end{equation*}
  For left multiplication of $A$ on $A$
  \begin{equation*}
    \norm{ax} = \tau(x^*a^*ax)^{1/2} \le \abs{a} \norm{x}.
  \end{equation*}
  We conclude that right multiplication and for $V=A^n$ also left  multiplication by
  $a$ give rise to
  bounded operators with operator norm $\le \abs{a}$. The
  corresponding maps are $*$-homomorphisms since
  \begin{equation*}
    \innerprod{va,w}_2=\tau(a^*\innerprod{v,w})=
    \tau(\innerprod{v,w}a^*)=\innerprod{v,wa^*}_2\qquad\forall a\in A,
    \;v,w\in V. 
  \end{equation*}
  \begin{equation*}
        \innerprod{ax,y}_2=\tau(x^*a^*y)=\innerprod{x,a^*y}_2; \qquad
        \forall a,x,y\in A
  \end{equation*}
  Left or right multiplication by $a$ on $A^n$ is the zero map only if
  $a=0$.
\end{proof}

\begin{remark}
  \begin{enumerate}
  \item In Lemma \ref{lem:A_acts_boundedly_on_l2A}, $l^2(A)$ and
    $l^2(V)$ depend of course on the chosen trace $\tau$. We will not
    indicate this in the notation since we adopt the convention that
    the trace $\tau$ is fixed throughout. Moreover, we will see in
    Section \ref{sec:an-equiv-categ} that one can recover $V$ from
    $l^2(V)$, such that the particular choice of $\tau$ does not play
    too much of a role.
  \item Lemma \ref{lem:A_acts_boundedly_on_l2A} contains the easy case of
    the representation theorem for $C^*$-algebras: if $A$ has a trace
    as in Definition \ref{def:good_L2_trace} then $A$ can be
    isometrically embedded into the algebra of bounded operators on
    the Hilbert space $l^2(A)$.
\end{enumerate}
\end{remark}

\begin{definition}\label{def:of_A_Hilbert_space}
  A finitely generated projective $A$-Hilbert space $V$ is a Hilbert
  space together with a right action of $A$ such that $V$ embeds
  isometrically preserving the $A$-module structure as a direct
  summand into $l^2(A)^n$ for some $n$, and such that the orthogonal
  projection $l^2(A)^n\onto V$ is given by left multiplication with an
  element of $M_n(A)$. 

  A (general) \emph{$A$-Hilbert space} $V$ satisfies the same
  conditions a finitely generated projective $A$-Hilbert space does,
  with the exception that $l^2(A)^n$ is replaced by $H\tensor l^2(A)$
  for some Hilbert space $H$ (the tensor product has to be
  completed), and $M_n(A)$ by $\boundedops(H)\tensor A$ (where $A$ is
  here understood to act by {left} multiplication). Observe that, if $H$ is
  separable, then $H\tensor l^2(A)\iso l^2(H_A)$, and
  $\boundedops(H)\tensor A\iso \Hom_A(H_A)$.
\end{definition}

\begin{remark}
  Assume that, in Definition \ref{def:of_A_Hilbert_space}, $A$ is a von
  Neumann algebra. Then the condition that the projection $H\tensor
  l^2(A)\onto V$ belongs to $\boundedops(H)\tensor A$ is automatically
  satisfied, since the commutant of the right multiplication of $A$ on
  $H\tensor l^2(A)$ is $\boundedops(H)\tensor A$ (and on $l^2(A^n)$ is
  $M_n(A)$), and the projection
  by definition commutes with the right multiplication of $A$.
\end{remark}

% \Kommentar{ This example is done a little later.}
%\begin{example}
%  Assume that $\Gamma$ is a countable group. Consider the Hilbert
%  space $l^2(\Gamma)$ of $L^2$-functions (counting measure) on
%  $\Gamma$. The algebra $\complexs\Gamma$ acts on $l^2(\Gamma)$ from
%  the left and from the right by convolution. The group von Neumann
%  algebra $\NeumannN\Gamma$ is defined equivalently either as the weak
%  closure of $\complexs\Gamma$, considered as a subset of
%  $\boundedops(l^2(\Gamma))$ by right convolution, or as the commutant
%  $\boundedops(\l^2(\Gamma))^\Gamma$ of the left convolution action
%  of $\complexs\Gamma)$ (the commutant is by definition the set of
%  operators which commutes with the designated subset). Then we define
%  the \emph{canonical trace}
%  \begin{equation*}
%    \tau\colon \NeumannN\Gamma\to\complexs;\; a\mapsto \innerprod{a(1),1}_{\l^2(\Gamma)},
%  \end{equation*}
%  where $1$ is the function on $\Gamma$ with value $1$ at the identity
%  element, and $0$ elsewhere. On the subalgebra $\complexs\Gamma$ of
%  $\NeumannN\Gamma$, $\tau$ picks the coefficient of the identity
%  element.

%  Then $l^2(\NeumannN\Gamma)$ constructed with this trace $\tau$ is
%  canonically isomorphic to $l^2(\Gamma)$. The isomorphism is induced
%  by the linear map
%  \begin{equation*}
%    \NeumannN\Gamma\to l^2(\Gamma); \; a\mapsto a(1).
%  \end{equation*}
%\end{example}

\subsection{$A$-Hilbert space bundles}
\label{sec:hilb-space-bundl}

\begin{definition}
 An $A$-Hilbert space morphism is a bounded $A$-linear map
    between two $A$-Hilbert spaces. If it is an isometry for the
    Hilbert space structure, it is called an $A$-Hilbert space isometry.

 An $A$-Hilbert space bundle $H$ on a space $X$ is a locally
    trivial bundle of $A$-Hilbert spaces, the transition functions
    being $A$-Hilbert space isometries. A smooth structure is given by
    a trivializing atlas where all the transition functions are
    smooth.

    If the fibers are finitely generated projective $A$-Hilbert space,
    the bundle is called a \emph{finitely generated projective
      $A$-Hilbert space bundle}.
\end{definition}

\begin{lemma}\label{lem:sections are A Hilbert space}
  The $L^2$-sections of an $A$-Hilbert space bundle $W$ on a
  Riemannian manifold $X$ form themselves an $A$-Hilbert space.
\end{lemma}
\begin{proof}
  The action of $A$ is given by pointwise multiplication. We want to show that $L^2(W)\iso
  L^2(M)\tensor V$, where $V$ is a typical fiber of $M$ (we assume for
  simplicity that $M$ is connected). Since $V$ embeds into $H\tensor
  l^2(A)$, the same is then true of $L^2(W)$.

  To prove that $L^2(W)\iso L^2(M)\tensor V$, choose a subset
  $U\subset M$ such that $M\setminus U$ has measure zero, and such
  that $W|_U$ is trivial ($U$ could e.g.~consist of the interiors of
  the top cells of a smooth triangulation of $M$). Then $L^2(W) \iso
  L^2(W|_U)\iso L^2(U)\tensor V\iso L^2(M)\tensor V$, since $U$ and
  $M$ differ only by a set of measure zero, and since $W|_U\iso
  U\times V$.
\end{proof}

As an example, we now want to give the most important $A$-modules,
$A$-Hilbert spaces and corresponding bundles. To do this, we have in
particular to specify the von Neumann algebra $A$. This is the
$A$-Hilbert space bundle featuring in Atiyah's $L^2$-index theorem and
its generalization by L{\"u}ck.

\begin{example}\label{ex:Hilbert_Gamma}
  Let $M$ be a smooth compact manifold and $\Gamma$
  its fundamental group. Let $\pi\colon \tilde M\to M$ be a universal
  covering of $M$, with $\Gamma$-action from the right by deck
  transformations.

  The Hilbert space $l^2(\Gamma)$ is the space of complex valued
  square summable functions on the discrete group
  $\Gamma$. $\complexs\Gamma$ acts through bounded operators on
  $l^2(\Gamma)$ by left as well as right convolution
  multiplication. By definition, the reduced $C^*$-algebra
  $C^*_r\Gamma$ of $\Gamma$ is the norm closure in
  $\boundedops(l^2(\Gamma))$ of $\complexs\Gamma$
  acting from the right, and $\NeumannN\Gamma$ is the weak closure of
  the same algebra. By the double commutant theorem, this is the set
  of all operators which commute with left convolution of
  $\complexs\Gamma$.

  On $\NeumannN\Gamma$ and therefore also on its subalgebra
  $C^*_r\Gamma$ we have the canonical faithful positive trace $\tau$
  with
  $\tau(f)=\innerprod{f(1),1}_{l^2\Gamma}$, where $1\in l^2(\Gamma)$ is
  by definition the characteristic function of the unit element.

  The construction of $l^2(C^*_r\Gamma)$
  and of $l^2(\NeumannN\Gamma)$ with respect to this trace yields
  precisely $l^2(\Gamma)$.

  Since the left $\Gamma$-action and the right
  $C^*_r\Gamma$ or $\NeumannN\Gamma$-action, respectively, on
  $l^2(\Gamma)$ and $C^*_r\Gamma$ or $\NeumannN\Gamma$, respectively,
  commute, the bundles 
  $\tilde M\times_\Gamma C^*_r\Gamma$ and  $\tilde M\times_\Gamma
  \NeumannN\Gamma$ are smooth finitely generated projective
  Hilbert $C^*_r\Gamma$ and Hilbert $\NeumannN\Gamma$ module bundle,
  and $\tilde M\times_\Gamma l^2(\Gamma)$ is a finitely generated
  projective $C^*_r\Gamma$-Hilbert space or $\NeumannN\Gamma$-Hilbert
  space bundle, all on $M$. Moreover, $\tilde M\times_\Gamma
  l^2(\Gamma)$ can be considered as the $A$-Hilbert space completion
  of the former bundles with respect to the canonical trace.

  To see that the bundles are smooth, observe that the canonical
  trivializations are obtained by choosing lifts to $\tilde M$, and
  the transition functions are then given by left multiplication with
  fixed elements $\gamma\in \Gamma$. Since these maps do not depend on
  the basepoint in $M$ they are smooth (the argument shows that these
  bundles are actually flat).

  The same construction works if $\Gamma$ is some homomorphic image of
  the fundamental group of $M$, and $\tilde M$ the corresponding
  normal covering space of $M$.

  The trivial connection on $\tilde M\times C^*_r\Gamma$ and $\tilde
  M\times \NeumannN\Gamma$ descents to a canonical flat connection on
  $\tilde M\times_\Gamma
  C^*_r\Gamma$ and $\tilde M\times_\Gamma \NeumannN \Gamma$, since left
  (as well as right) multiplication with an element $\gamma\in\Gamma$
  is parallel.
\end{example}

\subsection{Connections on $A$-Hilbert space bundles}
\label{sec:conn-a-hilb}

\begin{definition}
  Let $A$ be a von Neumann algebra with a trace $\tau$ as in
  \ref{not:special_tau}. 
  Assume that $M$ is a smooth manifold and $X$ is a smooth finitely generated
  projective $A$-Hilbert space bundle on $M$. 
A \emph{connection} $\nabla$ on $X$ is an $A$-linear map $\nabla\colon
  \Gamma(X)\to
  \Gamma(T^*M\tensor X)$ which is a derivation with respect to
  multiplication with
  sections of the trivial bundle $M\times A$, i.e.
  \begin{equation*}
    \nabla(sf ) = s df  +  \nabla(s)f\qquad\forall
    s\in\Gamma(X),\; f\in C^\infty(M;A).
  \end{equation*}
  Here we use the multiplication $X\tensor T^*M\tensor (M\times A)\to
  X\tensor T^*M\colon s\tensor \eta\tensor f\mapsto sf\tensor
  \eta$. (In particular, elements of $A$ are considered to be of
  degree zero.)

  We say that $\nabla$ is a \emph{metric connection} if
  \begin{equation*}
    d\innerprod{s_1,s_2} = \innerprod{\nabla s_1,s_2} +
    \innerprod{s_1,\nabla s_2}
  \end{equation*}
  for all smooth sections $s_1,s_2$ of $X$. Here, we consider
  $\innerprod{s_1,s_2}$ to be a section of the trivial bundle $M\times
  \complexs$.
\end{definition}

\begin{example}
  In the situation of Example \ref{ex:Hilbert_Gamma}, $\tilde
  M\times_\Gamma l^2(\Gamma)$ inherits a canonical flat connection,
  descending from $\tilde M\times l^2(\Gamma)$, which extends the
  corresponding flat connection on the subbundle
  $\tilde M\times_\Gamma \NeumannN\Gamma$.
\end{example}

\subsection{Operators twisted by $A$-Hilbert space bundles}
\label{sec:twist-oper-with}

In this paper, we will only twist ordinary Dirac type differential operators with
$A$-Hilbert space bundles. For a more complete theory of
(pseudo)differential operators on such bundles compare
e.g.~\cite[Section 2]{BFKM(1996)}.

\begin{definition}\label{def:of_H_Hilbert_twisted_Dirac}
  Let $D\colon\Gamma(E^+)\to\Gamma(E^-)$ be a generalized Dirac
  operator between sections of finite dimensional bundles on the
  Riemannian manifold $(M,g)$.

  Let $H$ be a smooth $A$-Hilbert space bundle with connection
  $\nabla_H$. Then we define (as usual) the twisted Dirac operator
  \begin{equation*}
    D_H\colon \Gamma(E_+\tensor H) \xrightarrow{\nabla\tensor 1
      +1\tensor\nabla_H} \Gamma(T^*M\tensor E_+\tensor H)
    \xrightarrow{g} \Gamma(TM\tensor E_+\tensor H)\xrightarrow{c}
    \Gamma(E_-\tensor H),
  \end{equation*}
  where $c$ stands for Clifford multiplication.

  This is an elliptic differential operator of order $1$ on
  $A$-Hilbert space bundles in the
  sense of \cite{BFKM(1996)}. In particular, it extends to an
  unbounded operator on $L^2(E\tensor H)$. 

  If $A$ is a von Neumann algebra, then the kernel as well as the
  orthogonal complement of the image are
  $A$-Hilbert spaces. The $A$-action is evident. The assertion about
  the projections follows from the fact that by measurable functional
  calculus, the projection onto the kernel of $A$ is given by
  $\chi_{\{0\}}(D_H^* D_H)$ ($\chi_{\{0\}}$ being the characteristic
  function of $\{0\}$), and similarly for the cokernel.
\end{definition}

\begin{remark}
  If $A$ is not a von Neumann algebra, kernel and cokernel are not
  necessarily $A$-Hilbert modules.
\end{remark}

\begin{definition}\label{def:dim_t}
  Assume that $A$ is a von Neumann algebra with a trace $\tau$ as in \ref{not:special_tau}.
  Let $t\colon A\to Z$ be a second trace which is required to be
  positive and normal (but not necessarily faithful or normalized), with values in a
  commutative von Neumann algebra $Z$ ($t=\tau$ is permitted). Given an
  $A$-Hilbert module $V$, we 
  define
  \begin{equation*}
    \dim_t(V):= t(\pr_V),
  \end{equation*}
  where $\pr_V\colon l^2(A)\tensor H\to l^2(A)\tensor H$ is the orthogonal
  projection onto $V$, and $t$ here also stands for the extension of
  the trace to $A\tensor\boundedops(H)$ (to do this, the fact that the
  trace $t$ is normal has to be used). We will discuss the definition
  and properties of these traces in Section \ref{sec:properties-traces}.
\end{definition}

\begin{definition}\label{def:ind_t}
  Let $A$ be a von Neumann algebra with traces $t$ and $\tau$ as in
  Definition \ref{def:dim_t}.

  Let $D_H$ be a generalized Dirac operator twisted by a finitely
  generated projective $A$-Hilbert space bundle $H$ as in Definition
  \ref{def:of_H_Hilbert_twisted_Dirac}. Assume that $M$ is compact
  without boundary. Ellipticity implies that
  $\chi_{\{0\}}(D_H^*D_H)$ and $\chi_{\{0\}}(D_HD_H^*)$ are of
  $t$-trace class (compare Section \ref{sec:properties-traces} for the
  definition and Section \ref{sec:other-poss-proofs} for a proof of
  this fact). Then define
  \begin{equation*}
    \ind_t(D_H) := t(\chi_{\{0\}}(D_H^*D_H)) - t(\chi_{\{0\}}(D_HD_H^*)).
  \end{equation*}
\end{definition}

Our goal now is to prove an index formula for $\ind_t(D_H)$ in the
general situation of Definition
\ref{def:of_H_Hilbert_twisted_Dirac}. One way to do this would be the
following: 
\begin{enumerate}
\item develop a theory of connections and curvature for $A$-Hilbert space
  bundles similar to what we have done for Hilbert $A$-module bundles.
  This is possible in exactly the same way as done above.
\item\label{item:lower_order_changes} Show that $\ind_t$ is unchanged by lower order perturbations of
  $D_H$ (in particular if the connection on $H$ is changed). One way
  to do this would be to prove that $\ind_t$ can be calculated from
  the remainder terms $S_0$ and $S_1$ in $D_HQ=1-S_0$ and
  $QD_H=1-S_1$, where $Q$ is a suitable parametrix (such that the
  remainder terms are of $t$-trace class), namely
  \begin{equation*}
    \ind_t(D_H) = t(S_1)-t(S_0).
  \end{equation*}
This step is
  already done by Atiyah \cite{Atiyah(1976)} (in his special
  situation), and his proof
  does only use a few general properties of the trace, in particular
  that it is normal, a trace, and that operators of order $-k$, for
  $k$ sufficiently big, are of trace class. Since all these properties
  are satisfied here, the proof goes through. A more formal discussion
  of this prove can be found in \cite{Schick(2001)}. For a lower order
  perturbation $D_H-a$ of $D_H$, we can then use the parametrix
  $Q'=Q+QaQ+QaQaQ+\cdots +QaQ\cdots aQ$. Then
  $(D_H-a)Q'=1-S_0-aQ\cdots aQ$, and $Q'(A_H-a)=1-S_1-Qa\cdots Qa$,
  and the trace property implies immediately that
  \begin{equation*}
    t(S_1')-t(S_0') = t(S_1)-t(S_0).
  \end{equation*}
\item Follow the proof of Theorem \ref{theo:index_theorem} to get a
  very similar formula for $\ind_t$.
\end{enumerate}

Although all this can be done, Step \ref{item:lower_order_changes} is
rather lengthy. Therefore, we prefer to show in Section
\ref{sec:other-poss-proofs} that the ``new''
situation can be reduced to the index theorem
\ref{theo:index_theorem} by directly showing that
\begin{equation}\label{eq:A Hilbert reduces to Hilbert A}
  \ind_t(D_H) = t(\ind(D_V))
\end{equation}
for a finitely generated projective Hilbert $A$-module bundle $V$
canonically associated to $H$ (in particular, $\ind(D_V)\in K_0(A)$). 

%\Kommentar{This can probably be removed: As a
%special case Theorem \ref{theo:Atiyah} will directly follow:
%\begin{equation*}
%  \ind_t(D_{\tilde M\times_\Gamma l^2(\Gamma)}) = t(\ind(D_{\tilde
%    M\times_\Gamma \NeumannN\Gamma})),
%\end{equation*}
%without using Atiyah's calculation of $\ind_{(2)}(\tilde D)$.}

\subsection{Flat $A$-Hilbert space bundles and coverings}
\label{sec:flat-a-hilbert}

  Assume that $A=\NeumannN\Gamma$ is the von Neumann algebra of the discrete
  group $\Gamma$ and $t=\tau$ is the canonical trace of Example
  \ref{ex:Hilbert_Gamma}. Let $H=\tilde M\times_\Gamma l^2(\Gamma)$ be
  the canonical flat $l^2(\Gamma)$-bundle of Example
  \ref{ex:Hilbert_Gamma}, and let $D\colon\Gamma(E^+)\to\Gamma(E^-)$ be a
  generalized Dirac operator
  on $M$. In this situation, we have defined
  $\ind_t(D_H)\in\reals$. Fix, more generally, an element $g\in\Gamma$
  which has only finitely many conjugates, and let $[g]$ be this
  finite conjugacy class. Then it is well known that
  $\sum_{\gamma\in[g]} f(\gamma)$ for $f\in \complexs[\Gamma]$ extends
  to a finite normal trace $t_g$ on $\NeumannN\Gamma$, a so called
  delocalized trace. The indices
  generated by these traces are studied by L{\"u}ck in\cite{MR1914619}.

  In \cite{Atiyah(1976)}, Atiyah is working with the lifted operator
  to $\tilde M$: lift the differential (and hence local) operator $D$
  to $\tilde D\colon \Gamma(\tilde E^+)\to \Gamma(\tilde E^-)$, where
  $\tilde E^{\pm}$ are the pullbacks of $E^{\pm}$ to the universal
  covering $\tilde M$.

  In this situation, there is a literal translation between
  spaces of sections and operators on them for $\tilde E^{\pm}$ on the
  one hand, and for $E^{\pm}\tensor H$ on the other hand. This is
  rather straigtforward (and well 
  known). For the sake of completeness we indicate the
  constructions. Other accounts (with more details) can be found e.g.~in \cite[Section
  3.1]{SchickICTP} and \cite[Example 3.39]{SchickBologna}.

  The translation is
  summarized in the following table

  \begin{tabular}[t]{c|c}\\ 
    $\tilde M$ & $\cdot \tensor H$\\ \hline
    $L^2(\tilde E^{\pm})$ & $L^2(E^{\pm}\tensor H)$\\
    $\{s\in\Gamma(E^{\pm})\mid \sum_{\gamma\in\Gamma} \abs{s(\gamma
      x)}^2<\infty\;\forall x\in \tilde M\}$ &  $\Gamma(E^{\pm}\tensor
    H)$ \\
    $\tilde D$ & $D_H$\\
   $ \tilde D/(1+\tilde D^2)^{1/2}$ & $D_H/(1+D_H^2)^{1/2}$\\
   $\phi(\tilde D)$ & %/(1+\tilde D^2)^{1/2}\right)$ &
   $\phi\left(D_H\right)$\\ %/(1+D_H^2)^{1/2}\right)$\\
   $\int_{\tilde M/\Gamma} \tr_x k(x,x)\;dx$ & $t$\\
   $\sum_{\gamma\in [g]}\int_{\tilde M/\Gamma} \tr_x k(x,\gamma
   x)\;dx$ & $t_g$\\
   $\ind_t(\tilde D)$ & $\ind_t(D_H)$\\
   $\ind_{t_g}(\tilde D)$ & $\ind_{t_g}(D_H)$
  \end{tabular}

  Some explanations are in order:
  \begin{enumerate}
  \item A section $s$ of $\tilde E$ corresponds to the section $\hat
    s$ of $E\tensor H$ with $\hat s(x) = \sum_{\gamma\in\Gamma}
    s(\gamma\tilde x)\tensor (\tilde x,\gamma)$, where $\tilde x$ is an arbitrary
    lift of $x$. Of course we identify the fibers $E_x$ and $\tilde
    E_{\gamma\tilde x}$, and $H_x= \Gamma\tilde x \times_\Gamma l^2(\Gamma)$. This construction is well defined by the
    definition of the twisted bundle $H$, with fiber identified with
    $l^2(\Gamma)$ using the chosen lift $\tilde x$.
  \item This identification defines an isometry of the spaces of
    $L^2$-sections. Moreover, it is compatible with the
    $\Gamma$-action, therefore an isometry of $A$-Hilbert spaces. In
    addition, it preserves smoothness and continuity,
    where the condition as given in the table is used to really get a
    section of $E\tensor H$.
  \item The operators $\tilde D$ and $D_H$ are conjugated to each
    other under the isomorphism of the section spaces. This follows
    from their local definition. Here we use that for a small
    connected neighborhood $U$ of $x\in M$ we can choose a lift
    $\tilde U$, a connected neighborhood of a lift $\tilde x$, such
    that there is a unique section $U\to \tilde U$ of the restriction
    of the covering $\tilde M\to M$ to $U$, and then $y\mapsto (\tilde
    y,\gamma)$ is a flat section of $H|_U$ for each $\gamma\in\Gamma$.
  \item Since the self-adjoint unbounded operators $\tilde D$ and $D_H$ are
    unitarily equivalent, the same is true for all bounded measurable
    functions of them, using functional calculus. In particular, this
    is the case for $\tilde D/(1+\tilde D^2)^{1/2}$, but also for any
    other bounded measurable function $\phi\colon\reals\to\reals$.
    As a particular example we will have to study the projections onto
    the kernels of the operators.
  \item Appropriate functions of $\tilde D$, e.g.~the projection onto
    the kernel, have by elliptic regularity a smooth integral kernel $k(x,y)$
    on $\tilde M\times\tilde M$. This kernel is invariant (in the
    appropriate sense) under the
    diagonal $\Gamma$-action, in particular, its restriction to the
    diagonal descends to the quotient by this action. On the diagonal,
    $k(x,x)$ is an endomorphism of the fiber $\tilde E_x$ and
    therefore has a finite dimensional trace $\tr k(x,x)$. Since the
    fiber $\tilde E_{\gamma x}$ can for each $\gamma\in\Gamma$ be
    canonically identified with $\tilde E_x$ (since they are both 
    identified with $E_{p(x)}$, $p\colon\tilde M\to M$), we
    can also take the finite dimensional trace $\tr k(x,\gamma x)$.

    The integrals in the tables define then certain traces which are
    the ones used by Atiyah and by L{\"u}ck.
  \item Choose a subset $U\subset M$ such that $M\setminus U$ has
    measure zero and such that the restriction of the covering $\tilde
    M\to M$ to $U$ is trivial. If we choose an appropriate lift
    of $U$
    then $\tilde M|U \iso U\times \Gamma$. This induces a trivialization
    $H|_U\iso U\times l^2(\Gamma)$. Using this, we identified in Lemma
    \ref{lem:sections are A Hilbert space} $L^2(E\tensor H) =
    L^2(E|_U)\tensor l^2(\Gamma)$, and this in turn was used to define
    $t$ and $t_g$ on (trace class) operators acting on $L^2(E\tensor
    H)$,  e.g.~the projection onto the kernel of $D_H$.

    On the other hand, using the corresponding trivialization of the
    covering $\tilde M|U\iso U\times \Gamma$ we get the identification
    $L^2(\tilde E|U)\iso
    L^2(E)\tensor l^2(\Gamma)$, and our unitary identification defined
    above becomes the identity under these identifications.

    It was proved by Atiyah in \cite{Atiyah(1976)} that the formula of
    the integral computes the tensor product of the ordinary Hilbert
    space trace on $L^2(E)$ with the trace $t$ on $l^2(\Gamma)$ under
    the last identification. This proof extends to the second
    integral, which corresponds to the tensor product of the Hilbert
    space trace on $L^2(E)$ with the delocalized trace $t_g$ on
    $l^2(\Gamma)$.

    On the other hand, we defined $t$ (or $t_g$, respectively) for
    operators on $L^2(E\tensor H)$ as tensor product of $t$ (or $t_g$)
    on $l^2(\Gamma)$ with the usual trace on $L^2(E)$, using the
    identification $L^2(E\tensor H)\iso L^2(E|_U)\tensor l^2(\Gamma)$.
    Since all these identifications coincide with each other, the
    traces also do so.
  \item From the discussion so far, it follows in particular that the
    unitary isomorphism described above induces $A$-Hilbert space
    isometries between $\ker(\tilde D^{\pm})$ and $\ker(D^{\pm}_H)$,
    such that the traces of the projectors onto these kernels coincide, defined
    either using the integral over the diagonal in $\tilde M\times
    \tilde M/\Gamma$ for the integral
    kernel, or using the recipe of Definition \ref{def:dim_t} with
    the Hilbert $A$-module structure given by Lemma \ref{lem:sections
      are A Hilbert space} on $L^2(E\tensor H)$.

  In particular $\ind_t(D_H)=\ind_t(\tilde D)$, and
  $\ind_{t_g}(D_H)=\ind_{t_g}(\tilde D)$, where the left hand side is
  defined in Definition \ref{def:ind_t}, and the right hand side is
  defined with the integrals of the table evaluated for the
  projection operators $k^{\pm}$ onto the kernels of $\tilde D^+$ and
  $\tilde D^-$:
  \begin{equation*}
    \begin{split}
      \ind_t(\tilde D) = &\int_{\tilde M/\Gamma} \tr_x k^+(x,x)\;dx -
      \int_{\tilde M/\Gamma} \tr_x k^-(x,x)\;dx,\\
      \ind_{t_g}(\tilde D) = &\sum_{\gamma\in[g]}\left(\int_{\tilde
          M/\Gamma} \tr_x k^+(x,\gamma x)\;dx -
      \int_{\tilde M/\Gamma} \tr_x k^-(x,\gamma x)\;dx\right).
  \end{split}
\end{equation*}
  \end{enumerate}

In particular, we have proved:
\begin{theorem}\label{theo:covering_equal_twist}
  The $L^2$-index defined in terms of a covering equals the
  $L^2$-index using the corresponding flat $A$-Hilbert space twisting
  bundle.
\end{theorem}

Therefore, we will have proved Theorem \ref{theo:Atiyah} and then recovered
Atiyah's $L^2$-index theorem as soon as we prove
the index formula for $\ind_t(D_H)$, which we will reduce to Theorem
\ref{theo:index_theorem} by proving Equation \eqref{eq:A Hilbert
  reduces to Hilbert A}. 

  Note that Atiyah defines the $L^2$-index for arbitrary elliptic
  differential operators on $M$, not necessarily of Dirac type. This
  is possible since $\tilde M\times_\Gamma l^2(\Gamma)$ is a flat
  bundle, and arbitrary differential operators can be twisted with
  every flat bundle. A corresponding construction is possible in our
  more general setting. Since all geometrically important operators
  are generalized Dirac operators, and since only those can be twisted
  with bundles with non-flat connections, we will stick to the latter
  more restricted class.

\subsection{Equivalences of categories}
\label{sec:an-equiv-categ}

In this section we show how one can go back and forth between Hilbert
$A$-modules and $A$-Hilbert spaces, and the corresponding bundles.

\begin{lemma}\label{lem:l2_V_remains_finitely_generated_projective}
  If $V$ is a finitely generated projective Hilbert $A$-module, then
  $l^2(V)$ is a finitely generated projective $A$-Hilbert
  space. 
%Conversely, every finitely generated projective $A$-Hilbert
%  space is isomorphic to $l^2(V)$ for a suitable finitely generated
%  projective Hilbert $A$-module $V$
\end{lemma}
\begin{proof}
  Let $V\oplus W\iso A^n$ be a decomposition into $V$ and an
  orthogonal complement $W$. Then $V$ and $W$ are orthogonal also with
  respect to the inner product $\innerprod{\cdot,\cdot}_2$, and
  therefore their completions add up to the completion $l^2(A)^n$ of
  $A^n$. Moreover, the projection $A^n\to A^n$ with image $V$ is given
  (as is any right $A$-linear map from $A^n$ to itself) by
  multiplication from the left with a matrix with entries in $A$. This
  same matrix will act on $l^2(A)^n$ (by Lemma
  \ref{lem:A_acts_boundedly_on_l2A}) with kernel containing $W$
  (i.e.~also its closure $l^2(W)$) and
  image containing $V$ and ---since the matrix is a projection--- also
  its closure $l^2(V)$). This shows that the orthogonal projection is
  given by multiplication with the matrix. This completes the proof
  that $l^2(V)$ is a finitely generated projective $A$-Hilbert space.
\end{proof}

\begin{lemma}\label{lem:module_maps_are_l2_bounded}
  Assume that $f\colon V\to W$ is an adjointable $A$-module
  homomorphism between Hilbert
  $A$-modules $V$ and $W$. Then $f$ extends to a bounded $A$-linear operator
  $f\colon l^2(V)\to l^2(W)$ with adjoint the extension of $f^*$.

  If $f\colon V\to W$ is a Hilbert $A$-module isometry, then $f$
  extends to an isometry $f\colon l^2(V)\to l^2(W)$.
\end{lemma}
\begin{proof}
  By \cite[Proposition 1.2]{Lance(1995)}, $\innerprod{f(x),f(x)}\le
  \norm{f}^2\innerprod{x,x}$ in $A$. Therefore, because of positivity
  and linearity of $\tau$
  \begin{equation*}
    \innerprod{f(x),f(x)}_2=\tau(\innerprod{f(x),f(x)}) \le \norm{f}^2
    \tau(\innerprod{x,x})=\norm{f}^2\innerprod{x,x}_2,\qquad\forall
    x\in V.
  \end{equation*}
This shows that $f$ is $l^2$-bounded.

For the adjoint observe that
\begin{equation*}
  \innerprod{f(x),y}_2=\tau(\innerprod{f(x),y})=\tau(\innerprod{x,f^*(y)}) = \innerprod{x,f^*(y)}_2\qquad\forall x\in V. 
\end{equation*}

  If $f\colon V\to W$ is an isometry, then in particular
  \begin{equation*}
\innerprod{f(v),f(v')}_2=\tau(\innerprod{f(v),f(v')})
  =\tau(\innerprod{v,v'})=\innerprod{v,v'}_2\qquad\forall v,v'\in V.
\end{equation*}
\end{proof}

\begin{definition}\label{def:l2_of_bundle}
    Let $W$ be a Hilbert $A$-module
    bundle on a space $X$. Fiberwise application
    of the construction of Lemma \ref{lem:A_acts_boundedly_on_l2A}
    produces an $A$-Hilbert space bundle on $X$ which we call
    $l^2(W)$. The transition functions are obtained as extensions of
    Hilbert $A$-module isometries to $A$-Hilbert space isometries as
    described in Lemma \ref{lem:A_acts_boundedly_on_l2A}. In
    particular, we define an induced smooth structure on $l^2(W)$ from
    a smooth structure on $W$.
\end{definition}

\begin{lemma}\label{lem:connection_on_l2_of_bundle}
  Assume that $W$ is a smooth Hilbert
  $A$-module bundle on a smooth manifold $M$. Let $\nabla$ be a
  connection on $W$ which is locally given by the $\End(W)$-valued
  $1$-form $\omega$ as in Proposition
  \ref{prop:curvature_and_its_properties}, with curvature $2$-form
  $\Omega$. Then the
  connection extends to $l^2(W)$, locally given by $\omega$ and with
  curvature $\Omega$, where we extend the endomorphisms of $W$ to
  endomorphisms of $l^2(W)$ using Lemma
  \ref{lem:module_maps_are_l2_bounded}.

  This extension still satisfies the Leibnitz rule for the right
  $A$-action. If $\nabla$ is a metric connection, the same is true for
  its extension (now with respect to the $l^2$-inner product).
\end{lemma}
\begin{proof}
  Recall that, if a trivialization $W|_U\iso V\times U$ is given, then
  $\nabla=\nabla_0+\omega$, where $\nabla_0$ is the trivial connection
  given by the trivialization. The latter one extends to the
  trivialized bundle $l^2(V)\times U$ as the trivial connection. By
  Lemma \ref{lem:module_maps_are_l2_bounded} $\omega$ extends to a
  $1$-form with values in $A$-Hilbert space endomorphisms of
  $l^2(V)$. Consequently, $\nabla_0+\omega$ defines the desired
  extension of $\nabla$. From the local formula for the curvature of
  Proposition \ref{prop:curvature_and_its_properties}, its curvature
  is the extensions of $\Omega$.

  The Leibnitz rule holds for the trivial connection on $l^2(V)\times
  U$ by the usual calculus proof of the Leibnitz rule (which only uses
  distributivity in both variables), and since $\omega$ is compatible
  with the $A$-module structure also for the extension of $\nabla$.

  If $\nabla$ is a metric connection of $W$, then $\omega$ has values
  in skew adjoint $A$-module endomorphisms. By Lemma
  \ref{lem:module_maps_are_l2_bounded} the extension has values in
  skew adjoint Hilbert space endomorphism and therefore the extension
  of $\nabla$ is a metric connection for the $l^2$-inner product.
\end{proof}

\begin{definition} Assume that $A$ is a von Neumann algebra.
  Let $X$ be an $A$-Hilbert space. Choose an embedding $X\into
  H_X\tensor l^2(A)$ for an appropriate Hilbert space $H_X$ (finite
  dimensional if $X$ is finitely generated projective), as in
  Definition \ref{def:of_A_Hilbert_space}. Let $p\in
  \boundedops(H_X)\tensor A$ be the corresponding orthogonal projection
  onto $X$. Set 
  \begin{equation*}
  A(X):= p(H_X\tensor A) \subset X,
\end{equation*}
where $H_X\tensor A\subset
  H_X\tensor l^2(A)$ is the canonical Hilbert $A$-module contained in
  $H_X\tensor l^2(A)$ (isomorphic to $H_A$ is $H_X$ is separable). Since
  $p$ is a projection in $\boundedops(H_X)\tensor
  A=\boundedops_A(H_X\tensor A)$, the image $p(H_X\tensor A)$ is itself a
  Hilbert $A$-module with the induced structure from the ambient space
  $H_X\tensor A$.

  If $X$ is a finitely generated projective $A$-Hilbert space, $H_X$ can
  be chosen finite dimensional, say~$H_X=\complexs^n$. Then $A(X)$ is a finitely generated
  projective Hilbert $A$-module, the image of the projection $p\in
  \boundedops_A(\complexs^n\tensor A) = M_n(A)$.
\end{definition}

Of course, the construction of $A(X)$ a priori depends on the choice
of the projection $p$. In the next lemma, we will see that this is not
the case.

\begin{lemma}\label{lem:functoriality_of_passing_to_A}
  A bounded $A$-linear operator $f\colon X\to Y$ between two
  $A$-Hilbert spaces induces by restriction an adjointable $A$-linear
  map $A(f)\colon A(X)\to A(Y)$, for every choice of projection
  $p_X\in \boundedops(H_X)\tensor A$ and
  $p_Y\in\boundedops(H_Y)\tensor A$ with image $X$ and $Y$,
  respectively. Moreover, $A(f)^*=A(f^*)$ and $A(\cdot)$ is a
  functor. If $f$ is a Hilbert space isometry, then
  $A(f)$ is an isometry of Hilbert $A$-modules.

  In particular, if we apply this to $\id_X\colon X\to X$, with $A(X)$
  defined using two different projections, we see that $\id_X$
  restricts to the identity map on $A(X)$, therefore $A(X)$ (with its
  structure as Hilbert $A$-module) is well defined.
\end{lemma}
\begin{proof}
  If $i_Y\colon Y\to H_Y\tensor l^2(A)$ is the inclusion, then
  \begin{equation*}
i_Y\circ f\circ p_X\colon H_X\tensor l^2(A)\to H_Y\tensor l^2(A)
\end{equation*}
is a bounded operator which commutes with right multiplication by
$A$. Since $A$ is a von Neumann algebra, by Lemma
\ref{lem:calculate_commutant} the composition belongs to
$\boundedops(H_X,H_Y)\tensor A$, where $A$ acts by right
multiplication on $l^2(A)$. In particular, the subspace $H_X\tensor A$
is mapped to the subspace $H_Y\tensor A$, and since $A(X)$ is the
intersection $X\cap (H_X\tensor A)$, and similarly $A(Y)=Y\cap
(H_Y\tensor A)$, $f$ maps these subspaces to each other.

Moreover, $\boundedops(H_X,H_Y)\tensor A$ is exactly the space of
adjointable operators from $H_X\tensor A$ to $H_Y\tensor A$. Since
$A(f) = p_Y\circ (i_Y f p_X) \circ i_X$, and $p_Y$, $i_X$ are also
adjointable, the same follows for $A(f)$.

$A(f)$ is functorial by construction, since it is just given by
restriction to the subspace $A(X)$. Since the representations of $A$ on $l^2(A)$
by left and right multiplication are both $C^*$-homomorphisms,
$\boundedops(H_X,H_Y)\tensor A\to \boundedops(H_X\tensor l^2(A),
H_Y\tensor l^2(A))$ is also adjoint preserving. It follows that
$A(f)^*=A(f^*)$.

Finally, $f$ is an isometry $\iff$ $ff^*=1=f^*f$ $\iff$
$A(f)A(f)^*=1=A(f)^* A(f)$ $\iff$ $A(f)$ is an isometry. 
\end{proof}

Note that for Lemma \ref{lem:functoriality_of_passing_to_A} it is
crucial that $A$ is a von Neumann algebra, the corresponding result
does not necessarily hold for arbitrary $C^*$-algebras.

We needed the following lemma.

\begin{lemma}\label{lem:calculate_commutant}
  Let $A$ be a von Neumann algebra with a trace $\tau$ as in
  \ref{not:special_tau}. Then $A$ acts by left and right
  multiplication on $l^2(A)$. The corresponding subalgebras of
  $\boundedops(l^2(A))$ are mutually commutants of each other, i.e.~the
  operators given by right multiplication with elements of $A$ are
  exactly those operators commuting with left multiplication by $A$.

  Let $H_1$ and $H_2$ be two Hilbert spaces. Then
  \begin{equation*}
\boundedops(H_1\tensor l^2(A),H_2\tensor
  l^2(A))^A=\boundedops(H_1,H_2)\tensor A,
\end{equation*}
where $\boundedops(H_1\tensor l^2(A),H_2\tensor
  l^2(A))^A$ is defined as those operators commuting with left
  multiplication by $A$, and the factor $A$ in
  $\boundedops(H_1,H_2)\tensor A$ acts by right multiplication on $l^2(A)$.
\end{lemma}
\begin{proof}
    The first assertion follows from Tomita modular theory. The vector
  $1\in l^2(A)$ is a separating and generating vector for left as well
  as right multiplication of $A$ on $l^2(A)$ since the trace is
  faithful, and since, by definition, $l^2(A)$ is the closure of the
  subspace $A$. The map 
  \begin{equation*}
J=S=F\colon A\to A;\; a\mapsto a^*
\end{equation*}
is a conjugate linear
  isometry of order $2$, in particular extends to all of $l^2(A)$.

  By \cite[Theorem 9.2.9]{Kadison2} the elements of the commutant of
  right multiplication $R_a$ with elements  $a\in A$ are given as
  operators $J R_a J = L_{a^*}$,
  $a\in A$ (where $L_a$ denotes left multiplication with $A$). The
  first statement follows.

  The second assertion follows since the commutant of $A_1\tensor A_2$
  acting on $H_1\tensor H_2$ is $A_1'\tensor A_2'$ (here
  $A_1=\complexs$, $A_1'=\boundedops(H_1,H_2)$).
\end{proof}

\begin{theorem}
  Let $A$ be a von Neumann algebra with a trace $\tau$ as in
  \ref{not:special_tau}. The category of finitely generated projective
  $A$-Hilbert spaces is
  equivalent to the category of finitely generated projective Hilbert
  $A$-modules, and the category of $A$-Hilbert spaces is equivalent
  to the category of projective Hilbert $A$-modules. The equivalence
  is given by $V\mapsto l^2(V)$ and $X\mapsto A(X)$ for any Hilbert
  $A$-module $V$ and $A$-Hilbert space $X$.
\end{theorem}
\begin{proof}
  This follows from Lemma \ref{lem:functoriality_of_passing_to_A} and
  Lemmas \ref{lem:module_maps_are_l2_bounded} and
  \ref{lem:l2_V_remains_finitely_generated_projective}. 
\end{proof}

\begin{proposition}\label{prop:eq_of_cat}
   Assume that $A$ is a von Neumann algebra with a
  trace $\tau$ as in \ref{not:special_tau}.
  The naturality of the construction of $A(X)$ for an $A$-Hilbert
  space $X$ implies that we get a corresponding functor which assigns
  to each finitely generated projective (smooth) $A$-Hilbert space bundle a
  finitely generated projective (smooth) Hilbert $A$-module
  bundle. Here we also use that the transition functions (in both
  cases isometries) are preserved since the functors map isometries to
  isometries. Together with the construction of Definition
  \ref{def:l2_of_bundle} this gives rise to an equivalence between
  finitely generated projective (smooth) $A$-Hilbert space bundles and
  finitely generated projective Hilbert $A$-module bundles.

  A connection on a smooth finitely generated projective $A$-Hilbert
  space bundle preserves the Hilbert $A$-module subbundle and
  therefore gives rise to a connection on the latter. In view of Lemma
  \ref{lem:connection_on_l2_of_bundle}, we
  also get an equivalence between smooth Hilbert $A$-module bundles
  with connection and smooth $A$-Hilbert space bundles with connection.
\end{proposition}
\begin{proof}
  We only have to check that a connection on an $A$-Hilbert space
  bundle indeed preserve the Hilbert $A$-module subbundle. This is
  clear for the trivial connection on a trivial bundle $U\times
  X$. Locally, an arbitrary connection differs from such a trivial
  connection by a one form with values in endomorphisms which commute
  with the right $A$-multiplication. Using Lemma
  \ref{lem:calculate_commutant} in the same way as in the proof of
  Lemma \ref{lem:functoriality_of_passing_to_A}, such endomorphisms preserve the
  Hilbert $A$-module subbundle, and therefore the same is true for the
  connection.
\end{proof}

\begin{corollary}\label{corol:find_A_Hilbert_bundle}
  Given any smooth finitely generated projective $A$-Hilbert space
  bundle $X$ with connection, we can assume that $X=l^2(V)$ for an
  appropriate smooth finitely generated projective Hilbert $A$-module
  bundle $V$ with connection, where the connection on $l^2(V)$ is
  obtained as described in Lemma
  \ref{lem:connection_on_l2_of_bundle}.
\end{corollary}

\subsection{The general version of Atiyah's $L^2$-index theorem}
\label{sec:gener-vers-atiy}

In view of Corollary \ref{corol:find_A_Hilbert_bundle} we can now
formulate our general version of the $L^2$-index theorem.

\begin{theorem}\label{theo:general_Atiyah}
  Let $M$ be a closed manifold, and $D\colon \Gamma(E^+)\to
  \Gamma(E^-)$ a generalized Dirac operator
  on $M$. Let $A$ be a von Neumann algebra with a normal trace $t$ and
  a faithful trace
  $\tau$ as in Definition \ref{def:dim_t}. Let $X$ be a smooth
  finitely generated projective $A$-Hilbert space
  bundle on $M$, obtained (by Corollary
  \ref{corol:find_A_Hilbert_bundle}) as $X=l^2(V)$ for a smooth finitely generated
  projective Hilbert $A$-module bundle $V$. Assume that $X$ has a
  connection which is extended from $V$ as in Lemma
  \ref{lem:connection_on_l2_of_bundle} and Proposition
  \ref{prop:eq_of_cat}.
   Then
   \begin{equation*}
     \ind_t(D_X) = t(\ind(D_V)),
   \end{equation*}
  where $\ind_t(D_X)$ is defined in Definition \ref{def:ind_t}, and
  $\ind(D_V)\in K_0(A)$ is defined in Definition \ref{def:ind_D_V}. In
  particular, by Theorem \ref{theo:index_theorem}
  \begin{equation*}
     \ind_t(D_X) = \innerprod{ \ch(\sigma(D)) \cup \Td(T_\complexs
      M) \cup
      \ch_t(V) , [TM]}.
  \end{equation*}
  We might as well define $\ch_t(X):=\ch_t(V)$ and observe that it can
  be obtained from the connection on $X$ (which gives rise to the
  connection on $V$ simply by restriction). In particular, if $X$ (or equivalently
  $V$) are flat, then
  \begin{equation*}
      \ind_t(D_X) = \innerprod{ \ch(\sigma(D)) \cup \Td(T_\complexs
      M), [TM]}\cdot \dim_t(X_p),
  \end{equation*}
  where $\dim_t(X_p)$ is the locally constant function (with values in
  $Z$) which assigns to $p\in M$ the value $\dim_t(X_p)=\dim_t(V_p)$,
  where $X_p$ and $V_p$ are the fibers over $p$ of $X$ and $V$, respectively.
\end{theorem}

\begin{corollary}\label{corol:l2-equals ordinary index}
  If $A$ in Theorem \ref{theo:general_Atiyah} is a finite von Neumann
  algebra with center valued trace $t\colon A\to Z$, then
  $\ind_t(D_X)$ and $\ind(D_V)$ can be obtained from each other.
\end{corollary}
\begin{proof}
  One direction follows from $t(\ind(D_V))=\ind_t(D_X)$. The converse
  is true because the center valued trace induces an \emph{injection}
  $K_0(A)\xrightarrow{t} Z$ by \ref{prop:explicit_formula_for_chern},
  applied to $X=\{*\}$.
\end{proof}

  Theorem \ref{theo:general_Atiyah} is a consequence of Corollary
  \ref{corol:ind_D_V_with_ker_projection} and of properties of the trace
  $t$ established in Section \ref{sec:properties-traces}. Therefore,
  we first establish these properties of $t$, before completing the
  proof of Theorem \ref{theo:general_Atiyah}. 

\subsection{Properties of traces}
\label{sec:properties-traces}

In Definition \ref{def:dim_t} we used the extension of the trace $t$
from $A$ to $\boundedops(H)\tensor A$. Here, we want to recall the
definition and the main properties (we are following \cite[I 6,
Exercise 7]{Dixmier}). Similar considerations can be found in
\cite[Section 2]{Schick(2001)}.

\begin{definition}\label{def:of_t_on_H_tensor_A}
  Let $A$ be a von Neumann algebra with a trace $\tau$ as in
  \ref{not:special_tau} and with a normal trace $t\colon A\to
  Z$, where $Z$ is a commutative  von Neumann algebra
  (e.g.~$Z=\complexs$). Let $H$ be a Hilbert
  space with orthonormal basis $\{e_i\mid i\in I\}$. For a positive
  operator $a\in \boundedops(H)\tensor A$ (acting on $H\tensor
  l^2(A)$) define
  \begin{equation*}
    t(a):= 
    \begin{cases}
      \sum_{i\in I} t(U_i^* a U_i)\in Z,& \text{if the sum is ultraweakly
        convergent}\\
      \infty &\text{otherwise}
\end{cases}
  \end{equation*}
  where $U_i\colon l^2(A)\to H\tensor l^2(A)$ is given by the
  decomposition of $H$ according to the orthonormal basis
  $\{e_i\}$. Note that $U_i^* aU_i\in A$, since the map $a\mapsto
  U_i^*aU_i$ is norm
  continuous from $\boundedops(H\tensor l^2(A))\to \boundedops(l^2(A))$ and maps
  elementary tensors  $T\tensor x\in \boundedops(H)\tensor A$ to
  elements of $A$. Note that $\sum_{i\in I} t(U^*_i
  aU_i)$ is an infinite sum of \emph{non-negative} elements. It is
  convergent if and only if the corresponding collection of finite
  sums has an upper bound in $Z$, in which case the least upper bound
  is the limit. In particular, convergence is independent of the
  ordering in the sum.

  The linear span of all positive operators $a$ with $t(a)<\infty$ is an
  ideal in $\boundedops(H)\tensor A$, and $t$ extends by linearity to
  this ideal.
\end{definition}

In the above definition, we must check that $t(a)$ does not depend on
the chosen orthonormal basis $\{e_i\}$. If $f_j$ is a second
orthonormal basis with induced unitary inclusions $V_j\colon l^2(A)\to
H$, then this follows from the
following calculation
\begin{equation*}
  \begin{split}
    \sum_{i\in I} t(U_i^* aU_i) =& \sum_{i\in I} t(U_i^* \sum_{j\in J}
    V_jV_j^* a U_i) \\
    =& \sum_{i\in I}\sum_{j\in J} t(U_i^*V_j V_j^* a U_i)\\
    =& \sum_{i\in I,j\in J} t(V_j^* a U_i U_i^* V_j)\\
    =& \sum_{j\in J} t(V_j^* a \sum_{i\in I} U_i U_i^* V_j) =
    \sum_{j\in J} t(V_j^* a V_j).
  \end{split}
\end{equation*}
  Here we used the fact that $\sum_{i\in I} U_i U_i^*= \sum_{j\in J}
  V_jV_j^* =\id_{H\tensor l^2(A)}$, where the convergence is in the
  ultraweak sense, and that $t$ is normal and a trace. 

  Moreover, we
  use that the linear map $a\mapsto U_i^* a V_j\colon \boundedops(H)\tensor
  A\to \boundedops(l^2(A))$ is norm continuous and maps elementary
  tensors $T\tensor x\in \boundedops(H)\tensor A$ to elements of $A$,
  such that the image is contained in $A$. In particular $U_i^*V_j
  =U_i^* 1 V_j\in A$ and $V_j^* a U_i\in A$, such that $t( (U_i^* V_j)
  (V_j^* a U_i)) = t((V_j^* aU_i) (U_i^* V_j))$ by the trace property
  for operators in $A$.

  Again, since all the summands in the above infinite sums are
  positive elements of $Z$, the ordering is not an issue, and the
  limit (if it exists) is the least upper bound.

\begin{definition}\label{def:general_trace_class}
  Let $A$ be a von Neumann algebra with traces $\tau$ and $t$ as above.
    
  Assume that $V_1$ and $V_2$ are $A$-Hilbert spaces and $f\colon
  V_1\to V_2$ is an $A$-linear bounded operator. Let $i_1\colon V_1\to
  H_1\tensor l^2(A)$ and $i_2\colon V_2\to H_2\tensor l^2(A)$ be
  inclusions as in Definition \ref{def:of_A_Hilbert_space}, and $p_1$,
  $p_2$ the corresponding orthogonal projections. We say that $f$ is a
  \emph{$t$-Hilbert Schmidt} operator, if $i_1 f^*f p_1$ is of
  $t$-trace class. We say that $f$ is of $t$-trace class, if there are
  $f_1\colon V_1\to V_3$ and $f_2\colon V_3\to V_1$ $t$-Hilbert
  Schmidt operators ($V_3$ an additional $A$-Hilbert space) such that
  $f=f_2f_1$.

  If $V_1=V_2$ and $f$ is of $t$-trace class, set $t(f):= f(i_1fp_1)$.

  If $\id_{V_1}$ is of $t$-trace class, define
  $\dim_t(V_1):=t(\id_{V_1})$, else set $\dim_t(V_1):=\infty$.
\end{definition}

Again, it is necessary to check that the definitions in
\ref{def:general_trace_class} are independent of the choices
made. Moreover, we have to check that the trace so defined has the
usual properties (which we are going to use later on). This is the
content of the following theorem. Essentially the same theorem, with
$t$ complex valued, is stated in \cite[Theorem 2.3]{Schick(2001)} and
\cite[9.13]{SchickDiss}. The
proof given there also applies to the more general situation here.

\begin{theorem}\label{theo:prop_of_t}
  Assume that $A$ is a von Neumann algebra with traces $\tau$ and $t$
  as above.

  Let $V_0, V_1$, $V_2$ and $V_3$ be $A$-Hilbert spaces and $f\colon V_1\to
  V_2$, $g\colon  V_2\to V_3$, $e\colon V_0\to V_1$ be bounded
  $A$-linear operators. Then:
\begin{enumerate}
\item\label{adjoint}\label{tr_descr} $f$ is of $t$-trace class $\iff$
  $f^*$ is of $t$-trace class  $\iff$ $\abs{f}$ if of $t$-trace class
\item  $f$ is a $t$-Hilbert-Schmidt operator $\iff$ $f^*$ is a
  $t$-Hilbert-Schmidt operator.
\item\label{HS_comp} If $f$ is a $t$-Hilbert-Schmidt operator then
  $gf$ and $fe$ are $t$-Hilbert-Schmidt operators.
\item\label{tr_comp} If $f$ is a $t$-trace class operator, then $gf$
  and $fe$ are $t$-trace class operators.
\item\label{C_weak} If $f$ is of $t$-trace class and $V_1=V_3$ then
  $g\mapsto t(gf)$ is ultra-weakly continuous.
\item\label{tr_b}\label{tr_HS} If $V_1=V_3$ and either $f$ if of
  $t$-trace class or $f$ and $g$ are $t$-Hilbert-Schmidt operators
  then $t(gf)=t(fg)$.
\item\label{HS_prod}\label{tr_prod} If $V_{1,2}=H\tensor l^2(A)$ for a
  Hilbert space $H$,
  $a$ is a $t$-Hilbert-Schmidt operator and $B\in\boundedops(H)$ is a
  Hilbert-Schmidt operator, then $f=a\tensor
  B$ is $t$-Hilbert-Schmidt operator.
If $a$ is of $t$-trace class and $B$ is of trace class, then $f$ is of
$t$-trace class with $t(f)=t(a)Sp(B)$, where $Sp$ is the ordinary
trace on the trace class ideal of $\boundedops(H)$.
\item \label{item:dim_t_topolo_inv} Assume that $u\colon V_1\to V_2$
  is bounded $A$-linear with a bounded (necessarily $A$-equivariant)
  inverse $u^{-1}$. Then $\dim_t(V_1)=\dim_t(V_2)$, i.e.~$\dim_t$ does
  not depend on the Hilbert space structure.
\end{enumerate}
\end{theorem}
\begin{proof}
  \begin{itemize}
  \item[\ref{item:dim_t_topolo_inv}] We have $\dim_t(V_1) =
    \tr_t(\id_{V_1}) = \tr_t(u^{-1}u \id_{V_1}) = \tr_t(u\id_{V_1}
    u^{-1}) = \tr_t(\id_{V_2})$ if either $\id_{V_1}$ or $\id_{V_2}$
    are of $t$-trace class, and the calculation shows that then the
    other one also is of $t$-trace class. Here we used \ref{tr_HS}.
  \end{itemize}
\end{proof}

\subsection{Trace class operators}
\label{sec:trace-class-oper}

\begin{definition}
  Assume that $f\in\End_A(H_A)$ is a self adjoint positive endomorphism of the
  standard countably generated Hilbert $A$-module $H_A$. We call $f$
  of \emph{$\tau$-trace class} if $\tau(f):=\sum_{n\in\naturals}
  \tau(\innerprod{f(e_n), e_n}_A) <\infty$. An arbitrary $f\in
  \End_A(H_A)$ is called a $\tau$-trace class operator if it is a (finite)
  linear combination of self adjoint positive $\tau$-trace class
  operators. Then $\tau(f)$ is defined as the corresponding linear combination.

  Let $V$, $W$ be countably generated Hilbert $A$-modules, $f\in
  \Hom_A(V,W)$. We call $f$ of \emph{$\tau$-trace class}, if $f\oplus
  0\colon V\oplus H_A\to W\oplus H_A$ is of $\tau$-trace class. Recall
  that by Kasparov's stabilization theorem \cite[Theorem
  15.4.6]{Wegge-Olsen} $V\oplus H_A\iso H_A\iso
  W\oplus H_A$ such that being of $\tau$-trace class is already
  defined for $f\oplus 0$. The normality of $\tau$ is used to prove
  that this concept and the extension of $\tau$ we get this way is
  well defined and that we can define traces with the usual properties
  in Proposition \ref{prop:trace_properties}.
\end{definition}

\begin{proposition}\label{prop:trace_properties}
  If $f\in \Hom_A(V,W)$ is of $\tau$-trace class and $g\in\Hom_A(W,V)$
  then $fg$ and $gf$ are both of $\tau$-trace class and
  $\tau(fg)=\tau(gf)$.

  If $g\colon l^2(\naturals)\to l^2(\naturals)$ is of trace class with
  trace $\tr(g)$ (in the sense of endomorphisms of the Hilbert space
  $l^2(\naturals)$), and
  $f\in\End_A(A)$ then $f\tensor g\in \End_A(H_A)$ is of $\tau$-trace
  class and
  \begin{equation*}
    \tau(f\tensor g) = \tau(f) \cdot \tr(g)
  \end{equation*}
\end{proposition}
\begin{proof}
  The trace on $\End_A(H_A)$ is the tensor product of $\tau$ on $A$
  and the standard trace on $l^2(\naturals)$ which both have the trace
  property. For a more detailed treatment of such results compare
  e.g.~\cite[Section 2]{Schick(2001)}.

  Recall that we define $f\tensor g(a e_n) := f(a) g(e_n)$, which extends by
  linearity and continuity to an element of $\End_A(H_A)$.
\end{proof}

\begin{definition}
  Exactly the same kind of definition was made for $A$-Hilbert space
  morphisms. Observe that the two constructions are compatible  in the sense that if $f\in\End_A(V)$ is
  of $\tau$-trace class then the same is true for its extension to
  $l^2(V)$ as in Lemma \ref{lem:module_maps_are_l2_bounded} with
  unchanged trace $\tau(f)$.
\end{definition}

\subsection{Proof of Theorem \ref{theo:general_Atiyah}}
\label{sec:other-poss-proofs}

Note first
  that, by definition, 
  \begin{equation*}
\ind_t(D_X)=t(\pr_{ker(D_X)}) -
t(\pr_{coker(D_X)})
,
\end{equation*}
where $\pr_{\ker(D_X)}$ is the orthogonal
  projection onto the kernel of $D_X$ inside the space of
  $L^2$-section $L^2(E^+\tensor X)$, and $\pr_{\coker(D_X)}$ is the
  projection onto the orthogonal complement of the image of $D_X$ in
  $L^2(E^-\tensor X)$. Here, we consider $D_X\colon L^2(E^+\tensor
  X)\to L^2(E^-\tensor X)$ as unbounded operator.

  $D_X$ also gives rise to a bounded operator between Sobolev spaces.
The following definition should be compared with Definition
\ref{def:of_Sobolev}.

\begin{definition}\label{def:of_Sobolev_2}
  Given a finitely
  generated smooth $A$-Hilbert space bundle $X$ over a compact smooth
  manifold $M$, Sobolev spaces $H^s(X)$ can be
  defined ($s\in\reals$), compare e.g.~\cite{BFKM(1996)}. One way to
  do this is to pick a trivializing
  atlas $(U_\alpha)$ with subordinate partition of unity $(\phi_\alpha)$
  and then define for smooth sections $u,v$ of $X$ the inner
  product
  \begin{equation*}
    (u,v)_s = \sum_\alpha \int_{U_\alpha} \innerprod{(1+\Delta_\alpha)^s
    \phi_\alpha u(x),\phi_\alpha v(x)}\;dx,
  \end{equation*}
  where $\Delta_\alpha$ is the ordinary Laplacian on $\reals^n$ acting
  on the trivialized bundle (in the notation, some diffeomorphisms are omitted).

  The inner product is $\complexs$-valued and the completion is an
  $A$-Hilbert space.
\end{definition}

\begin{theorem}
  Assume that $W$ is a smooth finitely generated projective Hilbert
  $A$-module bundle  over a compact manifold $M$, 
  For each $\epsilon>0$, the natural inclusion $H^s(W)\to
  H^{s-\epsilon}(W)$ is $A$-compact.

  If $r>\dim(M)/2$, then the natural inclusion $H^s(W)\to H^{s-r}(W)$
  is of $\tau$-trace class.

  The second assertion holds also if $W$ is a finitely generated projective $A$-Hilbert
  space bundle.
\end{theorem}
\begin{proof}  
  Using charts and a
  partition of unity, it suffices to prove the statement for the
  trivial bundle $A\times T^n$ on the $n$-torus $T^n$. In the latter
  case, one obtains isomorphisms $H^s(A\times T^n)\iso H^s(T^n)\tensor
  A$. In particular, the inclusion $H^s(A\times T^n)\to
  H^{s-r}(A\times T^n)$ is the tensor product of the inclusion of
  $H^s(T^n)\to H^{s-r}(T^n)$ with the identity on $A$. By Proposition
  \ref{prop:trace_properties}, the trace class property follows, and
  compactness is handled in a similar way.

  The same argument applies to $A$-Hilbert space bundles.
\end{proof}

  A twisted Dirac operator $D_H$ as in Definition
  \ref{def:of_H_Hilbert_twisted_Dirac} extends to a bounded operator
  between Sobolev spaces $D_H\colon H^1(W^+\tensor X) \to
  L^2(W^-\tensor X)$.

  Of course, the inner product on $H^s(W)$ depends on a number of
  choices, However, two different choices give rise to equivalent
  inner products and therefore isomorphic Sobolev spaces.

  Observe that if $V$ is a finitely generated projective Hilbert
  $A$-module bundle with corresponding 
  $A$-Hilbert module completion $X=l^2(V)$, the $A$-Hilbert space
  completion $l^2(H^s(V))$ and $H^s(l^2(V))$ are isomorphic. This
  follows since the trace $\tau$ used to define $l^2(V)$ is continuous
  by Lemma \ref{lem:continuity}. $l^2(H^s(V))$ is the completion of
  $\Gamma(V)$ with respect to the inner product $\sum
  \tau\int_{U_\alpha} \innerprod{(1+\Delta)^\alpha\cdot,\cdot}$,
  whereas $H^s(l^2(V))$ is the completion of $\Gamma(V)$ with respect
  to the inner product $\sum \int_{U_\alpha}
  \tau(\innerprod{(1+\Delta)^s\cdot,\cdot})$ and by continuity, $\tau$
  commutes with integration so that the two inner products coincide.

  Moreover, 
  \begin{equation*}
D_X=l^2(D_V)\colon H^1(E^+\tensor X)\to L^2(E^-\tensor
  X)
\end{equation*}
 under this identification (and is in particular a bounded
  operator). We can now look at $\chi_{\{0\}}(D_X^*D_X)$ and
  $\chi_{\{0\}}(D_XD_X^*)$. These are the projections onto the kernel
  of $D_X$ in $H^1(E^+\tensor X)$ and onto the orthogonal complement
  of the image of $D_X$ in $L^2(E^-\tensor X)$. Note that the second
  space is exactly the same one showing up in the definition of
  $\dim_t(D_X)$, since $H^1$ is exactly the domain of the closure of
  the unbounded operator $D_X$ on $L^2$.

  However, the kernels in $H^1$ and in $L^2$ strictly speaking are different. The
  inclusion $H^1(E^+\tensor X)\to L^2(E^+\tensor X)$ maps the kernels
  bijectively onto each other (by elliptic regularity), but the
  topologies are different. Note,
  however, that $\ker(D_X)\subset L^2(E^+\tensor X)$ is a closed
  subset, therefore complete. By the open mapping theorem, the
  bijection between the kernels has a bounded inverse (which is of
  course also $A$-linear). It follows from Theorem
  \ref{theo:prop_of_t} \ref{item:dim_t_topolo_inv} that
  $\dim_t(\ker(D_X))$ does not depend on the question whether we
  consider $D_x$ as unbounded operator on $L^2$ or as bounded operator
  from $H^1$ to $L^2$. In particular,
  \begin{equation*}
    \ind_t(D_X) = t(\chi_{\{0\}}(D_X^*D_X)) - t(\chi_{\{0\}}(D_XD_X^*)),
  \end{equation*}
  where $D_X$ is considered as bounded operator from $H^1$ to $L^2$.

  Note that, on the level of operators, the functor $l^2$ embeds for
  each Hilbert $A$-module $U$ the
  $C^*$-algebra $\Hom_A(U,U)$ into the $C^*$-algebra
  $\boundedops(l^2(U))$. Embeddings of $C^*$-algebras commute with
  functional calculus. In particular,
  $\chi_{\{0\}}(D_X^*D_X)=l^2(\chi_{\{0\}}(D_V^*D_V))$ and
  $\chi_{\{0\}}(D_XD_X^*)=l^2(\chi_{\{0\}}(D_VD_V^*))$. 

  Next, we must look at $t(\ind(D_V))$. This is defined as follows:
  after stabilization, $L^2(E^+\tensor V)\oplus H_A\iso H_A$. Then,
  there is a unitary $u\in \boundedops_A(H_A)$ such that
  \begin{equation*}
   p:=u^*(\chi_{\{0\}}(D_V^*D_V)\oplus 0_{H_A}) u\in
  \boundedops_A(H_A)
  \end{equation*}
 (using the above
  isomorphism) is a projection which is represented by a matrix with
  finitely many non-zero entries, where we understand $M_n(A)\subset
  \boundedops_A(H_A)$ using an orthonormal basis of $H$ in
  $H_A=H\tensor A$. Similarly, $\chi_{\{0\}}(D_VD_V^*)$
  gives rise to a projection $q$ in $M_n(A)\subset
  \boundedops_A(H_A)$. We can apply the functor $l^2$ to the whole
  construction, and therefore get elements $l^2(p)$ and $l^2(q)$,
  represented exactly by the same finite matrices $p$ and $q$ in
  $M_n(A)$ which are unitarily equivalent (by $A$-linear operators
  $l^2(u)$) to 
  \begin{equation*}
  \chi_{\{0\}}(D_X^*D_X)\oplus 0_{H\tensor l^2(A)} \quad\text{and} \quad
  \chi_{\{0\}}(D_XD_X^*)\oplus 0_{H\tensor l^2(A)}.
  \end{equation*}

  Then $\ind_t(D_V) = t(p)-t(q)$. Because $t$ is normal, we have
  Theorem \ref{theo:prop_of_t} \ref{tr_HS} which is valid for
  non-finitely generated $A$-Hilbert spaces and therefore 
  \begin{equation*}
    t(p) =t(l^2(p))= t(l^2(\chi_{\{0\}}(D_V^*D_V))),\quad
    t(q)=t(l^2(q))= t(l^2(\chi_{\{0\}}(D_VD_V^*))). 
  \end{equation*}
  For the first equal sign in both equations note that $t(p)$
  and $t(l^2(p))$ are by their very definitions exactly the same
  thing.

  This finally implies the assertion of Theorem \ref{theo:general_Atiyah}.

\section{Trace class subalgebras}
\label{sec:trace-class-subalg}

Throughout this paper, we have been dealing with a $C^*$-algebra $A$
with a trace $\tau\colon A\to Z$, $Z$ being a commutative
$C^*$-algebra. We were able to derive rather explicit index theorems
for $\tau(\ind(D_W))$, where $D$ is a Dirac type operator on a closed
manifold $M$ and $W$ is a Hilbert $A$-module bundle on $M$ (with
finitely generated projective fibers). Here $\tau\colon K_0(A)\to Z$
is the
induced map on K-theory and $\ind(D_W)\in K_0(A)$ is defined e.g.~by the
Mishchenko-Fomenko construction.

However, there are many other situations where trace-like maps on
$K_0(A)$ exist which do not come from a trace on $A$. One of the most
prominent is if $\mathcal{B}\subset A$ is a dense subalgebra which is closed
under holomorphic functional calculus in $A$ and which has a trace
$\tau\colon \mathcal{B}\to Z$. Since the inclusion $\mathcal{B}\to A$ induces
an isomorphism $K_0(\mathcal{B})\to K_0(A)$, we still get an induced map $\tau\colon
K_0(A)\to Z$. The most prominent example is the trace class ideal
inside the compact operators on a separable Hilbert space. In the
notation of \cite{MR99g:46104}, $\mathcal{B}$ is a local $C^*$-algebra
with completion $A$.

In this section we describe how to generalize the results obtained in
the rest of the paper to this situation. In particular we will get an
explicit index theorem.

Our goal is to describe how the contents of Sections
\ref{sec:hilb-modul-their} to \ref{sec:mishch-fomenko-index}
remain true almost entirely in the more general situation.

We note that the following concepts are relevant.
\begin{enumerate}
\item Finitely generated projective modules make sense in exactly the
  same way for the local $C^*$-algebra $\mathcal{B}$ as for the
  $C^*$-algebra $A$.
\item Finitely generated projective Hilbert $\mathcal{B}$-modules also
  make sense. Of course we can not assume that such a module is
  complete in any way. But being finitely generated projective, those
  modules are upto isomorphism described as an orthogonal summand of
  $\mathcal{B}^n$, and this makes sense for each
  $*$-algebra. Restricting the $\mathcal{B}$-valued inner product of $\mathcal{B}^n$
  induces a $\mathcal{B}$-valued inner product on the summands.
\item $\mathcal{B}$ is a dense subalgebra of $A$, and in a canonical
  way is each finitely generated projective
  (Hilbert)$\mathcal{B}$-module $\mathcal{V}$ a dense subspace of the
  finitely generated projective (Hilbert) $A$-module
  $V:=\mathcal{V}\tensor_\mathcal{B} A$. In the same way, for any two
  such modules $\mathcal{V}_1$, $\mathcal{V}_2$,
  $\Hom_\mathcal{B}(\mathcal{V}_1,\mathcal{V}_2)$ is a dense 
 subspace of $\Hom_A(V_1,V_2)$. If $V_2=V_1$ and therefore these
 spaces are algebras, then the smaller one is holomorphically closed
 in the bigger one. This follows since
 $\Hom_\mathcal{B}(\mathcal{B}^n,\mathcal{B}^n)=M(n,\mathcal{B})$ is
 holomorphically closed in $\Hom_A(A^n,A^n)=M(n,A)$ by definition of being
 holomorphically closed.
\item Now the definition of a bundle with finitely generated projective
  (Hilbert) $\mathcal{B}$-module fibers makes perfectly sense. This is
  also true for smooth bundles, there a map to
  $\Hom_\mathcal{B}(\mathcal{V}_1,\mathcal{V}_2)$ is smooth if and
  only if it is smooth when composed with the inclusion into
  $\Hom_A(V_1,V_2)$ as above.

  In particular, we can consider each such bundle as included in an
  induced finitely generated projective (Hilbert) $A$-module bundle,
  and a smooth structure on the $\mathcal{B}$-bundle induces a smooth
  structure on the $A$-bundle.
\item The results in section \ref{sec:hilb-modul-their}, in particular
  in Section \ref{sec:struct-finit-gener} remain true for finitely
  generated projective (Hilbert) $\mathcal{B}$-module bundles. This
  follows by carrying out the constructions for the induced $A$-module
  bundle and then observing that all the constructions, which only
  involve the $*$-operator, the algebra structure, and taking holomorphic
  functions of elements in $\Hom_\mathcal{B}(\mathcal{V},\mathcal{V})$
  remain inside this subset of $\Hom_A(V,V)$ by the very fact that
  $\mathcal{B}$ is holomorphically closed in $A$.

  The most important such function takes the inverse of an
  invertible element. One should note that the set of invertible
  elements of $\mathcal{B}$, being the intersection of $\mathcal{B}$
  with the corresponding open subset of $A$ is itself open in
  $\mathcal{B}$. This property is also used occasionally. 
\item It is now possible to define K-theory groups
  $K_c^0(X;\boundedops)$ of finitely
  generated projective $\mathcal{B}$-modules in the same way as we
  define  such K-groups for $A$-module bundles in Section
  \ref{sec:k-theory-with}. We can also consider the normed $*$-algebra
  $C_0(X;\mathcal{B})$ of continuous $\mathcal{B}$-valued functions on
  a locally compact space $X$ which vanish at infinity., and the same
  proofs as for $A$-bundles implies that we get a commutative diagram
  \begin{equation}\label{eq:K-diagram}
    \begin{CD}
      K_c^0(X;\mathcal{B}) @>{\iso}>> K_0(C_0(X;\mathcal{B}))\\
      @V{\iso}VV @V{\iso}VV\\
      K_c^0(X;A) @>{\iso}>> K_0(C_0(X;A))
\end{CD}
\end{equation}
for any locally compact space $X$. Moreover, $C_0(X;\mathcal{B})$ is a
$*$-subalgebra of $C_0(X;A)$ which is closed under holomorphic
functional calculus. This is true since $f(\phi)$ for a holomorphic
function $f$ and $\phi\colon X\to A$ is given by
$(f(\phi))(x)=f(\phi(x))$ whenever it is defined, so that the
statement reduces to the fact that $\mathcal{B}$ is closed under
holomorphic functional calculus in $A$. This implies that the
vertical maps in \eqref{eq:K-diagram} are isomorphisms.
\item Note that we do not use $C(X)\tensor \mathcal{B}$ here, which
  will in general not be isomorphic to $C(X;\mathcal{B})$.
\item We can then also define $K_c^1(X;\mathcal{B}):=
  K^0_c(X\times\reals;\mathcal{B})$, such that
  $K^1_c(X;\mathcal{B})\xrightarrow{\iso} K^1_c(X;A)$ is an immediate
  consequence of the corresponding result for $K^0$. Moreover, we get
  a commutative Bott periodicity diagram
  \begin{equation*}
    \begin{CD}
      K^0(X;\mathcal{B}) @>{\beta}>{\iso}>
      K^0_c(X\times\reals^2;\mathcal{B})\\
      @VV{\iso}V @VV{\iso}V\\
      K^0(X;A) @>{\beta}>{\iso}> K^0_c(X\times\reals^2;A),
    \end{CD}
  \end{equation*}
  where we might define the upper horizontal map such that the diagram
  commutes.
\item Given a trace $\tau\colon\mathcal{B}\to Z$ on $\mathcal{B}$, the
  constructions of Section \ref{sec:traces-hilbert-a} immediately
  generalize. In particular, for each finitely generated projective
  Hilbert $\mathcal{B}$-module bundle $\mathcal{W}$ on a manifold $M$ we get an induced
  map
  \begin{equation*}
    \tau\colon\Omega^*(M;\End_{\mathcal{B}}(\mathcal{W})) \to \Omega^*(M;Z).
  \end{equation*}
\item Connections on smooth Hilbert $\mathcal{B}$-module bundles are
  defined in exactly the same way as they are defined on $A$-bundles,
  and the properties proved in Section \ref{sec:connections}
  generalize immediately to Hilbert $\mathcal{B}$-module bundles. Moreover,
  each $\mathcal{B}$-module connection induces a connection on the
  induced Hilbert $A$-module bundle which restricts to the given one
  on the $\mathcal{B}$-subbundle. 
\item If $\nabla$ is a connection on a finitely generated projective Hilbert $\mathcal{B}$-module
  bundle over $M$ and $\Omega$ is its curvature, and if $\tau\colon
  \mathcal{B}\to Z$ is a trace, then the Chern character
  \begin{equation*}
    \ch_\tau(\Omega)\in \Omega^{2*}(M;Z)
  \end{equation*}
  with its corresponding de Rham cohomology class are defined as in
  Section \ref{sec:chern-weyl-theory}. Moreover, all the properties
  proved there generalize to this situation.
\item\label{item:twistsituation} Let $M$ is a closed manifold, $D$ a Dirac type operator on $M$
  and $\mathcal{W}$ a smooth finitely generated projective Hilbert
  $\mathcal{B}$-module bundle with connection and with curvature
  $\Omega$ on $M$ with induces Hilbert $A$-module
  bundle $W$. Let $\tau\colon \mathcal{B}\to Z$ be a trace with
  induced homomorphism $\tau\colon K_0(\mathcal{B})\iso K_0(A)\to
  Z$. Commutativity in the K{\"u}nneth diagram (compare \eqref{eq:kuenneth})
  \begin{equation*}
    \begin{CD}
      K^*(\mathcal{B})\tensor K_*(M)\tensor \rationals @>>>
      K_0(M;\mathcal{B})\tensor \rationals\\
      @VV{\iso}V @VV{\iso}V\\
      K^*(A)\tensor K_*(M)\tensor\rationals @>{\iso}>> K_0(M;A)
  \end{CD}
\end{equation*}
implies that the upper horizontal map is an isomorphism. Recall that
the map sends $[\mathcal{P}]\tensor [E]$, where $\mathcal{P}$ is a
finitely generated projective $\mathcal{B}$-module and $E$ a finite
dimensional vector bundle to the fiberwise tensor product bundle
$[\mathcal{P}\tensor E]$, and
$[\mathcal{W},\mathcal{W_2},\phi_{\mathcal{W}},\phi_{\mathcal{W}_2}]\tensor
[E,E_2,\phi_E,\phi_{E_2}]$
to 
\begin{multline*}
\beta^{-1}[\mathcal{W}\tensor E\oplus
\mathcal{W_2}\tensor E_2, \mathcal{W}\tensor E_2\oplus
\mathcal{W_2}\tensor E,\\
 \phi_{\mathcal{W}}\tensor\phi_E\oplus
\phi_{\mathcal{W}_2}\tensor \phi_{E_2}, \phi_{\mathcal{W}}\tensor
\phi_{E_2}\oplus \phi_{\mathcal{W}_2}\tensor \phi_E],
\end{multline*}
where the
objects are now corresponding tuples of bundles on $\reals$,
$M\times\reals$ and $M\times \reals^2$ as in the description of
compactly supported K-theory of Proposition
\ref{prop:K_theory_of_non_compacts_via_bundles}. 
\end{enumerate}
Given this, the proof of Theorem \ref{theo:index_theorem} now
  goes through for Hilbert $\mathcal{B}$-module bundles as it does for
  $A$-bundles, and we arrive, in the situation and with the notation
  of \ref{item:twistsituation}, at the formula we want to prove.
  \begin{theorem}
    \begin{equation*}
    \ind_\tau(D_W) = \innerprod{\ch(\sigma(D))\cup \Td(T_\complexs
      M)\cup \ch_\tau(\mathcal{W}), [TM]},
  \end{equation*}
\end{theorem}
  and if $M$ is oriented of dimension $n$, we get similarly
    \begin{equation*}
    \ind_\tau(D_W) = (-1)^{n(n-1)/2}\innerprod{\pi_!\ch(\sigma(D))\cup \Td(T_\complexs
      M)\cup \ch_\tau(\mathcal{W}), [M]},
  \end{equation*}
  compare also Theorem \ref{theo:index_theorem}. Note that
  $\ch_\tau(\mathcal{W})=\ch_\tau(\Omega)$ is only defined in terms of
  the curvature of the $\mathcal{B}$-bundle, because the trace is only
  defined on $\mathcal{B}$. Only after passage to K-theory can it also
  be used for $K_0(A)$.

The extensions of the index theory to $\mathcal{B}$-module bundles
were inspired by conversations with John Roe and his proof of the
relative index theorem in \cite{MR92g:58125}, which also uses traces
on densely defined subalgebras.

{\small
  \bibliographystyle{amsxport}
\bibliography{L2index}
}

\end{document}